\documentclass[12pt]{article}

\usepackage{amsmath,amssymb,mathabx}
\usepackage{mathrsfs,upgreek,calligra}
\usepackage{amsthm}
\usepackage{dsfont}
\usepackage[utf8]{inputenc}
\usepackage[T1]{fontenc}
\usepackage{lmodern} 
\usepackage{ucs}
\usepackage{fullpage}

\usepackage{graphicx}
\usepackage{caption}
\usepackage{subcaption}

\usepackage{array}
\usepackage{multicol}
\usepackage{multirow}
\usepackage{color}
\usepackage{authblk}
\usepackage{esint}
\usepackage{enumitem,hyperref}

\usepackage{accents}

\usepackage{tikz}
\usetikzlibrary{decorations}
\usetikzlibrary{decorations.pathreplacing}
\usetikzlibrary{arrows}

\newtheorem{lemma}{Lemma}[section]
\newtheorem{theo}[lemma]{Theorem}
\newtheorem{rmk}[lemma]{Remark}
\newtheorem{proposition}[lemma]{Proposition}

\DeclareMathAlphabet{\mathpzc}{OT1}{pzc}{m}{it}
\DeclareMathAlphabet{\mathcalligra}{T1}{calligra}{m}{n}

\makeatletter
\def\namedlabel#1#2{\begingroup
    #2%
    \def\@currentlabel{#2}%
    \phantomsection\label{#1}\endgroup
}
\makeatother

\makeatletter
\renewcommand*{\eqref}[1]{%
  \hyperref[{#1}]{\textup{\tagform@{\ref*{#1}}}}%
}
\makeatother

\newcommand{\QED}{\mbox{}\hfill \raisebox{-0.2pt}{\rule{5.6pt}{6pt}\rule{0pt}{0pt}} \medskip\par}
\newcommand{\N}{\mathbb{N}}

\newcommand{\R}{\mathbb{R}}
\newcommand{\C}{\mathbb{C}}

\newcommand{\Ker}{\mathrm{Ker}}
\newcommand{\ds}{\displaystyle}
\newcommand{\ud}{\, {\mathrm{d}}}

\title{On the dynamics of $N$  classical particles interacting with a complex environment}

\author[1]{Thierry~Goudon\thanks{ {\tt thierry.goudon@univ-cotedazur.fr}}}
\author[1]{Bruno Marcos\thanks{ {\tt thierry.goudon@univ-cotedazur.fr}}}
\author[1]{Thomas Oudeville\thanks{ {\tt thomas.oudeville@univ-cotedazur.fr}}}
\author[1,2]{\\ Simona Rota Nodari\thanks{ {\tt simona.rotanodari@univ-cotedazur.fr}}}
\date{}

\affil[1]{\small Universit\'e C\^ote d'Azur, CNRS, LJAD,  

Parc Valrose, F-06108 Nice, France}

\affil[2]{\small Institut Universitaire de France}

\begin{document}
\maketitle

\begin{abstract}
The interactions of a particle with its environment eventually lead to the damping of its motion.
This phenomenom is well understood when considering a single particle exchanging 
energy with the environment.
The situation is definitely more complex when many particles are considered: 
the energy dissipation mechanisms are not clear whenever particles create a common well in which they keep interacting.

In this framework, we identify the equilibrium solutions, that consist in clusters of particles at rest, and we establish their spectral stability.
We show that the mean position/velocity pair is damped  exponentially fast at the linearized level. 

Coming back to the non linear problem, we can show, 
at least in some simplified cases, that the dissipation mechanisms also involve far slower 
rates, that explain the difficulty in observing 
the relaxation behavior on the individual trajectories.

\end{abstract}

\vspace*{.5cm}
{\small
\noindent{\bf Keywords.}
Open  systems. Particles interacting with a vibrational field.}  
\\[.4cm]

\noindent{\bf Math.~Subject Classification.} 
82C22 
82C70 
37K05 

\section{Introduction}

It is a common experience to realize that the speed of a particle decays exponentially fast: the decay is attributed to a force exerted by the  environment of the particle, which is proportional to the velocity and opposed to the motion. This is a phenomenological 
description of the transfer 
of energy between the particle and the environment.  A more accurate treatment should take into account the evolution of both the particle and its environment, from which the dissipation mechanism should emerge.
 
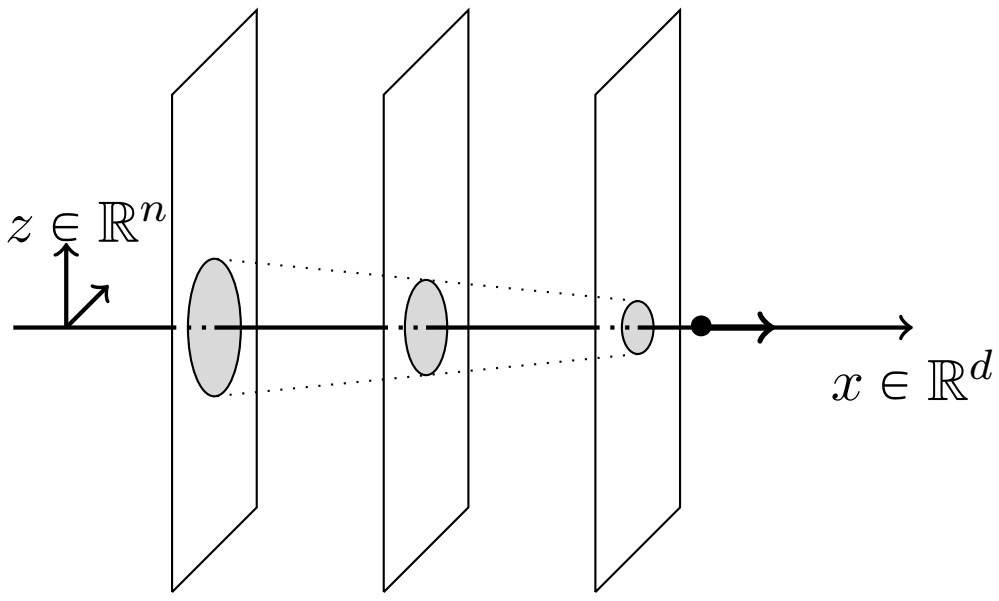
\begin{figure}[!htpb]
\centering
\begin{tikzpicture}[scale=.7]
\draw (8,0) node[below] {{\small{$x\in\mathbb R^d$}}};
\draw (0.2,1) node {{\small{$z\in\mathbb R^n$}}};
\draw [thick,->] (0,0) -- (0,0.8);
\draw [thick,->] (0,0) -- (0.4,0.4);
\draw [thick] (-0.5,0) -- (1,0) ;
 \fill[color=gray!30](5.4,0) ellipse (.15cm and .25cm);  
 \draw(5.4,0) ellipse (.15cm and .25cm);  
\fill[color=gray!30] (3.4,0) ellipse (.20cm and .45cm);
\draw (3.4,0) ellipse (.20cm and .45cm);
\draw [thick,dotted] (3.0,0) -- (3.4,0);
\draw [thick]  (3.4,0) -- (5,0) ;
\fill[color=gray!30] (1.4,0) ellipse (.25cm and .65cm);
\draw (1.4,0) ellipse (.25cm and .65cm);
\draw [thick,dotted] (5.0,0) -- (5.4,0);
\draw [thick,->]  (5.4,0) -- (8,0) ;
\draw [thin,dotted] (5.4,.25) -- (1.4, .65);
\draw [thin,dotted] (5.4,-.25) -- (1.4, -.65);
\draw [thick,dotted] (1.0,0) -- (1.4,0);
\draw [thick]  (1.4,0) -- (3,0) ;
\draw [thin] (1,-2.5) -- (1,2.2) -- (1.8,3) -- (1.8,-1.7) -- (1,-2.5);
\draw [thin] (3,-2.5) -- (3,2.2) -- (3.8,3) -- (3.8,-1.7) -- (3,-2.5);
\draw [thin] (5,-2.5) -- (5,2.2) -- (5.8,3) -- (5.8,-1.7) -- (5,-2.5);
\draw (6,0) node[] {\small{$\bullet$}};
\draw [very thick,->]  (6,0)--(6.7,0) ;
\end{tikzpicture}
\caption{A schematic picture of the particle-environment interactions} 
\label{membranes}
\end{figure}

The starting point of this work is the following coupled system, introduced by L.~Bruneau \& S.~de Bi\`evre \cite{BdB}: 
here,  $t\mapsto q(t)$ stands for  the position of a classical  particle,  while the environment is described by a vibrating field $(t,x,z)\mapsto \psi(t,x,z)$
and we have
\begin{subequations}
    \begin{alignat}{1}
        \label{particle}&\ds \ddot{q}(t)=
       -\iint_{\R^{d}\times\R^{n}}\sigma_{1}(q(t)-y)\sigma_{2}(z)\nabla_{y}\psi(t,y,z)\ud y\ud z,
       \\
        \label{wave}&\ds\frac{1}{c^{2}}\partial_{tt}^{2}\psi-\Delta_{z}\psi=-\sigma_{2}(z)\sigma_{1}(x-q(t)),
    \end{alignat}
\end{subequations}
for $t\geq 0$, $x\in\R^{d}$, $z\in\R^{n}$.
System \eqref{particle}--\eqref{wave} is completed by initial data 
\begin{equation}\label{cond_init}
(q,\dot{q})\big|_{t=0}=(q_{0},p_{0}), \quad
 (\psi,\partial_{t}\psi)\big|_{t=0}=(\psi_{0},\psi_{1}).
\end{equation}
Figure~\ref{membranes}  guides the intuition on the derivation of the model. The particle evolves in the physical  space $
\mathbb R^d$, while at each space position $x\in \mathbb R^d$ is attached a membrane, oscillating on a \emph{transverse} direction $z\in \mathbb R^n$.
The motion space and the vibration space are distincts and there is no a priori relation between $n$ and $d$.
The environment is this  thought of as a (continuum) set of membranes,
 activated by the passage of the particle.
The details of the interaction are thus embodied into:
\begin{itemize}
\item two form functions $x\in \mathbb R^d\mapsto \sigma_1(x)$ and $z\in \mathbb R^n\mapsto \sigma_2(z)$
that determine the interaction  domain, in the physical and the transverse directions respectively,  between the particle and the waves; they are both non negative, spherically symmetric, infinitely smooth and compactly supported; 
\item the (uniform) wave speed $c>0$.
\end{itemize}
This model is reminiscient to the framework introduced in \cite{CL} or \cite{JP2}.

As time becomes large, the remarkable fact brought out in \cite{BdB} is that 
 the membranes eventually act as a friction force on the particle, as described by the following statement
 (see \cite[Theorem 2]{BdB} for further details, we also refer to \cite{KKS,KKSb} for the analysis of similar models).

\begin{theo}{\cite[Theorem 2]{BdB}}\label{friction}
Let $n=3$. There exists 
a constant $\upgamma>0$, determined by the parameters of the equations, such that 
for any $\eta>0$, we can find 
a critical wave speed $c_{*} =c_{*}(\eta)>0$ 
and a constant $ K>0$   
 such that, for any  $c\geq c_{*}$, we have
\[
    \left|q_{\infty}-q(t)\right|\leq Ke^{-\frac{\upgamma(1-\eta)}{c}t},
    \qquad 
 \left|\dot q(t)\right|\leq Ke^{-\frac{\upgamma(1-\eta)}{c}t}
\]
for a certain $q_{\infty}\in\R^{d}$.
\end{theo}

As revealed by this statement,  the transverse dimension $n$ 
and the wave speed $c$  should be large enough to provide enough dissipation, so that the energy is efficiently  evacuated at infinity by the vibrational field.
(We warn the reader that we adopt here a scaling of the problem that slightly differs from the one in \cite{BdB}; this choice has been motivated in \cite{dBGV,Vi3}.) 

It is natural to extend the model \eqref{particle}--\eqref{wave} by considering a set of $P$ particles which all interact with the membranes \cite{GV}.
Let $q_{i}$ stand for the  position of the  $i$th particle. The system is now governed by  
\begin{subequations}
    \begin{alignat}{11}
        \label{N_particles}&\ds \ddot{q}_{i}(t)=-\iint_{\R^{d}\times\R^{n}}\sigma_{1}(q_{i}(t)-y)\sigma_{2}(z)\nabla_{y}\psi(t,y,z)\ud y\ud z,\\
        \label{N_wave}&\Big(\ds\frac1{c^{2}}\partial_{tt}^{2}\psi-\Delta_{z}\psi\Big)(t,x,z)=-\sigma_{2}(z)\left(\sum\limits_{j=1}^{P}\sigma_{1}(x-q_{}j(t))\right),    \end{alignat}
\end{subequations}
for $t\geq 0$, $x\in\R^{d}$, $z\in\R^{n}$, $i\in\{1,...,P\}$.
The mean field regime leads to a Vlasov-like equation for the 
particle distribution function \cite{GV}, and the interaction eventually lead to Landau damping phenomena \cite{Vi2}.
The behavior of the model \eqref{particle}--\eqref{wave}  for a single particle, and  
\eqref{N_particles}--\eqref{N_wave}  for $P$ particles has been investigated on numerical grounds on \cite{Vi1}. Theorem~\ref{friction} is numerically verified 
but the case $P>1$ is definitely less clear: several scenario emerge  depending 
on the initial state and situations where the damping is far from neat have been observed, with particles bouncing in a common potential well.

We wish to analyse in further details the situation when there are $P>1$ interacting particles.
As it will be explained below,  and already pointed out in \cite{BdB}, we need the  dimension  of the membranes 
to be large enough $n\ge3$, and in fact the case $n=3$ is specific for a part of our results.
Assuming that $n$ 
is odd and that the initial data $(\psi,\partial_t\psi)\big|_{t=0}$ are compactly supported simplifies the analysis. In this case, by Huygens' principle, the action of these initial data on the particles vanishes for sufficiently large times. In the context of statistical physics, this situation corresponds to zero temperature. At positive temperature, the initial field would instead be random, distributed according to the Gibbs measure of the free wave equation, and the free field would then act on the particles as a stationary random force which never vanishes. In that case one expects a Brownian-type motion of the particles, and relaxation in distribution towards thermal equilibrium, rather than convergence of each trajectory. We will go back to this issue elsewhere.

To perform the analysis, we combine two complimentary 
formulations of the problem, discussed in Section~\ref{prelim}.
The former relies on the underlying energetic and Hamiltonian structures of the system  that couple the particles and the environment.
The latter can be seen as a ``non-autonomous'' formulation of the problem,  
obtained by plugging into the definition of the potential the solution of the wave equation sourced by the presence of particles. 
This expression makes memory effects appear, that produce the friction effects.
In Section~\ref{En_Min}, we identify equilibrium solutions, which can be interpreted as (local) minimizers of the energy: they 
 appear as clusters of particles at rest (Proposition~\ref{prop_energymin}).
The spectral analysis of the linearized problem brings out typical effects of the many particle
framework, with the existence of oscillating modes (Proposition~\ref{prop_spectrum_lin}).
We further study the linearized equations in Section~\ref{Tlarge_Lin}, showing that the mean position/velocity 
are damped, with exponential rate, in line with numerical observations (Proposition~\ref{prop:meanLin}).
Finally, we pay attention in Section~\ref{TwoPart} to the specific case of two particles in dimension $d=1$, where the initial data are prepared so that 
the solution keeps zero mean velocity. We are then able to show that the relative position/velocity
also tend to 0 as time runs to infinity, with a rate of the order of $\mathscr O\left(\frac{1}{\sqrt {1+t/c}}\right)$ (Theorem~\ref{th_wp}).
This result exhibits the coexistence of different time scales at which 
the relaxation phenomena occurs when there are many interacting particles.

\section{Preliminary observations and reformulation of the problem}
\label{prelim}

Throughout the paper, the form functions $\sigma_1,\sigma_2$ are required to satisfy the following requirements:
\begin{itemize}
\item [h1)] $\sigma_1:\mathbb R^d\to [0,\infty)$ is $C^\infty$ smooth, radially symmetric and compactly supported in $B(0,R_1)$. With a slight abuse of notation, 
denoting $\sigma_1(x)=\sigma_1(|x|)$,
we assume that $r\geq 0\mapsto \sigma_1(r)$ is non increasing. 
\item [h2)] $\sigma_2:\mathbb R^n\to [0,\infty)$ is $C^\infty$ smooth, radially symmetric and compactly supported in $B(0,R_2)$.
\end{itemize}

In what follows we will deal with 
\[\Sigma_1(x)=\sigma_1*\sigma_1(x)=\ds\int_{\mathbb R^d} \sigma_1(x-y)\sigma_1(y)\ud y\]
which is also radially symmetric, with support in $B(0,2R_1)$.
With the same abuse of notation, the function  $r\mapsto \Sigma_1(r)$ can be shown to be  decreasing on its support, thus reaching its maximum at $x=0$.
A closely related statement appeared in \cite{Lieb}. 

A fundamental feature of the model is its Hamiltonian structure. In particular, we can set $X=(q_1,...,q_P,p_1,...,p_P,\psi,\varpi)$
and  rewrite the problem as 
\begin{equation}\label{aut0}
\partial_t X=\mathscr J\nabla_X \mathscr H(X).\end{equation}
Here, $\mathscr J$ is a skew-symmetric operator which defines the symplectic structure while $\mathscr H$ is an energy functional which is a conserved quantity thanks to Noether's Theorem.  
A natural choice for $\mathscr J$ is given by 
\begin{equation}\label{aut1}
\mathscr J=\begin{pmatrix}
0 & \mathbb I & 0 & 0
\\
-\mathbb I & 0& 0 & 0
\\
0 & 0 & 0 & c
\\
0 & 0& -c & 0
\end{pmatrix}
.\end{equation}
In this  case the corresponding energy functional $\mathscr H$ is 
\begin{equation}\label{aut2}\begin{array}{lll}
\mathscr H(X)&=&\ds\sum_{j=1}^P \ds\frac{|p_j|^2}{2}+
\ds\frac{1}{2}\ds\iint_{\mathbb R^d\times \mathbb R^n} |\varpi|^2\ud z\ud x+ \ds\frac12\ds\iint_{\mathbb R^d\times\mathbb R^n} |\nabla_z \psi|^2\ud z\ud x
\\&&+\ds\sum_{j=1}^P\ds\iint_{\mathbb R^d\times\mathbb R^n} 
\sigma_1(x-q_j)\sigma_2(z)\psi(x,z)
\ud z\ud x\end{array}\end{equation}
which is conserved by the dynamics. Moreover, we can easily see that the system~\eqref{N_particles}-\eqref{N_wave} is invariant under translations of the form 
\begin{equation*}
  \begin{pmatrix}
    q(t)\\
    p(t)\\
    \psi(t,x,z)\\
    \varpi(t,x,z)\\
  \end{pmatrix}
  \to
  \begin{pmatrix}
    q(t)+h\\
    p(t)\\
    \psi(t,x+h,z)\\
    \varpi(t,x+h,z)\\
  \end{pmatrix}.
\end{equation*}
This leads to the conservation of the total momentum
\begin{equation}\label{defmt}
 \mathscr P(X)=\ds\sum_{j=1}^P {p}_j -\frac{1}{c}\ds\iint_{\R^{d}\times\R^{n}}\varpi\nabla_x \psi\ud z\ud x.
\end{equation}
These conservation properties will be important in 
the analysis of the equilibria.

A complementary viewpoint focuses on the ``resolution'' of the wave equation.
To this end, we proceed by using the fact that  the Fourier transform\footnote{Throughout the paper we use the convention $\widehat f(\xi)=\int_{\mathbb R^n} f(x)e^{-ix\cdot \xi}\ud x$} of $\psi$ is given by 
\begin{equation}
\label{Fourrier}\begin{array}{lll}
\widehat{\psi}(t,x,\xi) &=& 
\ds \widehat{\psi_0}(x,\xi) \cos(c|\xi|t) + \widehat{\psi_1}(x,\xi) \frac{\sin(c|\xi|t)}{c|\xi|}
\\[.4cm]&&\qquad\qquad
-\ds\sum_{j=1}^P \ds\int_0^t  \sigma_1 (x-q_j(t-s))\  \widehat{\sigma_2}(\xi) \frac{c\sin(c s|\xi|)}{|\xi|} \ud s
. \end{array}\end{equation}
This corresponds to the splitting $\psi=\psi_I+\psi_S$, where $\psi_I$ satisfies the homogeneous Cauchy problem 
$\frac{1}{c^2}\partial^2_{tt}\psi_I-\Delta_z \psi_I=0$ with $(\psi_I,\partial_t\psi_I)\big|_{t=0}=(\psi_0,\psi_1)$ and
$\psi_S$ is  the solution of the wave equation with source term 
$\frac{1}{c^2}\partial^2_{tt}\psi_S-\Delta_z \psi_S=-\sigma_2(z)\sum_{j=1}^P \sigma_1(x-q_j(t))$ 
and as initial data $(\psi_S,\partial_t\psi_S)\big|_{t=0}=0$.

Let us denote by 
\begin{equation}
  \label{def_potentialPhi}
  \Phi(t,x)= \iint_{\R^{d}\times\R^{n}}\sigma_{1}(x-y)\sigma_{2}(z)\psi(t,y,z)\ud y\ud z,
\end{equation}
 the  potential describing the interaction with the environment. Note that 
\begin{equation*}
  \Phi(t,x)= \frac{1}{(2\pi)^n}\iint_{\R^{d}\times\R^{n}}\sigma_{1}(x-y)\widehat\sigma_{2}(\xi)\widehat\psi(t,y,\xi)\ud y\ud \xi,
\end{equation*} 
and equation~\eqref{N_particles} reads as
\begin{equation*}
  \dot q_j(t)=p_j(t),\qquad
 \ds\dot p_j(t)=-\nabla\Phi(t,q_j(t)).
\end{equation*}
As a consequence, using the splitting $\psi=\psi_I+\psi_S$, we can decompose the potential as 
\[\Phi=\Phi_I+\Phi_S\]
where the former term only depends on the initial data 
\[\Phi_I(t,x)=
 \ds\frac{1}{(2\pi)^n}
\ds\int_{\mathbb R^d}\sigma_1 (x-y)\left(\ds \int_{\mathbb R^n} \left(\widehat{\psi_0}(y,\xi) \cos(c|\xi|t) + \widehat{\psi_1}(y,\xi) \ds\frac{\sin(c|\xi|t)}{c|\xi|} \right)\ 
\widehat{\sigma_2}(\xi) \ud\xi\right)\ud y\]
and the latter corresponds to a self-consistent potential 
\[ \Phi_S(t,x)   
=- \ds\sum_{j=1}^P
\ds\int_0^t \left(
\ds\int_{\mathbb R^d} \sigma_1(x-y)\sigma_1(y-q_j(t-s))\ud y
\int_{\mathbb R^n} \frac{c\sin(c s|\xi|)}{|\xi|} \frac{\left| \widehat{\sigma_2}(\xi) \right|^2}{(2\pi)^n} \ud \xi \right)\ud s 
.\]
Going back to the 
 equation of motion for the particles, 
 we write
 \begin{equation}\label{non_aut}
 \ds\frac{\ud}{\ud t} q_j(t)=p_j(t),\qquad
 \ds\frac{\ud}{\ud t} p_j(t)=-\nabla\Phi(t,q_j(t))=-\nabla\Phi_I(t,q_j(t))-\nabla\Phi_S(t,q_j(t)).\end{equation}
This formula makes the following kernel appear
 \[\mathpzc{p}_c(t)=
 \int_{\mathbb R^n} \frac{c\sin(c t|\xi|)}{|\xi|} \frac{\left| \widehat{\sigma_2}(\xi) \right|^2}{(2\pi)^n} \ud \xi
 \]
 which contains many information about the dissipation properties of the system.
 Key properties of the kernel are summarized in the following statement (see \cite[Lemma~14]{dBGV},
 \cite[Section~2.4]{Vi2}).
 
 \begin{lemma}\label{lemAV}
 Let $n\geq 3$ and let $\sigma_2$ satisfy (h2).
 Then, $ \mathpzc{p}_c\in L^1([0,\infty))$, with 
 $\int_0^\infty  \mathpzc{p}_c(t)\ud t=\kappa=\frac{1}{(2\pi)^n}\int_{\mathbb R^n}\frac{|\widehat \sigma_2(\xi)|^2}{|\xi|^2}\ud \xi\in (0,\infty)$.
  If, moreover, $n$ is odd, then $t\mapsto  \mathpzc{p}_c(t)$ is compactly supported in $[0,2R_2/c]$, and
 it lies in $L^\infty([0,\infty))$.
 \end{lemma}

 For further purposes, we will also need to evaluate the integral $\int_0^\infty t\mathpzc{p}_c(t)\ud t$.
 As already observed in \cite{BdB}, the value of this integral makes the case $n=3$ specific
 (and this statement identifies the friction coefficient given in Theorem~\ref{friction}).

 \begin{lemma}\label{mtpc}
 If $n>3$, then $\int_0^\infty t\mathpzc{p}_c(t)\ud t=0$, while for $n=3$, we get
\begin{equation}
  \label{def_gammaconstant}
   \ds\int_0^\infty t\mathpzc{p}_c(t)\ud t=\ds\frac
1{4 \pi c}\left( \ds\int_{\mathbb R^3} \sigma_2(z)\ud z\right)^2:=\ds\frac{\upgamma}{c}.
\end{equation}
 \end{lemma}

\noindent
{\bf Proof.}
Let $n\ge 3$. Since $\sigma_2$ is assumed to be radially symmetric, its Fourier transform is radially symmetric too and we denote $\widehat\sigma_2(\xi)=\tilde\sigma_2(|\xi|)$. 
Using spherical coordinates, we get
\[
\mathpzc{p}_c(t)=\ds\frac{c|\mathbb S^{n-1}| }{(2\pi)^n }
\ds\int_0^\infty \sin(ctr) |\tilde\sigma_2(r)|^2 r^{n-2}\ud r
=c\mathpzc{p}_1(ct).\]

Now, consider the function $f: r\geq 0 \mapsto |\tilde\sigma_2(r)|^2 r^{n-2}\in \R$. Since it vanishes for $r=0$, we can extend it 
as an odd function defined on $\mathbb R$ and we have 
\begin{equation*}
  \widehat f(\tau)=\int_{\R}f(r)e^{-i\tau r}\ud r=-2i\int_0^{\infty}  f(r)\sin(\tau r)\ud r.
\end{equation*}
As a consequence,
\[
\mathpzc{p}_c(t)=\ds\frac{c|\mathbb S^{n-1}| }{(2\pi)^n }\ds\int_{\mathbb R} \sin(ctr) f(r)\ud r
=
-\ds\frac{c|\mathbb S^{n-1}| \widehat f(ct)}{2i(2\pi)^n }=\ds i\frac{c|\mathbb S^{n-1}| \widehat f(ct)}{2(2\pi)^n }.\]
Note that the obtained expression is  an odd function of the variable $t$.
Then, we are led to compute
\[\begin{array}{lll}
\ds\int_0^\infty t\mathpzc{p}_c(t)\ud t
&=&\ds i\frac{|\mathbb S^{n-1}| }{2(2\pi)^n c}\ds\int_0^\infty  ct\widehat f(ct)\, c\ud t
=\ds\frac{|\mathbb S^{n-1}| }{4(2\pi)^n c}\ds\int_{\mathbb R}  i\tau\widehat f(\tau)\ud \tau
\\[.4cm]
&=&\ds\frac{|\mathbb S^{n-1}| }{4(2\pi)^n c}\ds\int_{\mathbb R}
\widehat {f'}(\tau)\ud\tau 
=\ds\frac{|\mathbb S^{n-1}| }{4(2\pi)^n c}{2\pi}
f'(0).
\end{array}\]
It remains to evaluate 
$f'(r)= |\tilde\sigma_2(r)|^2 (n-2)r^{n-3}+ 2r^{n-2}\tilde\sigma_2(r)\tilde\sigma'_2(r)$ for $r=0$. This
leads to $f'(0)=0$ if $n>3$ and $f'(0)=|\tilde\sigma_2(0)|^2=(\widehat \sigma_2(0))^2$ if $n=3$.
As a conclusion, $\int_0^\infty t\mathpzc{p}_c(t)\ud t=0$ for $n>3$. For $n=3$, we obtain
\begin{equation*}
  \int_0^\infty t\mathpzc{p}_c(t)\ud t= \ds\frac{|\mathbb S^{2}| }{4(2\pi)^3 c}{2\pi}
  |\widehat \sigma_2(0)|^2=\ds\frac{1}{4\pi c}
  \left(\int_{\R^3} \sigma_2(x)\ud x\right)^2.
\end{equation*}
\QED

In some sense, we have at hand a non autonomous formulation
\eqref{non_aut}
and an autonomous formulation \eqref{aut0}-\eqref{aut2}
of the problem, and going back and forth between the two will be fruitful 
in this study.

\section{Energy minimizers, equilibrium solutions, linearization, spectral stability}
\label{En_Min}

\subsection{Equilibrium solutions}

An important family of solutions to the problem~\eqref{N_particles}-\eqref{N_wave} are the equilibrium solutions $X^*$ which are such that $\partial_t X^*(t)=0$. This implies
\begin{equation*}
  p_j^*=0 \qquad \varpi^*=0
\end{equation*}
and $(q_j^*,\psi^*)$ is solution to the system 
\begin{equation*}
  \left\{
    \begin{aligned}
      &\iint_{\R^{d}\times\R^{n}}\sigma_{1}(q_{j}^*-y)\sigma_{2}(z)\nabla_{y}\psi^*(y,z)\ud y\ud z=0, \\
      &-\Delta_z\psi^*(x,z)=-\sigma_2(z)\left(\sum\limits_{k=1}^{P}\sigma_{1}(x-q^*_{k})\right).
    \end{aligned}
  \right.
\end{equation*}
Now, let us introduce the solution $z\mapsto \Gamma(z)$ of 
\[-\Delta_z \Gamma=\sigma_2\]
and  set 
\[\kappa=\ds\int_{\mathbb R^n}\Gamma(z)\sigma_2(z)\ud z=\ds\int_{\mathbb R^n}\ds\frac{|\widehat\sigma_2(\xi)|^2}{|\xi|^2}\ds\frac{\ud \xi}{(2\pi)^n}>0.\]
(Note that the definition is meaningful provided $n\geq 3$.)
We get
$$\psi^*(x,z)=-\Gamma(z)\sum_{k=1}^P \sigma_1(x-q^*_k).$$
Finally, the force term should vanish, for any $j\in\{1,...,P\}$. This leads to
\begin{align}\label{eq_defequilibrium}
  0&=\iint_{\R^{d}\times\R^{n}}\sigma_{1}(q_{j}^*-y)\sigma_{2}(z)\nabla_{y}\psi^*(y,z)\ud y\ud z\nonumber\\
  &=-\iint_{\R^{d}\times\R^{n}}\sigma_{1}(q_{j}^*-y)\sigma_{2}(z)\Gamma(z)\sum_{k=1}^P \nabla\sigma_1(y-q^*_k)\ud y\ud z\nonumber\\
  &=-\kappa\sum_{k=1}^P \int_{\R^{d}}\sigma_{1}(q_{j}^*-y)\nabla\sigma_1(y-q^*_k)\ud y\ud z=-\kappa\sum_{k=1}^P (\sigma_{1}\star \nabla\sigma_1)(q^*_j-q^*_k)\nonumber\\
  &=-\kappa\sum_{k=1}^P \nabla(\sigma_{1}\star \sigma_1)(q^*_j-q^*_k)=-\kappa\sum_{k=1}^P \nabla\Sigma_1(q^*_j-q^*_k)
\end{align}
with $\Sigma_1=\sigma_1\star\sigma_1$.

Intuitively, we can then deduce that equilibrium solutions are such that
\begin{itemize}\item either the particles are far away from each other and do not interact: $|q^*_j-q^*_k|\ge 2R_1$, with $\mathrm{supp}(\sigma_1)=B(0,R_1)$,
\item or the particles occupy the same position.
\end{itemize}
We actually claim that this is the only possible equilibrium configuration. In other words, we have the following proposition.

\begin{proposition}\label{prop_eq_solutions}
Suppose $n\geq 3$.
  Let $q^*=(q_1^*,\ldots, q_P^*)\in (\R^d)^P$ a configuration that satisfies~\eqref{eq_defequilibrium}. Fix $j\in\{1,\dots,P\}$. For any $k\in\{1,\dots,P\}$, either $|q^*_j-q^*_k|\ge 2R_1$, or $q_j^*=q_k^*$.
\end{proposition}

\noindent
{\bf Proof.}  First of all, we recall that the function $\Sigma_1$ is spherically symmetric 
  with $\mathrm{supp}(\Sigma_1)\subset B(0,2R_1)$. Moreover, $\nabla \Sigma_1(x)\neq 0$ for any $x\in B(0,2R_1)\setminus \{0\}$. In particular, for any $x\neq 0$, 
  $$\nabla\Sigma_1(x)=\frac{\Sigma_1'(|x|)}{|x|}x.$$

  Now, let $q^*=(q_1^*,\ldots, q_P^*)\in (\R^d)^P$ a configuration that satisfies~\eqref{eq_defequilibrium}. For any $j\in\{1,\dots,P\}$, we denote by $\Pi_j=\{k\in \{1,\dots,P\},\ q^*_j=q_k^*\}$, $m_j=\#\Pi_j$ and, for any $k\in \Pi_j$, we identify $q^*_k$ with the pseudo-particle $\tilde q_j^*$. If $m_j=P$, we are done. Otherwise, we obtain a configuration of (pseudo-)particles $\tilde q^*=(\tilde q_1^*,\ldots, \tilde q_{P'}^*)$ with $1<P'\le P$ and such that $\tilde q^*_k\neq \tilde q_\ell^*$ for any $k\neq \ell$. 
  This configuration satisfies the relation~\eqref{eq_defequilibrium} that can be written as 
  \begin{equation}\label{eq_defequilibrium_uniq}
    0=-\sum_{k=1}^{P'}\sum_{\ell\in \Pi_k}\nabla\Sigma_1(\tilde q_j^*-q^*_\ell)=-\sum_{k=1}^{P'}m_k \nabla\Sigma_1(\tilde q_j^*-\tilde q^*_k)= -\sum_{k=1}^{P'}m_k \frac{\Sigma_1'(|\tilde q_j^*-\tilde q^*_k|)}{|\tilde q_j^*-\tilde q^*_k|}(\tilde q_j^*-\tilde q^*_k)
  \end{equation}
  for any $j\in \{1,\ldots,P'\}$. Let us denote 
  \begin{equation}\label{defdjk}
    d_{j,k}=-m_k \frac{\Sigma_1'(|\tilde q_j^*-\tilde q^*_k|)}{|\tilde q_j^*-\tilde q^*_k|}
  \end{equation}
  for any $j,k\in\{1,\ldots, P'\}$. Owing to h1), we remark that $d_{j,k}\ge 0$ for any $j,k\in\{1,\ldots, P'\}$ and $d_{j,k}= 0$ if and only if $|\tilde q_j^*-\tilde q^*_k|\ge 2R_1$. 
  Our goal is now to prove that for any fixed $j\in\{1,\dots,P'\}$, $|\tilde q^*_j-\tilde q^*_k|\ge 2R_1$ for any $k\in\{1,\dots,P'\}$, $k\neq j$. 
  
  We proceed as follows. Let $\mathcal G$ a graph with $P'$ vertices, labeled $\tilde q_1^*, \ldots, \tilde q_{P'}^*$, each vertex representing one of the (pseudo-)particle $\tilde q_j^*$. 
  Two vertices $\tilde q_k^*, \tilde q_\ell^*$ are connected by an edge if and only if $|\tilde q_k^*- \tilde q_\ell^*|<2R_1$. Fix $j\in\{1,\dots,P'\}$ and let $\mathcal G_j=(V_j, E_j)$ the connected component of $\mathcal G$ that contains the vertex $\tilde q^*_j$. If the order of $\mathcal G_j$ is $1$, we are done since this means that $|\tilde q^*_j-\tilde q^*_k|\ge 2R_1$ for any $k\in\{1,\dots,P'\}$, $k\neq j$. Assume, by contradiction, that the order of $\mathcal G_j$ is strictly bigger than $1$. After possibly relabelling the vertices of $\mathcal G_j$, we may assume that the set of vertices $V_j$ is given by $(\tilde q^*_1,\ldots,\tilde q^*_{P_j})$ with $1<P_j\le P'$ the order of $\mathcal G_j$. Since $P_j>1$, for any $\tilde q^*_\ell\in V_j$, there exists at least one $k_{\ell}\in\{1,\ldots,P_j\}$ such that $d_{\ell, k_\ell}>0$. As a consequence,~\eqref{eq_defequilibrium_uniq} can be rewritten as 
  \begin{equation}
    \label{eq_defmatriceDj}
    \tilde q^*_\ell=\sum_{k=1}^{P_j}\tilde d_{\ell,k}\tilde q_k^*,\text{ with } \tilde d_{\ell,k}=\frac{d_{\ell,k}}{\sum_{k=1}^{P_j}d_{\ell,k}}.
  \end{equation}
  Note that, if $\ell,k\in \{1,\ldots,P_j\}$ are such that $d_{\ell, k}>0$, then $d_{k,\ell}>0$ and the directed graph $\mathcal G_j$ is actually strongly connected.
  Finally relation~\eqref{eq_defmatriceDj} implies that, for any $i\in\{1,\ldots, d\}$, the vector $(\tilde q_1^*(i),\ldots,\tilde q_{P_j}^*(i))\in \R^{P_j}$ is a solution to 
  \[\begin{pmatrix} q_1\\\vdots \\q_{P_j}\end{pmatrix}=\underbrace{
\begin{pmatrix}
0 & \tilde d_{1,2}  & \dots & \tilde d_{1,P_j}
\\
\tilde d_{2,1} &0 &   \dots & \tilde d_{2,P_j}
\\
\vdots & \vdots & \ddots & \vdots
\\
\tilde  d_{P_j,1} & \tilde  d_{P_j,2} & \dots & 0
\end{pmatrix}}_{\tilde D_j(i)}\begin{pmatrix} q_1
\\
\vdots \\
q_{P_j}\end{pmatrix}.\]
  namely an eigenvector of the $P_j\times P_j$ matrix $\tilde D_j(i)$ associated to the eigenvalue $1$. Here, for any $k\in\{1,\ldots, P_j\}$, $\tilde q_k^*(i)$ denotes the $i$-th component of the vector $\tilde q_k^*\in \R^d$. The matrix $\tilde D_j(i)$ is a non-negative stochastic irreducible matrix. Indeed, for any $k,\ell\in \{1,\ldots,P_j\}$ $\tilde d_{\ell,k}\ge 0$, $\sum_{k=1}^{P_j}\tilde d_{\ell,k}=1$ and the matrix is irreducible since $\mathcal G_j$ is strongly connected. As a consequence, the Perron-Frobenius theorem implies that the eigenvalue $1$ is simple and $\Ker(\tilde D_j(i)-\mathbb I)$ is spanned by the vector $(1,\ldots,1)\in \R^{P_j}$. This gives $\tilde q_1^*(i)=\ldots=\tilde q_{P_j}^*(i)$ for any $i\in\{1,\ldots, d\}$ which contradicts the fact that $\tilde q_k^*\neq \tilde q_\ell^*$ for any $k,\ell\in \{1,\ldots,P'\}$.
\QED

As a consequence, the equilibrium solutions $X^*=(q_1^*,...q^*_P,p_1^*,...,p^*_P,\psi^*,\varpi^*)$ are defined by 
\begin{equation}\label{EqSol_bis}
  \begin{aligned}
    &p^*_1=...=p^*_P=0,\qquad \psi^*(x,z)=-\Gamma(z)\ds\sum_{j=1}^P \sigma_1(x-q^*_j), \qquad \varpi^*=0,\\
    &\left\{
      \begin{aligned}
        &\{1,...,P\}=\ds\bigcup_{\ell=1}^{L} \Pi_\ell, \textrm{ for some } L\in\{1,...,P\} \textrm{ with } \ds\sum_{\ell=1}^{L} \#\Pi_\ell=P,\\
        &\Pi_\ell\cap \Pi_{\ell'}=\emptyset \textrm{ if }\ell\neq \ell',\\
        &|q^*_k-q^*_j|{\ge}2R_1 \textrm{ if } k\in \Pi_\ell,\,j\in \Pi_{\ell'},\,\ell\neq \ell',\\
        &q^*_k=q^*_j \textrm{ if } k,j\in \Pi_\ell.
      \end{aligned}
    \right.
  \end{aligned}
\end{equation}

The associated energy reads
\[\begin{array}{lll}
\mathscr H(X^*)&=&\ds\frac12
\ds\iint_{\mathbb R^d\times\mathbb R^n}|\nabla_z\Gamma(z)|^2\Big|
\ds\sum_{j=1}^P \sigma_1(x-q_j^*)\Big|^2\ud z\ud x
\\
&&
\hspace*{3cm}-
\ds\iint_{\mathbb R^d\times\mathbb R^n}\Gamma(z)\sigma_2(z)
\ds\sum_{j=1}^P \sigma_1(x-q_j^*)\ds\sum_{k=1}^P \sigma_1(x-q_k^*)\ud z\ud x
\\
&=&
-\ds\frac\kappa2 \ds\int_{\mathbb R^d}\Big|
\ds\sum_{j=1}^P \sigma_1(x-q_j^*)\Big|^2\ud x.
\end{array}\]
Expanding the square, it becomes
\begin{equation}\label{exp_min}
\begin{array}{lll}
\mathscr H(X^*)&=&-\ds\frac\kappa 2\left(
\ds\sum_{j=1}^P\ds\int_{\mathbb R^d}
 |\sigma_1(x-q_j^*)|^2\ud x+2 \ds\sum_{ k>j}\ds\int_{\mathbb R^d}
 \sigma_1(x-q_j^*)\sigma_1(x-q_k^*)\ud x\right)
 \\&=&
- \ds\frac\kappa 2\Big(P+\ds\sum_{\ell=1}^{L} \#\Pi_\ell(\#\Pi_\ell-1)\Big)\|\sigma_1\|_{L^2}^2
=- \ds\frac\kappa 2\|\sigma_1\|_{L^2}^2\ds\sum_{\ell=1}^{L} \#\Pi_\ell^2,
 \end{array}\end{equation}
 since $\sum_{\ell=1}^{L}\#\Pi_\ell=P$.
 All these equilibrium solutions  are local minimizers of the energy $\mathscr H$, except if there exists $\bar k, \bar j\in \{1,\ldots,P\}$,
 such that $|q_{\bar k}^*-q_{\bar j}^*|=2R_1$. In this latter case, we can always find a small perturbation which lowers the energy. Indeed, let $\xi \in (\R^d)^P$ such that $\xi_i=0$ for all $i\notin \Pi_{\bar k}$ and $\xi_{k}=q_{\bar j}^*-q_{\bar k}^*$ for all $k\in \Pi_{\bar k}$, 
 and let $X=(q,p,\psi,\varpi):=(q^*+\varepsilon \xi,0,-\Gamma(z)\sum_{k=1}^P\sigma_1(x-q^*_k -\varepsilon \xi_k) ,0)$.
 Observe that
 \[\begin{array}{l}\ds\iint_{\mathbb R^n\times \mathbb R^d}\Big |\nabla\Gamma(z)\sum_{k=1}^P\big(\sigma_1(x-q^*_k -\varepsilon \xi_k )
 -\sigma_1(x-q^*_k  )\big)\Big|^2\ud z\ud x
\\[.4cm]
\hspace*{3.5cm} \leq \kappa P
 \ds\sum_{k=1}^P\ds\int_{\mathbb R^d}\Big |\sigma_1(y -\varepsilon \xi_k )
 -\sigma_1(y)\Big|^2\ud y
\xrightarrow[\varepsilon \to 0]{}0
 \end{array}\]
 (by continuity of translation in $L^2$), so that $X$ lies in a neighborhood of $X^*$ for small $\varepsilon$'s.
 For all $\varepsilon\in (0,1)$, we get $$|q_{\bar j}-q_{k}|=|q_{\bar j}^*-q_{\bar k}^*-\varepsilon(q_{\bar j}^*-q_{\bar k}^*)|=|1-\varepsilon|2R_1<2R_1,\ \forall k\in \Pi_{\bar k}$$ 
 and
 \begin{align*}
  \mathscr{H}&(X)-\mathscr{H}(X^*)\\
  &=-\frac{\kappa}{2}\left(P\|\sigma_1\|_{L^2}^2+\sum_{j=1}^P\sum_{k=1,k\neq j}^P\Sigma_1(q_j-q_k) -P\|\sigma_1\|_{L^2}^2-\sum_{j=1}^P\sum_{k=1,k\neq j}^P\Sigma_1(q_j^*-q_k^*)\right)\\
  &=-\frac{\kappa}{2}\left(\sum_{j=1}^P\sum_{k=1}^P\left(\Sigma_1(q_j-q_k) -\Sigma_1(q_j^*-q_k^*)\right)\right)\\
  &=-\frac{\kappa}{2}\left(\sum_{j\notin \Pi_{\bar k}}\sum_{k\in \Pi_{\bar k}}\left(\Sigma_1(q_j-q_k) -\Sigma_1(q_j^*-q_k^*)\right)+\sum_{j\in \Pi_{\bar k}}\sum_{k\notin \Pi_{\bar k}}\left(\Sigma_1(q_j-q_k) -\Sigma_1(q_j^*-q_k^*)\right)\right)\\
  &=-{\kappa}\left(\sum_{j\notin \Pi_{\bar k}}\sum_{k\in \Pi_{\bar k}}\left(\Sigma_1(q_j-q_k) -\Sigma_1(q_j^*-q_k^*)\right)\right)=-{\kappa}\left(\sum_{j\notin \Pi_{\bar k}}\sum_{k\in \Pi_{\bar k}}\Sigma_1(q_j-q_k)\right)\\
  &\le -{\kappa} \Sigma_1(q_{\bar j}-q_{\bar k})<0.
\end{align*}
 
Hence, we consider a subset of equilibrium solutions given by 
\begin{equation}\label{EqSol}
  \begin{aligned}
    &p^*_1=...=p^*_P=0,\qquad \psi^*(x,z)=-\Gamma(z)\ds\sum_{j=1}^P \sigma_1(x-q^*_j), \qquad \varpi^*=0,\\
    &\left\{
      \begin{aligned}
        &\{1,...,P\}=\ds\bigcup_{\ell=1}^{L} \Pi_\ell, \textrm{ for some } L\in\{1,...,P\} \textrm{ with } \ds\sum_{\ell=1}^{L} \#\Pi_\ell=P,\\
        &\Pi_\ell\cap \Pi_{\ell'}=\emptyset \textrm{ if }\ell\neq \ell',\\
        &|q^*_k-q^*_j|>2R_1 \textrm{ if } k\in \Pi_\ell,\,j\in \Pi_{\ell'},\,\ell\neq \ell',\\
        &q^*_k=q^*_j \textrm{ if } k,j\in \Pi_\ell
      \end{aligned}
    \right.
  \end{aligned}
\end{equation}
and we deduce the following proposition. 

\begin{proposition}\label{prop_energymin}
Suppose $n\geq 3$.
  Let $X^*=(q_1^*,...q^*_P,p_1^*,...,p^*_P,\psi^*,\varpi^*)$ be defined by \eqref{EqSol}. 
  It defines a stationary solution of \eqref{N_particles}-\eqref{N_wave}; it is a local minimizer of the energy functional $\mathscr H$. 
  
  Moreover, for any $\tilde q^*\in \R^d$, $X^*_{m}=(q_m^*,0,\psi^*_m,0)$ with $$q^*_m=(\tilde q^*,\ldots,\tilde q^*) \quad \textrm{  and } \quad 
  \psi_m^*(x,z)=-\Gamma(z)\ds P \sigma_1(x-\tilde q^*)$$ is a global minimizer of $\mathscr H$.

  Finally, for all $P_0\in \R$,
  \begin{equation}
    \label{min_fixed_momentum}
    \mathscr H(X^*_m)=\inf\{\mathscr{H}(X),\ X\in (\R^d)^P\times (\R^d)^P\times L^2(\R^d,\dot{H}^1(\R^d))\times L^2(\R^d\times \R^n), \mathscr{P}(X)=P_0\},
  \end{equation}
  with $\mathscr{P}$ defined by \eqref{defmt}, and the infimum is attained if and only if $P_0=0$.
\end{proposition}

\noindent
{\bf Proof.}  By definition, $X^*=(q_1^*,...q^*_P,p_1^*,...,p^*_P,\psi^*,\varpi^*)$ given by \eqref{EqSol} is a stationary solution of \eqref{N_particles}-\eqref{N_wave}. As a consequence, using the Hamiltonian structure of the problem, we deduce that $X^*$ is a critical point of $\mathscr H$. Let us show that $X^*$ is a local minimizer.

  First of all, we remark that, for a given $(q_1,...,q_P)\in (\mathbb R^d)^P$, the function
  \[
  \begin{array}{llll}
  &(\mathbb R^d)^P\times L^2(\mathbb R^d; \dot H^1(\mathbb R^n))\times L^2(\mathbb R^d\times\mathbb R^n)
  &\longrightarrow & \mathbb R
  \\
  &(p_1,...p_P,\psi,\varpi)
  &  \longmapsto 
  &\hat {\mathscr H}_{q_1,...,q_P}(p_1,...,p_P,\psi,\varpi)
  \\
  &&&\qquad=\mathscr H(q_1,...,q_P,p_1,...,p_P,\psi,\varpi)
  \end{array}\]
  is convex. A straigthforward computation gives that the minimizer satisfies
  \begin{equation*}
    p_1=...=p_P=0,\qquad \psi(x,z)=-\Gamma(z)\sum_{j=1}^P\sigma_1(x-q_j) \text{ and } \varpi=0
  \end{equation*}
  and the minimum value is
  \begin{equation}\label{Hm}
  \mathscr H_m(q_1,...,q_P)=-\ds\frac\kappa2
  \ds\int_{\mathbb R^d} \Big|\ds\sum_{j=1}^P \sigma_1(x-q_j)\Big|^2\ud x=-\frac{\kappa}{2}\left(P\|\sigma_1\|_{L^2}^2+\sum_{j=1}^P\sum_{k=1,k\neq j}^P\Sigma_1(q_j-q_k)\right)
  .\end{equation}
  As a consequence, for any $X\in (\R^d)^P\times (\R^d)^P\times L^2(\R^d,\dot{H}^1(\R^d))\times L^2(\R^d\times \R^n)$,
  \begin{equation*}
    \mathscr H(X)\ge \mathscr H_m(q_1,...,q_P).
  \end{equation*}
  When there is only one particle $P=1$, we deduce that 
   $\mathscr H_m(q)=-\frac\kappa 2\|\sigma_1\|^2_{L^2}$ does not depend on $q$ and 
   for any $X\in \R^d\times \R^d\times L^2(\R^d,\dot{H}^1(\R^d))\times L^2(\R^d\times \R^n)$,
  \begin{equation*}
    \mathscr H(X)\ge \mathscr H_m(q)=\mathscr H(X^*)
  \end{equation*}
  thanks to~\eqref{exp_min}.
    The case with many particles 
    is more delicate.
   Critical points of \eqref{Hm}
   are requested to satisfy $\nabla_{q} \mathscr H_m(q_1,...,q_P)=0$. More precisely, using the properties of $\Sigma_1$, we obtain, for any $\xi=(\xi_1,\ldots,\xi_P)\in (\R^d)^P$,
   \begin{align*}
    D\mathscr{H}_m(q_1,\ldots,q_P)(\xi)&=-\frac{\kappa}{2}\sum_{j=1}^P\sum_{k=1}^P\nabla\Sigma_1(q_j-q_k)\cdot (\xi_j-\xi_k)=-\kappa\sum_{j=1}^P\sum_{k=1}^P\nabla\Sigma_1(q_j-q_k)\cdot \xi_j\\
    &= \langle \nabla_{q} \mathscr H_m(q_1,...,q_P),\xi\rangle
   \end{align*}
   where $\cdot$ denotes the scalar product in $\R^d$, $\langle q, \xi\rangle= \sum_{j=1}^P q_j\cdot\xi_j$ and 
   \begin{equation*}
    \nabla_{q_j} \mathscr H_m(q_1,...,q_P)=-\kappa\sum_{k=1}^P \nabla\Sigma_1(q_j-q_k).
   \end{equation*}
   
   As discussed above (see Proposition~\ref{prop_eq_solutions}), for any $j\in \{1,\ldots,P\}$, the solutions to the equation
   \begin{equation*}
    -\kappa\sum_{k=1}^P \nabla\Sigma_1(q_j-q_k)=0
   \end{equation*} 
   can be described as sets of  $L\in \{1,...,P\}$ clusters
   containing, for  $\ell\in \{1,..,L\}$, $\#\Pi_\ell$ particles all with the same location. Hence $q^*=(q_1^*,\ldots,q_p^*)$ given by \eqref{EqSol} is a critical point of $\mathscr{H}_m$.
   
   The second order derivative of $\mathscr H_m$ is given, for any $\xi=(\xi_1,...,\xi_P)\in (\mathbb R^d)^P$, by
   \begin{equation*}
    D^2\mathscr H_m(q_1,...,q_P)(\xi,\xi)=-\frac{\kappa}{2}\sum_{j=1}^P\sum_{k=1}^P (\xi_j-\xi_k)\cdot (H_{\Sigma_1}(q_j-q_k)(\xi_j-\xi_k))
   \end{equation*}
   where $H_{\Sigma_1}(q_j-q_k)$ is the Hessian matrix of the function $\Sigma_1$ evaluated at the point $q_j-q_k$.

  Considering a critical point given by~\eqref{EqSol} and using the fact that $H_{\Sigma_1}(q_j-q_k)=0$ if $q_j$ and $q_k$ are not in the same cluster, we get
  \begin{align*}
    D^2\mathscr H_m(q^*)(\xi,\xi)&=-\frac{\kappa}{2}\sum_{\ell=1}^{L}\sum_{j\in \Pi_\ell}\sum_{m=1}^{L}\sum_{k\in \Pi_m} (\xi_j-\xi_k)\cdot (H_{\Sigma_1}(q_j-q_k)(\xi_j-\xi_k))\\
    &= -\frac{\kappa}{2}\sum_{\ell=1}^{L}\sum_{j,k\in \Pi_\ell}(\xi_j-\xi_k)\cdot (H_{\Sigma_1}(0)(\xi_j-\xi_k))
  \end{align*}
  Since $\sigma_1$ is radially symmetric, a straigthforward computation allows us to conclude that $$(-H_{\Sigma_1}(0))=\alpha \mathbb I_d$$ with 
  \begin{equation}\label{def_alpha}
    \alpha=|B(0,1)|\int_{0}^{+\infty}(\sigma_1'(r))^2r^{d-1}\ud r.
  \end{equation}
  As a consequence, 
  \begin{align*}
    D^2\mathscr H_m(q^*)(\xi,\xi)&=\frac{\alpha \kappa}{2} \sum_{\ell=1}^{L}\sum_{j,k\in \Pi_\ell}(\xi_j-\xi_k)\cdot (\xi_j-\xi_k)
    \\
    &=\alpha \kappa\sum_{\ell=1}^{L}\#\Pi_\ell\sum_{j\in \Pi_\ell}(\xi_j\cdot\xi_j-\xi_j\cdot\langle \xi \rangle_{\ell})=\alpha \kappa\sum_{\ell=1}^{L}\#\Pi_\ell\sum_{j\in \Pi_\ell}\xi_j\cdot(\xi_j-\langle \xi \rangle_{\ell})\\
    &=\alpha \kappa \sum_{\ell=1}^{L}\#\Pi_\ell\sum_{j\in \Pi_\ell}\left((\xi_j-\langle \xi \rangle_{\ell})\cdot(\xi_j-\langle \xi \rangle_{\ell})+ \langle \xi \rangle_{\ell}\cdot(\xi_j-\langle \xi \rangle_{\ell})\right)\\
    &=\alpha \kappa \sum_{\ell=1}^{L}\#\Pi_\ell\sum_{j\in \Pi_\ell}(\xi_j-\langle \xi \rangle_{\ell})\cdot(\xi_j-\langle \xi \rangle_{\ell})
  \end{align*}
  where $\langle \xi \rangle_{\ell}=\frac{1}{\# \Pi_\ell}\sum_{k\in \Pi_\ell}\xi_k\in \R^d$ for all $\ell\in\{1,\ldots,L\}$.
  This implies that $D^2\mathscr H_m(q^*)(\xi,\xi)$ 
  is non-negative for any $\xi\in (\mathbb R^d)^P$ and it vanishes if and only if, for all $\ell\in \{1,\ldots, L\}$, $\xi_j=\langle \xi \rangle_{\ell}$ for all $j\in \Pi_\ell$, \emph{i.e.} $\xi$ is blockwise constant according to  the $\Pi_\ell$'s. Let $\mathcal P:(\R^d)^P\to (\R^d)^P$ the application defined by $\mathcal P\xi=((\mathcal P\xi)_1,\ldots,(\mathcal P\xi)_P)$ and 
  \begin{equation*}
    (\mathcal P\xi)_j=\langle \xi \rangle_{\ell}
  \end{equation*}
  for all $j\in \Pi_\ell$, $\ell\in \{1,\ldots, L\}$. Clearly $\mathcal P$ is a symmetric and idempotent endomorphism. Hence, the map $\mathcal P$ is an orthogonal projection on $(\R^d)^P$, and $D^2\mathscr H_m(q^*)(\xi,\xi)$ vanishes if and only if $\xi\in \ker(I-\mathcal P)$. Moreover, any $\xi\in (\R^d)^P$ can be written as 
  \begin{equation*}
    \xi= \mathcal P\xi+(I-\mathcal P)\xi:=\xi_\parallel+\xi_\perp
  \end{equation*} 
  with $\xi_\parallel\in \mathrm{im}(\mathcal P)$ and $\xi_\perp \in \ker(\mathcal P)$.
  Since $q^*$ is such that $|q_j^*-q_k^*|>2R_1$ whenever $j$ and $k$ are not in the same cluster, we deduce that $\mathscr{H}_m(q^*)=\mathscr{H}_m(q^*+\varepsilon\xi_\parallel)$ for any $\xi_\parallel \in \mathrm{im}(\mathcal P)$, provided $\varepsilon$ is small enough. As a consequence, 
  \begin{equation*}
    D^2\mathscr H_m(q^*)(\xi_\perp,\xi_\perp)=D^2\mathscr H_m(q^*+\varepsilon\xi_\parallel)(\xi_\perp,\xi_\perp)=\alpha \kappa \sum_{\ell=1}^{L}\#\Pi_\ell\sum_{j\in \Pi_\ell}(\xi_\perp)_j\cdot(\xi_\perp)_j
  \end{equation*}
  for any $\xi_\perp\in \ker(\mathcal P)$. Therefore, for any $\varepsilon$ small enough and for any $\xi\in (R^d)^P$, we get
  \begin{align*}
    \mathscr H_m(q^*+\varepsilon\xi)-\mathscr H_m(q^*)&=\mathscr H_m(q^*+\varepsilon\xi)-\mathscr H_m(q^*+\varepsilon\xi_\parallel)\\
    &=\varepsilon^2\frac{1}{2}D^2\mathscr H_m(q^*+\varepsilon\xi_\parallel)(\xi_\perp,\xi_\perp)+o(\varepsilon^2\|\xi_\perp\|^2)\ge 0
  \end{align*}
  and $q^*$ is a local minimimizer for $\mathscr H_m$. We conclude that  $X^*$ is a local minimizer for $\mathscr H$.

Moreover, since $\Sigma_1$ is radially symmetric and attains its maximum at $x=0$, we have, for any $X\in (\R^d)^P\times (\R^d)^P\times L^2(\R^d,\dot{H}^1(\R^d))\times L^2(\R^d\times \R^n)$,
\begin{align*}
  \mathscr H(X)\ge \mathscr H_m(q_1,...,q_P)\ge -\frac{\kappa}{2}\left(P\|\sigma_1\|_{L^2}^2+\sum_{j=1}^P\sum_{k=1,k\neq j}^P\Sigma_1(0)\right)=-\frac{\kappa}{2}P^2\|\sigma_1\|_{L^2}^2=\mathscr H(X^*_{m})
\end{align*}
with $X^*_{m}=(q_m^*,0,\psi^*_m,0)$ with $q^*_m=(\tilde q^*,\ldots,\tilde q^*)$ and $\psi_m^*(x,z)=-\Gamma(z)\ds P \sigma_1(x-\tilde q^*)$ for any $\tilde q^*\in \R^d$.
As a consequence, $X^*_m$ is a global minimizer of $\mathscr{H}$.

It remains to show that~\eqref{min_fixed_momentum} is satisfied for any $P_0\in \R$ and the infimum is attained if and only if $P_0=0$.

Let $X^*_m$ be a global minimizer of $\mathscr{H}$. First of all, we observe that the momentum  associated to these solutions vanishes, \emph{i.e.} $\mathscr{P}(X_m^*)=0$ Next, to show the result, we adapt the reasoning in   \cite[Prop.~1.2]{LB_PhD}, which relies on the fact that
$\psi^*$ is supported in $\mathscr C= B(\tilde q^*,R_1)\times \mathbb R^n$ and $\sigma_2$ is supported in $B(0,R_2)$.

The idea is to construct, for $\epsilon>0$, a perturbation $\widetilde X_\epsilon=(0,0,\widetilde \psi_\epsilon, \widetilde
\varpi_\epsilon)$ such that 
\begin{equation*}
  \mathscr H(X^*_m+\widetilde X_\epsilon)=\mathscr H(X_m^*)+\varepsilon \text{ and } \mathscr P(X^*_m+\widetilde X_\epsilon)=P_0.
\end{equation*} 
In particular, we choose $\widetilde \psi_\epsilon, \widetilde \varpi_\epsilon$ such that their supports in $x$ are localized outside $B(\tilde q^*,R_1)$ and in $z$ are localised outside $B(0,R_2)$.
As consequence, 
\begin{equation*}
  \mathscr H(X^*_m+\widetilde X_\epsilon)=\mathscr H(X_m^*)+\mathscr{H}(\widetilde X_\epsilon)=\mathscr H(X_m^*)+\frac{1}{2}\iint |\widetilde\varpi_\epsilon|^2\ud z\ud x+ \frac{1}{2}\iint |\nabla_z\widetilde\psi_\epsilon|^2\ud z\ud x
\end{equation*}
and 
\begin{equation*}
  \mathscr P(X^*_m+\widetilde X_\epsilon)=\mathscr P(X_m^*)+\mathscr{P}(\widetilde X_\epsilon)=-\frac{1}{c}\iint \widetilde\varpi_\epsilon \nabla_x \widetilde\psi_\epsilon\ud z\ud x.
\end{equation*}
Hence, for some $\ell_\epsilon>0$ and $x_\epsilon\in \mathbb R^d$, we set 
\[\widetilde\psi_\epsilon(x,z)=
\ell_\epsilon^{(1-d)/2}\zeta_1((x-x_\epsilon)/\ell_\epsilon)\zeta_2(z)
,\qquad
\widetilde\varpi_\epsilon(x,z)=\ell_\epsilon^{(1-d)/2}\tilde \zeta_1((x-x_\epsilon)/\ell_\epsilon)\zeta_2(z)\]
where 
$\zeta_1, \tilde \zeta_1,\zeta_2$ are $C^\infty$-smooth and compactly supported with $\mathrm{supp}(\zeta_1)\cap B(\tilde q^*,R_1)=\emptyset$, $\mathrm{supp}(\tilde \zeta_1)\cap B(\tilde q^*,R_1)=\emptyset$, $\mathrm{supp}(\zeta_2)\cap B(0,R_2)=\emptyset$ and  $\int \zeta_2^2(z)\ud z=1$. 
Next, we can choose $\zeta_1$, $\tilde \zeta_1$ to impose
\[\begin{array}{lll}
\ds\iint \widetilde\varpi_\epsilon\nabla_x \widetilde\psi_\epsilon\ud z\ud x
&=&\ell_\epsilon^{1-d}\ds\iint 
|\zeta_2(z)|^2\tilde\zeta_1((x-x_\epsilon)/\ell_\epsilon)
\ds\frac{(\nabla_x \zeta_1)((x-x_\epsilon)/\ell_\epsilon)}{\ell_\epsilon}\ud z\ud x
\\&=&\ds\int \tilde \zeta_1(y)
\nabla_x \zeta_1(y)\ud y=-c P_0,\end{array}\]
while we get
\[\begin{array}{l}
\ds\frac{1}{2}\iint |\widetilde\varpi_\epsilon|^2\ud z\ud x+ \frac{1}{2}\iint |\nabla_z\widetilde\psi_\epsilon|^2\ud z\ud x
\\\qquad=
\ds\frac{\ell_\epsilon^{1-d}}{2}\left(
\ds\int 
|\tilde \zeta_1((x-x_\epsilon)/\ell_\epsilon)|^2\ud x+\|\nabla\zeta_2\|^2_{L^2}\ds\int 
|\zeta_1((x-x_\epsilon)/\ell_\epsilon)|^2\ud x\right)
\\\qquad=
\ds\frac{\ell_\epsilon}{2}\Big(
\ds\|\tilde \zeta_1\|_{L^2}^2+\|\nabla\zeta_2\|^2_{L^2}\|\zeta_1\|_{L^2}^2
\Big).
\end{array}\]
We thus pick  $\ell_\epsilon$ small enough, and $x_\epsilon$ so that the support requirement is fulfilled.
\QED

It is worth pointing out that the 
set of minimizers on a given energy level
is not isolated, nor bounded: we can modify an equilibrium
state by any translation or rotation without changing
the value of the energy.

\subsection{Perturbations of least energy equilibria}

In what follows, we will naturally focus on perturbations of an equilibrium state with least energy, \emph{i.e.} a state of the form $X^*_{m}=(q_m^*,0,\psi^*_m,0)$ with
\begin{equation}
  \label{eq_defequilibrium_least}
   q^*_m=(\tilde q^*,\ldots,\tilde q^*)\in \big(\mathbb R^d\big)^P \text{ and } \psi_m^*(x,z)=-\Gamma(z)\ds P \sigma_1(x-\tilde q^*)
\end{equation}
for any $\tilde q^*\in \R^d$. We are going to prove that if we start from a small enough perturbation of this state, then the distance between particles remains bounded for ever.

It is convenient to introduce the following notation: to $q=(q_1,...,q_P)\in (\mathbb R^d)^P$, we associate 
$$
\psi_q(x,z)=-\Gamma(z)\ds\sum_{j=1}^P\sigma_1(x-q_j).$$
We have seen that, for any $X\in (\R^d)^P\times (\R^d)^P\times L^2(\R^d,\dot{H}^1(\R^d))\times L^2(\R^d\times \R^n)$, we have
\begin{equation}\label{ineqHs}
\mathscr H(X)
\geq \mathscr H(q,0,\psi,0)
\geq \mathscr H_m(q)=\mathscr H(q,0,\psi_q,0)
\geq \mathscr H(q_m^*,0,\psi_m^*,0)
=-\ds\frac\kappa 2 \|\sigma_1\|_{L^2}^2 P^2,\end{equation}
where $\mathscr H_m(q)$ is explicitly expressed in \eqref{Hm}.
Let us now introduce the set 
\[
\mathscr E=\big\{q\in (\mathbb R^d)^P \textrm{ such that } \exists (j,k)\in \{1,...,P\},\ |q_j-q_k|\geq 2R_1\big\}.\]
Let $q\in \mathscr E$. Without loss of generality, suppose that there exists $j_0\in \{1,...,P\}$ such that $|q_1-q_{j_0} |\geq  2R_1$.
From \eqref{Hm}, we get
\[
\mathscr H_m(q)=-\frac{\kappa P}{2}\|\sigma_1\|_{L^2}^2-
\frac{\kappa }{2}\sum_{k\notin\{1, j_0\}}\Sigma_1(q_1-q_k)
-\frac{\kappa }{2}\sum_{j=2}^P\sum_{k=1,k\neq j}^P\Sigma_1(q_j-q_k)
\]
which is thus always larger than
\begin{equation}
  \label{eq_defH2}
  \mathscr H^*_2:=
-\frac{\kappa }{2}\|\sigma_1\|_{L^2}^2(P+P-2+(P-1)^2)
=-\frac{\kappa }{2}\|\sigma_1\|_{L^2}^2(P^2-1)>
\mathscr H(q_m^*,0,\psi_m^*,0).
\end{equation}
Using, $-\Delta_z \Gamma=\sigma_2$, we can write
\begin{align*}
\ds\sum_{j=1}^P \ds\iint_{\mathbb R^d\times\mathbb R^n} \sigma_2(z)\sigma_1(x-q_j)\psi(x,z)\ud z\ud x
&=\ds\iint_{\mathbb R^d\times\mathbb R^n} \Delta_z \psi_q(x,z) \psi(x,z)\ud z\ud x
\\&=-\ds\iint_{\mathbb R^d\times\mathbb R^n} \nabla_z \psi_q(x,z)\cdot \nabla_z\psi(x,z)\ud z\ud x
\end{align*}
so that 
\begin{align*}
&\mathscr H(q,0,\psi,0)- \mathscr H_m(q)
=\ds\frac{1}{2}
\ds\iint_{\mathbb R^d\times\mathbb R^n} |\nabla_z \psi|^2\ud z\ud x
+\ds\sum_{j=1}^P \ds\iint_{\mathbb R^d\times\mathbb R^n} \sigma_2(z)\sigma_1(x-q_j)\psi(x,z)\ud z\ud x
\\&\hspace*{4cm} -\ds\frac{1}{2}
\ds\iint_{\mathbb R^d\times\mathbb R^n} |\nabla_z \psi_q|^2\ud z\ud x
-\ds\sum_{j=1}^P \ds\iint_{\mathbb R^d\times\mathbb R^n} \sigma_2(z)\sigma_1(x-q_j)\psi_q(x,z)\ud z\ud x
\\
&\ =\ds\frac{1}{2}
\ds\iint_{\mathbb R^d\times\mathbb R^n} |\nabla_z \psi|^2\ud z\ud x
-\ds\iint_{\mathbb R^d\times\mathbb R^n} \nabla_z \psi_q(x,z) \cdot\nabla_z\psi(x,z)\ud z\ud x
+\ds\frac{1}{2}
\ds\iint_{\mathbb R^d\times\mathbb R^n} |\nabla_z \psi_q|^2\ud z\ud x
\\
&\ =\ds\frac{1}{2}
\ds\iint_{\mathbb R^d\times\mathbb R^n} |\nabla_z (\psi-\psi_q)|^2\ud z\ud x.
\end{align*}
Therefore, we have
in fact
a more precise version of \eqref{ineqHs}
\[
\mathscr H(X)=
\mathscr H(q,0,\psi_q,0)
+
\underbrace{\ds\sum_{j=1}^P \ds\frac{|p_j|^2}{2}
+\ds\frac{1}{2}
\ds\iint_{\mathbb R^d\times\mathbb R^n} |\varpi|^2\ud z\ud x
+\ds\frac{1}{2}
\ds\iint_{\mathbb R^d\times\mathbb R^n} |\nabla_z (\psi-\psi_q)|^2\ud z\ud x}_{:=
\mathscr G(X)\geq 0}.
\]

The energetic properties can be used to derive a priori estimates on the solutions.
In particular we can quantify energy gaps in order to  guaranty 
that the particles remain in a single cluster for all times.

\begin{proposition}\label{bdd_dist}
Suppose $n\geq 3$.
Let $$X_0=(q_0,p_0,\Psi_0,\Psi_1) \in (\mathbb R^d)^P\times (\mathbb R^d)^P\times L^2(\R^d,\dot{H}^1(\R^d))\times L^2(\R^d\times \R^n)$$
be an initial data for the system \eqref{N_particles}-\eqref{N_wave}.
Then, the solution of \eqref{N_particles}-\eqref{N_wave} is such that the quantities
$t\mapsto p_j(t)$, $t\mapsto \|\nabla_z\psi(t,\cdot)\|_{L^2}$,
 $t\mapsto \|\partial_t\psi(t,\cdot)\|_{L^2}$ remain bounded for any $t\geq 0$. 

Moreover, if we assume that $X_0$ is a small perturbation 
of a least energy equilibrium state $X^*_{m}=(q_m^*,0,\psi^*_m,0)$ of the form~\eqref{eq_defequilibrium_least},
in the sense that 
\[
-\ds\frac\kappa 2 \|\sigma_1\|_{L^2}^2 P^2
\leq \mathscr H(X_0)
<\mathscr H^*_2\]
with $\mathscr H^*_2$ defined by~\eqref{eq_defH2},
then
 $|q_j(t)-q_k(t)|<2R_1$ holds for any $t\geq 0$, $j,k\in\{1,...,P\}$.
\end{proposition}

\begin{rmk}
The condition on the data can be rephrased as a smallness condition under the form
\[
0\leq \ds\frac{\kappa}2\ds\sum_{k=1}^P\ds\sum_{j\neq k} \big(\Sigma_1(0)-\Sigma_1(q_{0,j}-q_{0,k})\big)
+\mathscr G(X_0)< 
\ds\frac{\kappa}2\|\sigma_1\|^2_{L^2}.
\]
\end{rmk}

\noindent
{\bf Proof.}
Existence-uniqueness of a globally defined solution to~ \eqref{N_particles}-\eqref{N_wave} has been justified in \cite{dBGV} 
The asserted estimates are obtained as consequences of the energy conservation.
In the energy conservation, the only term which does not have a definite sign is the coupling energy. However, it satisfies 
\begin{align*}
&\left|\ds\sum_{k=1}^P\ds\iint_{\mathbb R^d\times\mathbb R^n}
\sigma_2(z)\psi(t,x,z)\sigma_1(x-q_k(t))\ud z\ud x\right|
\\&=
\left|\ds\sum_{k=1}^P\ds\iint_{\mathbb R^d\times\mathbb R^n}
\ds\frac{\widehat \sigma_2(\xi )}{|\xi|} |\xi| \widehat \psi(t,x,\xi)\sigma_1(x-q_k(t))\ds\frac{\ud \xi}{(2\pi)^n}\ud x\right|
\\
&\leq
\left(\ds\int_{\mathbb R^n}
\ds\frac{|\widehat \sigma_2(\xi )|^2}{|\xi|^2}\ds\frac{\ud \xi}{(2\pi)^n}
\right)^{1/2}
\ds\sum_{k=1}^P\ds\int_{\mathbb R^d}
\left(\ds\int_{\mathbb R^n}
 |\xi|^2 |\widehat \psi(t,x,\xi)|^2\ds\frac{\ud \xi}{(2\pi)^n}\right)^{1/2} \sigma_1(x-q_k(t))\ud x
\\&\leq
\sqrt\kappa \ds\sum_{k=1}^P
\left(\ds\int_{\mathbb R^d} |\sigma_1(x-q_k(t))|^2\ud x\right)^{1/2}
\left(\ds\iint_{\mathbb R^d\times \mathbb R^n}
 |\xi|^2 |\widehat \psi(t,x,\xi)|^2\ds\frac{\ud \xi}{(2\pi)^n}\ud x\right)^{1/2}
\\& \leq 
 P\|\sigma_1\|_{L^2}\sqrt{\kappa }
 \left(\ds\iint_{\mathbb R^d\times \mathbb R^n}
|\nabla_z\psi(t,x,z)|^2\ud z\ud x\right)^{1/2}.
 \end{align*}
 By energy conservation, we have
 \begin{align*}
 \mathscr H(X(0))&=\mathscr H(X(t))
=\ds\sum_{j=1}^P \ds\frac{p_j(t)}{2}
 +\ds\frac{1}{2c^2}\ds\iint_{\mathbb R^d\times\mathbb R^n}|\partial_t \psi(t,x,z)|^2\ud z\ud x
\\&\quad +\ds\frac{1}{2}\ds\iint_{\mathbb R^d\times\mathbb R^n}|\nabla_z \psi(t,x,z)|^2\ud z\ud x
 +\ds\sum_{k=1}^P\ds\iint_{\mathbb R^d\times\mathbb R^n}
\sigma_2(z)\psi(t,x,z)\sigma_1(x-q_k(t))\ud z\ud x
\\
&\geq 
\ds\sum_{j=1}^P \ds\frac{p_j(t)}{2}
 +\ds\frac{1}{2c^2}\ds\iint_{\mathbb R^d\times\mathbb R^n}|\partial_t \psi(t,x,z)|^2\ud z\ud x
 \\&
 \hspace*{4cm}+\ds\frac{1}{4}\ds\iint_{\mathbb R^d\times\mathbb R^n}|\nabla_z \psi(t,x,z)|^2\ud z\ud x
 -P^2\|\sigma_1\|^2_{L^2}\frac{\kappa}{2} 
 \end{align*}
which provides the uniform estimate on $t\mapsto \sum_{j=1}^P |p_j(t)|^2$.
We also know that $t\mapsto \partial_t \psi(t,\cdot)$
and $t\mapsto \nabla_z \psi(t,\cdot)$ are bounded in $L^\infty((0,\infty);L^2(\mathbb R^d\times\mathbb R^n))$.
Next, the energy conservation and the smallness condition tell us that
\[
-\ds\frac{\kappa P^2}2\|\sigma_1\|^2_{L^2}\leq \mathscr H(X(t))=\mathscr H(X_0)<\mathscr H^*_2.\]
Suppose, by contradiction, that there exists $\tilde t\in (0,\infty)$ and $(k,j)\in \{1,\ldots,P\}$, such that $|q_j(\tilde t)-q_k(\tilde t)|>2R_1$. From~\eqref{ineqHs} and~\eqref{eq_defH2}, we deduce that
\begin{align*}
  \mathscr H(X(\tilde t))\ge \mathscr H_m(q(\tilde t))\ge \mathscr H^*_2.
\end{align*}
This leads to a contradiction and we conclude that,
under the smallness condition, $|q_j(t)-q_k(t)|$ cannot exceeds $2R_1$ for any $j,k$ and for all $t\in (0,\infty)$.
\QED

\subsection{Linearized problem and spectral stability}

A starting point to understand the behavior of solutions to~\eqref{N_particles}--\eqref{N_wave} is the analysis of the linearized system around $X^*$ given by~\eqref{EqSol}, which reads as
\begin{equation}\label{Lin}
  \left\{
\begin{array}{l}
\dot{q}_j(t)=p_j(t),
\\
\dot{p}_j(t)=\ds\iint_{\mathbb R^d\times\mathbb R^n}\sigma_2(z)
\big(\nabla_x\sigma_1(y-q_j^*)\cdot q_j\nabla_x \psi^*(y,z)
-
\sigma_1(y-q_j^*)\nabla_x \psi(t,y,z)\big)\ud z\ud y
\\
\partial_t \psi(t,x,z)=\varpi(t,x,z)
\\
\partial_t \varpi(t,x,z)=c^2\left(
\Delta_z \psi(t,x,z) + 
\sigma_2(z)\ds\sum_{k=1}^P\nabla_x\sigma_1(x-q^*_k)\cdot q_k\right).
\end{array}
  \right.
\end{equation}
It will  be referred to with the shorthand notation $\partial_ t X=\mathbb LX$. Our goal is to investigate the spectrum of $\mathbb L$.

In fact, it is convenient at this stage to change variables by setting:
\[
\varphi=(-\Delta)^{1/2}\psi.\]
As observed in \cite{SRN3} this change of variable has the advantage of unifying the functional framework and to ease the identification 
of the key spectral properties.
Denoting $Y=(q_1,...,q_P,p_1,...,p_P, \varphi,\varpi)$, and $\varphi^*=(-\Delta)^{1/2}\psi^*$,
the system under consideration recasts as $\partial_ t Y(t,x,z)=\tilde L Y$ with
\[
\tilde L Y=
\begin{pmatrix}
p_1 (t) \\
\vdots
\\
p_P(t)
\\
\ds\iint_{\mathbb R^d\times\mathbb R^n} 
\tilde \sigma_2(z)\big(\nabla_x\sigma_1(y-q_1^*)\cdot
q_1\nabla_x \varphi^*(y,z)
+\nabla_x\sigma_1(y-q_1^*)\varphi(t,y,z)
\big)\ud z\ud y
\\
\vdots
\\
\ds\iint_{\mathbb R^d\times\mathbb R^n} 
\tilde \sigma_2(z)\big(\nabla_x\sigma_1(y-q_P^*)\cdot
q_P\nabla_x \varphi^*(y,z)
+\nabla_x\sigma_1(y-q_P^*)\varphi(t,y,z)
\big)\ud z\ud y\\
(-\Delta)^{1/2}\varpi(t,x,z)
\\
-c^2(-\Delta)^{1/2}\varphi(t,x,z) + c^2 \sigma_2(z)\ds\sum_{k=1}^P \nabla_x\sigma_1(x-q_k^*)\cdot q_k
\end{pmatrix},
\]  
where $\tilde \sigma_2=(-\Delta)^{-1/2}\sigma_2$.
Note that  we can write
\[\begin{array}{l}
\ds\iint_{\mathbb R^d\times\mathbb R^n} 
\nabla_x\sigma_1(x-q^*_j)\cdot q_j \tilde \sigma_2(z)\nabla_x\varphi^*(x,z)\ud z\ud x
\\\qquad=
-\ds\iint_{\mathbb R^d\times\mathbb R^n} 
\nabla_x\sigma_1(x-q^*_j)\cdot q_j \tilde \sigma_2(z)(-\Delta)^{1/2}\Gamma(z)
\ds\sum_{k=1}^P
\nabla_x\sigma_1(x-q^*_k)\ud z\ud x
\\\qquad=
-\kappa\alpha r_j q_j,
\end{array}\]
with $\alpha$ defined by~\eqref{def_alpha}.
We thus work with the functional space
$$\tilde{\mathbb H}=(\mathbb R^d)^{2P}\times L^2(\mathbb R^d\times\mathbb R^n)
\times L^2(\mathbb R^d\times\mathbb R^n).$$
Of course, the relation $\mathbb LX=\lambda X$ and $\tilde LY=\lambda Y$ are equivalent, up to the change of variables.

\begin{proposition}\label{prop_spectrum_lin}\footnote{In this statement we call ``point spectrum'' the set of 
$\lambda\in \mathbb C$ such that $\mathrm{Ker}(\mathbb L-\lambda)\neq \{0\}$ is finite dimensional
and $\sigma_{\mathrm{ess}}(\mathbb L)$ is defined by means of Weyl sequences.
We refer to \cite{EdEv} for further details in spectral theory.
}
Let us suppose that $n\geq 3$ and $\widehat \sigma_2(\xi)\neq 0$ holds for any $\xi\in \mathbb R^n$.
The pure point spectrum of $\tilde  L$ is given by 
$$
\{0\}\cup \big\{\pm i\sqrt{\kappa \alpha\#\Pi_\ell}, \ \ell\in\{1,...,L\},\ \#\Pi_\ell>1\big\}$$  
and  $\sigma_{\mathrm{ess}}(\tilde L)=i\mathbb R$.  
The equilibrium solutions $X^*=(q_1^*,...q^*_P,p_1^*,...,p^*_P,\psi^*,\varpi^*)$ defined by \eqref{EqSol} are thus spectrally stable. 
\end{proposition}

\noindent
{\bf Proof.}
Let us start by considering the associated eigenvalue problem
\begin{equation}\label{EV}
  \left\{
\begin{array}{l}
\lambda q_j=p_j,
\\
\lambda p_j=\ds\iint_{\mathbb R^d\times\mathbb R^n}\tilde\sigma_2(z)
\big(\nabla_x\sigma_1(y-q_j^*)\cdot q_j\nabla_x \varphi^*(y,z)+
\nabla_x\sigma_1(y-q_j^*) \varphi(y,z)\big)\ud z\ud y
\\
\lambda \varphi=(-\Delta)^{1/2}\varpi,
\\
\lambda \varpi(x,z)=-c^2(-\Delta)^{1/2}\varphi(x,z)+c^2\sigma_2(z)
\ds\sum_{k=1}^P\nabla_x\sigma_1(x-q^*_k)\cdot q_k
\end{array}
  \right.
\end{equation}
It yields
\begin{equation}\label{EVqj}
\lambda^2 q_j=\ds\iint_{\mathbb R^d\times\mathbb R^n}
\tilde\sigma_2(z)
\big(\nabla_x\sigma_1(y-q_j^*)\cdot q_j\nabla_x \varphi^*(y,z)
+\nabla_x\sigma_1(y-q_j^*) \varphi(y,z)\big)\ud z\ud y
\end{equation}
together with 
\[
\Big(\ds\frac{\lambda^2}{c^2}-\Delta_z\Big)\varphi(x,z)=(-\Delta)^{1/2}\sigma_2(z)
\ds\sum_{k=1}^P\nabla_x\sigma_1(x-q^*_k)\cdot q_k.
\]
The latter relation can be interpreted by means of Fourier transform with respect to the $z$ variable
\begin{equation}\label{FouPsi}
\widehat \varphi(x,\xi)=\ds\frac{|\xi|\, \widehat \sigma_2(\xi)}{\lambda^2/c^2+\xi^2}
\ds\sum_{k=1}^P\nabla_x\sigma_1(x-q^*_k)\cdot q_k
,\end{equation}
which can be plugged into \eqref{EVqj}
\[\begin{array}{lll}
\lambda^2 q_j&=&
\ds\int_{\mathbb R^d} \nabla_x\sigma_1(y-q_j^*)\cdot q_j
\left(\ds\int_{\mathbb R^n} \tilde\sigma_2(z)\nabla_x \varphi^*(y,z)\ud z\right)\ud y
\\&&+
\ds\sum_{k=1}^P\ds\iint_{\mathbb R^d\times\mathbb R^n}
\ds\frac{ |\widehat \sigma_2(\xi)|^2}{\lambda^2/c^2+\xi^2}
\nabla_x\sigma_1(y-q_j^*)\nabla_x \sigma_1(y-q^*_k)\cdot q_k\ds\frac{\ud \xi}{(2\pi)^n}\ud y.
\end{array}\]
Let us introduce the function
\[z\in\mathbb C
\longmapsto \kappa_z=\ds\int_{\mathbb R^n}
\ds\frac{|\widehat \sigma_2(\xi)|^2}{z+\xi^2}
\ds\frac{\ud \xi}{(2\pi)^n}.\]
Hence, if $\lambda$ is an eigenvalue of $\tilde L$, there exists $q=(q_1,\ldots,q_P)\neq 0$ such that 
   \[\begin{array}{lll}
\lambda^2 q_j&=&
-\kappa \ds\sum_{k=1}^P 
\ds\int \nabla_x\sigma_1(y-q^*_k)\otimes 
\nabla_x\sigma_1(y-q^*_j)
\ud y\ q_j
\\&&+\kappa_{\lambda^2/c^2}\ds\sum_{k=1}^P\ds\int_{\mathbb R^d}
\nabla_x\sigma_1(y-q_j^*)\otimes \nabla_x \sigma_1(y-q^*_k)\ud y\ q_k.
\end{array}\]
Given $j\in \{1,...,P\}$, there exists a unique $\ell_j\in \{1,...,L\}$ such that $q_j^*\in \Pi_{\ell_j}$, and we denote $r_j=\#\Pi_{\ell_j}$.
Since, when $k\in \Pi_{\ell'}$, $\ell'\neq \ell_j$, $\nabla_x\sigma_1(y-q_j^*)\otimes \nabla_x\sigma_1(y-q^*_k)$ vanishes for all $y\in \R^d$ and, when
$k\in \Pi_{\ell_j}$, $\int \nabla_x\sigma_1(y-q_j^*)\otimes \nabla_x\sigma_1(y-q^*_k)\ud y=\alpha \mathbb I_d$, we get
\begin{equation}\label{eigen_q}
\lambda^2 q_j= r_j\alpha \left(
-\kappa   q_j +\kappa_{\lambda^2/c^2}\langle q \rangle_j\right), \text{ with } \langle q \rangle_j=\frac{1}{r_j}\ds\sum_{k\in \Pi_{\ell_j}}q_k
.
\end{equation}
The same relation holds for any $k\in \Pi_{\ell_j}$
\[
\lambda^2 q_k= \alpha r_j \left(
-\kappa   q_k +\kappa_{\lambda^2/c^2}\langle q \rangle_j\right)
.\]
It follows that 
\begin{equation}\label{mean}
\lambda^2 \langle q \rangle_j=
-\alpha  r_j(\kappa  -\kappa_{\lambda^2/c^2})\langle q \rangle_j,\end{equation}
where we bear in mind that $\langle q \rangle_k$ takes the same value for any $k\in \Pi_{\ell_j}$.
From now on, we will thus denote them as $\langle q \rangle_1,...,\langle q \rangle_L$.
Let us now distinguish several cases:
\begin{itemize}
\item If $\lambda =0$. Since $\kappa_0=\kappa$, we obtain  from \eqref{eigen_q}, for any $\ell\in \{1,...,L\}$ 
and any $j\in \Pi_\ell$,
$q_j=\langle q \rangle_\ell$. 
We conclude that 
$\mathrm{Ker}(\tilde L)$
is made of vectors $X$ where $q$ is blockwise constant according to  the $\Pi_\ell$'s, 
$p=0$, $\varpi=0$ and $\psi(x,z)=\Gamma(z)\sum_{k=1}^P \nabla_x\sigma_1(x-q^*_k)\cdot q_k$.
\item If $\lambda=i\mu$, with $\mu \in \mathbb R\setminus\{0\}$.
The kernel $\xi\mapsto \frac{\widehat \sigma_2(\xi)}{-\mu^2/c^2+\xi^2}$
presents a singularity which prevents the solution $\psi$ for being in $L^2(\mathbb R^n\times\mathbb R^d)$, unless $\psi$ is identically $0$, \emph{i.e.} 
\begin{equation*}
  \sum_{k=1}^P\nabla_x\sigma_1(x-q^*_k)\cdot q_k=0.
\end{equation*}
(This is where we use the fact that  $\widehat \sigma_2$ never vanishes.) 
For any $\ell \in \{1,\ldots,L\}$, let $j\in \Pi_\ell$. Using the relation above, we get
\begin{align*}
  0=&\sum_{k=1}^P\int_{\R^d}\nabla_x\sigma_1(y-q^*_j)\nabla_x\sigma_1(y-q^*_k)\cdot q_k\ud y\\
  &=\sum_{k=1}^P\int_{\R^d}\nabla_x\sigma_1(y-q^*_j)\otimes\nabla_x\sigma_1(y-q^*_k)\ud y q_k= \alpha\sum_{k\in \Pi_\ell}q_k=\alpha \#\Pi_\ell \langle q \rangle_\ell
\end{align*}
so that  $\langle q \rangle_\ell=0$ for all $\ell\in \{1,\ldots, L\}$. 

If $\ell \in \{1,\ldots,L\}$ is such that $\#\Pi_\ell=1$, then $q_\ell=\langle q \rangle_\ell=0$. In particular, if $\#\Pi_\ell=1$ for all $\ell \in \{1,\ldots,L\}$, $\lambda=i\mu$ with $\mu \in \mathbb R\setminus\{0\}$ cannot be an eigenvalue of $\mathbb L$. We will see in a while that $i\mathbb R=\sigma_{\mathrm{ess}}(\tilde L)$. 

If $\ell \in \{1,\ldots,L\}$ is such that $\#\Pi_\ell>1$, the equation for the equilibrium positions \eqref{eigen_q} becomes 
$-\mu^2 q_j=-\kappa \alpha\#\Pi_\ell q_j$, for any $j\in \Pi_\ell$ and $\lambda=\pm i \sqrt{\kappa \alpha\#\Pi_\ell}$ is an eigenvalue of $\mathbb L$. The associated eigenvectors have the form $(q,p,0,0)$, 
where $p_j=\lambda q_j$ for all $j\in \{1,\ldots,P\}$, and the $q_j$'s are  defined such that 
$q_j=0$ if $\#\Pi_j=1$ and $\sum_{k\in\Pi_j} q_k=0$ if $\#\Pi_j>1$.

\item If $\lambda\in \mathbb R\setminus \{0\}$.
In such a case, we have $\frac{\lambda^2}{c^2}>0$ and thus $0<\kappa_{\lambda^2/c^2}<\kappa$. Hence, $\lambda^2+\alpha\#\Pi_\ell(\kappa  -\kappa_{\lambda^2/c^2})>0$, \eqref{mean} leads to $\langle q \rangle_\ell=0$ for any $\ell\in \{1,...,L\}$. Using \eqref{eigen_q}, we get 
\begin{equation*}
  (\lambda^2+\kappa\alpha\#\Pi_\ell)q_j=0
\end{equation*}
for all $j\in \Pi_\ell$, which implies $q_j=0$. As a conclusion, $\lambda$ cannot be an eigenvalue of $\mathbb L$.

 \item If $\lambda=a+ib$ with $a\neq 0$, $b\neq 0$.
 We have $\lambda^2=a^2-b^2+2iab$, where $ab\neq 0$.
 Going back to \eqref{mean}, we get 
 \begin{equation*}
  (\lambda^2+\alpha \#\Pi_\ell(\kappa-\kappa_{\lambda^2/c^2}) )|\langle q \rangle_\ell |^2=0
 \end{equation*}
 for all $\ell\in \{1,\ldots,L\}$. If $\langle q \rangle_\ell\neq 0$, then $ |\langle q \rangle_\ell |^2>0$ and, by
 identifying the imaginary part in the above relation, we are led to 
 \[
 \frac{1}{c^2}=-\ds\alpha\#\Pi_\ell \ds\int_{\mathbb R^n} \ds\frac{|\widehat \sigma_2(\xi)|^2}
 {((a^2-b^2)+c^2\xi^2)^2+4a^2b^2}\ds\frac{\ud\xi}{(2\pi)^n}
 \ds,
\]
 a contradiction, since both side of the equality have opposite sign. Hence, $\langle q \rangle_\ell= 0$ for all $\ell\in \{1,\ldots,L\}$ and, from \eqref{eigen_q}, we deduce 
 \begin{equation*}
  (\lambda^2+\kappa\alpha\#\Pi_\ell)q_j=0
\end{equation*}
for all $j\in \Pi_\ell$. As before, this leads to 
\begin{equation*}
  0=(\lambda^2+\kappa\alpha\#\Pi_\ell)|q_j|^2=(a^2-b^2+\kappa\alpha\#\Pi_\ell)|q_j|^2+2iab |q_j|^2,
\end{equation*} 
which implies $|q_j|^2=0$ and $q_j=0$ for all $j\in \{1,\ldots,P\}$. As a conclusion, $\lambda$ cannot be an eigenvalue for $\mathbb L$.
 \end{itemize}
 
 We turn to the identification of the essential spectrum of $\tilde L$. We start by showing $i\mathbb R\setminus\{0\}\subset \sigma_{\mathrm {ess}}(\tilde L)$.
  To this end, we adapt the construction of  a Weyl sequence for $(-c^2\Delta -\mu^2)$, $\mu\in \R$.
 For any $\mu \in \mathbb R\setminus\{0\}$, $m\in \N^*$, let $Y_m=(0,0,\varphi_m,\varpi_m)$ with $\varphi_m=\tfrac{1}{i\mu}(-\Delta)^{1/2} \varpi_m$. The construction is completed by setting, 
 given $\zeta\in C^\infty_c(\mathbb R^d)$, 
 \begin{equation*}
  \varpi_m(x,z)=\zeta(x)\int_{\R^n}\widehat{\varpi}_m(\xi)e^{i\xi\cdot z}\frac{\ud \xi}{(2\pi)^n},\quad \widehat{\varpi}_m(\xi)=(2\pi)^n m^{n/2}e^{-m^2|\xi-\xi_0|^2}\text{ with } c^2|\xi_0|^2=\mu^2.
\end{equation*}
It is clear that $\varpi_m\in \mathscr S(\mathbb R^d\times \mathbb R^n)$ and that
 $\widehat{\varphi}_m(x,\xi)=\tfrac{\zeta(x)}{i\mu}|\xi| \widehat{\varpi}_m(\xi)$
 also defines an element of 
 $ \mathscr S(\mathbb R^d\times\mathbb R^n)$. 
 We check that
\begin{equation*}
  \|Y_m\|_{\tilde{\mathbb  H}}^2=\frac{1}{(2\pi)^n}\left(\|\widehat{\varphi}_m\|^2_{L^2}+\|\widehat{\varpi}_m\|^2_{L^2}\right)\ge 2^{n/2}\pi^{3n/2}
  \|\zeta\|_{L^2(\mathbb R^d)}^2
  >0
\end{equation*}
for any $m\in \N^*$.
 We have
\[
(\tilde L-i\mu)Y_m=\begin{pmatrix}0\\
P'_m
\\
0
\\
\frac{i}{\mu}(-c^2\Delta -\mu^2)\varpi_m
\end{pmatrix},\] 
where the components of the vector $P'_m$ are given by
\[
\frac{1}{i \mu}\ds\iint_{\mathbb R^d\times\mathbb R^n} \sigma_2(z) \nabla_x\sigma_1(y-q^*_j) \varpi_m(y,z)\ud z\ud y.\]  
By construction, we have
\begin{align*}
  \|(-c^2\Delta -\mu^2)\varpi_m\|^2_{L^2}&=
  (2\pi m)^n \|\zeta\|_{L^2}^2\int_{\R^n}(c^2|\xi|^2 -\mu^2)^2e^{-2m^2|\xi-\xi_0|^2}\ud \xi\\
  &=(\sqrt 2 \pi)^n \|\zeta\|_{L^2}^2\int_{\R^n}\left(c^2\left|\frac{\tilde \xi}{\sqrt 2 m}+\xi_0\right|^2 -\mu^2\right)^2e^{-|\tilde \xi|^2}\ud \tilde \xi\\
  &=c^4(\sqrt 2 \pi)^n \|\zeta\|^2_{L^2} \int_{\R^n}\left(\frac{|\tilde \xi|^2}{2 m^2}+\frac{\sqrt 2}{m}\tilde \xi\cdot \xi_0\right)^2 e^{-|\tilde k|^2}\ud \tilde k
  \xrightarrow [m\to \infty]{}0.
\end{align*}
 Eventually, we establish the weak convergence to 0 of $X_m$: for any $f\in \mathscr S(\mathbb R^d\times \mathbb R^n)$, we have 
\[\begin{array}{lll}
  \langle \varpi_m,f\rangle&=&\ds\frac{1}{(2\pi)^n}\langle \widehat{\varpi}_m,\hat f\rangle=m^{n/2}\iint_{\mathbb R^d\times \mathbb R^n}e^{-m^2|\xi-\xi_0|^2}\hat f(x,\xi)\zeta(x)\ud \xi\ud x
  \\&=&\left(\ds\frac{1}{m}\right)^{n/2}\ds\iint_{\mathbb R^d\times \mathbb R^n}e^{-|\tilde \xi|^2}\zeta(x)\hat f\left(x,\frac{\tilde \xi}{m}+\xi_0\right)\ud \tilde \xi\ud x
\end{array}\]
so that 
  $|\langle \varpi_m,f\rangle|\lesssim  \left(\frac{1}{m}\right)^{n/2}$. Similarly, we show that 
$|\langle \varphi_m,f\rangle|\lesssim \left(\frac{1}{m}\right)^{n/2}$. This implies that $X_m$ converges weakly to $0$ as $m\to +\infty$.
Furthermore, since $\varpi_m \rightharpoonup 0$ and $\|(-c^2\Delta -\mu^2)\varpi_m\|^2_{L^2}\to 0$ as $m\to +\infty$, we deduce that $\|(\tilde L-i\mu)Y_m\|_{\tilde{\mathbb  H}}\to 0$ as $m\to +\infty$. We conclude that $(Y_m)_{m\in  \N^*}$ is a Weyl sequence and $i\mu\in \sigma_{\mathrm{ess}}(\mathbb L)$ for any $\mu\in \R\setminus\{0\}$.

We treat 
with a similar idea the case $\lambda=0$:
we construct an appropriate sequence under the form  $Y_m=(0,0,\varphi_m,0)$ with $
\varphi_m(x,z)=
\zeta(x)(-\Delta)^{1/2}(m^{1-n/2}\chi(z/m))$, where $\zeta\in C^\infty_c(\mathbb R^d)$, $\chi\in C^\infty_c(\mathbb R^n)$.
In particular, $\|\varphi_m\|_{L^2}=\|\zeta\|_{L^2(\mathbb R^d)} \|\nabla \chi\|_{L^2(\mathbb R^n)}$ is constant.
We have
\[\tilde LY_m=
\begin{pmatrix}
0 \\
P'_m\\0
\\
c^2\zeta(x)\Delta(m^{1-n/2}\chi(z/m))
\end{pmatrix}
\]
where the components of $P'_m$ read
\[\ds\iint_{\mathbb R^d\times\mathbb R^n} \sigma_2(z)\nabla_x\sigma_1(y-q^*_j)
\zeta(y)m^{1-n/2}\chi(z/m)
 \ud z\ud y.\]
 As a matter of fact, the latter can be dominated by $$m^{1-n/2}\|\sigma_2\|_{L^1(\mathbb R^n)} \|\zeta\|_{L^1(\mathbb R^d)}
 \|\nabla\sigma_1\|_{L^\infty(\mathbb R^d)}\|\chi\|_{L^\infty(\mathbb R^n)}$$
 and it thus tends to 0 as $m\to \infty$, since $n\geq 3$.
 The last component of $\tilde LY_m$ recasts as 
 $c^2\zeta(x)m^{-n/2-1}\Delta\chi(z/m)$: its $L^2$ norm is 
 $\frac{c^2}{m}\|\zeta\|_{L^2(\mathbb R^d)} \|\Delta \chi\|_{L^2(\mathbb R^n)}$.
 We conclude that $\tilde LY_m$  converges to 0 in $\tilde{\mathbb H}$.
 Moreover, for $f\in C^\infty_c(\mathbb R^d\times \mathbb R^n$),  we have
 $$\big|\iint_{\mathbb R^d\times\mathbb R^n} f(x,z)\varphi_m(x,z)\ud z\ud x\big|
 \leq m^{1-n/2} \|\chi\|_{L^\infty(\mathbb R^n)} \|\zeta\|_{L^2(\mathbb R^d)}\|\nabla f\|_{L^2(\mathbb R^n\times\mathbb R^d)},$$ 
and by density we deduce that $\varphi_m\rightharpoonup 0$ in $L^2(\mathbb R^n\times\mathbb R^d)$.
Hence, $\big(Y_m\big)_{m\in \mathbb N}$ is a Weyl sequence for $\tilde L$ and $0\in \sigma_{\mathrm{ess}}(\tilde L)$.

Conversely, let $\lambda \in \C\setminus i\R$ (in particular $\lambda \neq 0$) and suppose, by contradiction, that there exists a Weyl sequence $(Y_m)_{m\in \mathbb N}$: $Y_m\in \tilde{\mathbb H}$ with  $\|Y_m\|_{\tilde{\mathbb H}}=1$, 
and $Y_m\rightharpoonup 0$ while
\begin{equation*}
  (\tilde  L-\lambda) Y_m=Y'_m\xrightarrow[m\to \infty]{} 0 \text{ strongly in } \tilde{\mathbb H}.
\end{equation*}
From $(-\Delta)^{1/2} \varpi_m -\lambda \varphi_m=\varphi'_m$, we deduce 
\begin{equation*}
  (\lambda^2-c^2\Delta)\varpi_m= c^2(-\Delta)^{1/2}\varphi'_m+
 \lambda c^2
   \sigma_2\ds\sum_{k=1}^P\nabla_x\sigma_1(x-q_k^*)\cdot q_{m,k}  -\lambda\varpi'_m
\end{equation*}
and, since $\lambda^2+c^2|\xi|^2\neq 0$ for all $\xi\in \R^n$, we can write
\begin{equation}\label{def_phi_essspectrum}
  \widehat{\varpi}_m(x,\xi)= \lambda c^2 \frac{  \hat \sigma_2(\xi) }{\lambda^2+c^2|\xi|^2}
\ds\sum_{k=1}^P\nabla_x\sigma_1(x-q_k^*)\cdot q_{m,k}+ \underbrace{\frac{c^2|\xi|\widehat{\varphi'}_m(x,\xi) -\lambda\widehat{\varpi'}_m(x,\xi)}{\lambda^2+c^2|\xi|^2}}_{:=h_m(x,\xi)},
\end{equation}
where we check that
\begin{align*}
  \iint_{\R^d\times \R^n}|h_m(x,\xi)|^2\ud \xi\ud x \lesssim  \iint_{\R^d\times \R^n}
  |\widehat{\varphi'}_m(x,\xi)|^2\ud \xi\ud x+ \iint_{\R^d\times \R^n}|\widehat{\varpi'}_m(x,\xi)|^2\ud \xi\ud x \xrightarrow[m\to +\infty]{}0.
\end{align*}
 Similarly, for all $\xi\in \R^n$,
\begin{equation}\label{def_pi_essspectrum}
\begin{array}{lll}
  \widehat{\varphi}_m(x,\xi)&=&\ds\frac{1}{\lambda}|\xi|\widehat{\varpi}_m(x,\xi)-\frac{1}{\lambda}\widehat{\varphi'}_m(x,\xi)
  \\&=&c^2\ds\frac{ |\xi| \hat \sigma_2(\xi) }{\lambda^2+c^2|\xi|^2}
  \ds\sum_{k=1}^P\nabla_x\sigma_1(x-q_k^*)\cdot q_{m,k}
  +\ds\frac1\lambda \underbrace{\frac{-\lambda^2\widehat{\varphi'}_m(x,\xi) -\lambda|\xi|\widehat{\varpi'}_m(x,\xi)}{\lambda^2+c^2|\xi|^2}}_{:=g_m(x,\xi)}
\end{array}\end{equation} 
where $g_m$ converges strongly to $0$ in $L^2(\R^d\times \R^n)$. 
This yields the following relations for the components    $(q_{m,j},p_{m,j})$ 
\begin{equation}\label{sys_pq_essspectrum}
  \left\{
    \begin{aligned}
      &\lambda q_{m,j}= p_{m,j}-q'_{m,j}
      \\
      &\lambda p_{m,j}
 =-\alpha\kappa r_j q_{m,j} 
 +\ds\iint_{R^d\times\R^n} \ds\frac{c^2|\hat\sigma_2(\xi)|^2}{\lambda^2+c^2\xi^2}\ds\sum_{k=1}^P \nabla_x\sigma_1(y-q^*_j)\otimes \nabla_x\sigma_1(y-q^*_k)\cdot q_k\ds\frac{\ud \xi}{(2\pi)^n}\ud y     
 \\
 &\hspace*{1.5cm}-p'_{m,j}+\ds\frac1\lambda \ds\iint_{\R^d\times\R^n} \ds\frac{\hat\sigma_2(\xi)}{|\xi|}g_m(y,\xi)
 \nabla_x\sigma_1(y-q^*_j)\ds\frac{\ud \xi}{(2\pi)^n}\ud y
 \\&\hspace*{1cm}
 =-\alpha\kappa r_j q_{m,j}  +\alpha \kappa_{\lambda^2/c^2} r_j\langle q_m\rangle_j
 -
 p'_{m,j}\\     
 &\hspace*{1.5cm}+\ds\frac1\lambda \ds\iint_{\R^d\times\R^n} \ds\frac{\hat\sigma_2(\xi)}{|\xi|}g_m(y,\xi)
 \nabla_x\sigma_1(y-q^*_j)\ds\frac{\ud \xi}{(2\pi)^n}\ud y.
    \end{aligned}
  \right.
\end{equation}
which leads to 
\begin{align*}
\lambda^2 q_{m,j}=&\alpha r_j\big(-\kappa q_{m,j}+\kappa_{\lambda^2/c^2}\langle q_m\rangle_j\big)
-\lambda q'_{m,j} -p'_{m,j} \\
  &+\ds\frac1\lambda \ds\iint_{\R^d\times\R^n} \ds\frac{\hat\sigma_2(\xi)}{|\xi|}g_m(y,\xi)
 \nabla_x\sigma_1(y-q^*_j)\ds\frac{\ud \xi}{(2\pi)^n}\ud y.
\end{align*}
It follows that 
\begin{align*}
\big(\lambda^2+\alpha r_j(\kappa-\kappa_{\lambda^2/c^2})\big)
\langle q_m\rangle_j
=&-\lambda \langle q'_{m}\rangle_{j} -\langle p'_{m}\rangle_{j}\\ 
&+\ds\frac\lambda{r_j}\ds\sum_{k\in\Pi_j} \ds\iint_{\R^d\times\R^n} \ds\frac{\hat\sigma_2(\xi)}{|\xi|}g_m(y,\xi)
 \nabla_x\sigma_1(y-q^*_k)\ds\frac{\ud \xi}{(2\pi)^n}\ud y.
\end{align*}
The right hand side of this expression converges to $0$ as $m\to +\infty$. 
Moreover, we have seen above that 
$\lambda^2+\alpha r_j(\kappa-\kappa_{\lambda^2/c^2})$ does not vanish when $\lambda\in \mathbb R\setminus \{0\}$ and when $\lambda=a+ib$ with $ab\neq 0$.
As a consequence, we deduce, firstly, that $\langle q_m\rangle_j\to 0$ 
and, secondly, that $q_{m,j}\to 0$
as $m\to +\infty$.
Therefore, coming back to \eqref{def_phi_essspectrum}, \eqref{def_pi_essspectrum} and \eqref{sys_pq_essspectrum}, we conclude that  $X_m$ converges strongly to $0$ in $\tilde{\mathbb H}$ which  contradicts $\|X_m\|_{\tilde{\mathbb H}}=1$. Hence, we get
 $\sigma_{\mathrm{ess}}(\tilde L)=i\R$. 
\QED

\begin{rmk}
  The existence of  eigenvalues in $i\mathbb R$ is a purely ``many particles'' effect: it does not hold if all $\Pi_\ell$ contain only a single particle.
\end{rmk}

It would be natural to discuss the  stability of the equilibrium solutions based on the properties of the energy functional $\mathscr H$. This is a well established approach; however, in the present context, 
the translational invariance of $\mathscr H$ leads to a lack of coercivity 
and the final result would have a quite poor expression, with the impossibility of locating the asymptotic clusters. This difficulty disappears   when an additional confining potential is introduced in the equations: this is further discussed in Appendix~\ref{App}.

\section{Large time behavior of the linearized problem}
\label{Tlarge_Lin}

 The behavior of the solutions of the linearized problem can be made more precise,  by using a formulation inspired from  the non autonomous expression presented in Section~\ref{prelim}.
 By solving the wave equation in \eqref{Lin}, we are led to 
 \begin{equation}\label{Lin_nonaut}
  \left\{
    \begin{aligned}
      \dot{q}_j^{\mathrm{lin}}(t)=&p_j^{\mathrm{lin}}(t),
\\
\dot{p}^{\mathrm{lin}}_j(t)=&-\ds\iint_{\mathbb R^d\times\mathbb R^n}
\ds\frac{|\widehat \sigma_2(\xi)|^2}{|\xi|^2}
\nabla_x\sigma_1(y-q_j^*)\cdot q_j^{\mathrm{lin}}(t)
\ds\sum_{k=1}^P \nabla_x \sigma_1(y-q_k^*)
\ds\frac{\ud \xi}{(2\pi)^n}\ud y
\\&+
\ds\int_{\mathbb R^d}
\nabla_x\sigma_1(y-q_j^*)
\ds\int_0^t \mathpzc p_c(t-s) \ds\sum_{k=1}^P \nabla_x \sigma_1(y-q_k^*)\cdot q_k^{\mathrm{lin}}(s)
\ud s
\ud y
\\&+\ds\iint_{\mathbb R^d\times\mathbb R^n}
\nabla_x\sigma_1(y-q_j^*)
\widehat \sigma_2(\xi) 
\Big(\cos(c|\xi|t) \widehat\Psi_0^{\mathrm{lin}}(y,\xi)
+\ds\frac{\sin(c|\xi|t)}{c|\xi|}  \widehat\Psi_1^{\mathrm{lin}}(y,\xi)\Big)\ds\frac{\ud \xi}{(2\pi)^n}\ud y.
\end{aligned}
  \right.
\end{equation}
The last integral does not depend on the solution $q_j^{\mathrm{lin}}(t)$, and it will be considered as a given source term, that we denote $S_{I,j}(t)$, namely
\begin{align}
  \label{def_S_Ij}
  S_{I,j}&(t)=\nonumber\\
  &\frac{1}{(2\pi)^n}\ds\int_{\mathbb R^d}\nabla_x\sigma_1(y-q^*_j)
  \left(\ds \int_{\mathbb R^n} \left( \cos(c|\xi|t) \widehat{\Psi_0^{\mathrm{lin}}}(y,\xi)+  \ds\frac{\sin(c|\xi|t)}{c|\xi|} \widehat{\Psi_1^{\mathrm{lin}}}(y,\xi)
  \right)\ 
  \widehat{\sigma_2}(\xi) \ud\xi\right)\ud y
\end{align}
for all $t\in \R^+$.

\begin{lemma}\label{lem:SI}
Let $\Psi_0^{\mathrm{lin}}\in L^2(\mathbb R^d\times\mathbb R^n)$ and $\Psi_1^{\mathrm{lin}}\in L^2(\mathbb R^d;\dot H^1(\mathbb R^n))$. Then, $$S_I=(S_{I,1},\ldots, S_{I,P})\in L^{\infty}(\R^+;\R^P).$$
If, moreover $n\ge 3$ is odd and $\Psi_0^{\mathrm{lin}}$ and $\Psi_1^{\mathrm{lin}}$ are supported in $\mathbb R^d\times B(0,R_0)$, then $S_I$ is compactly supported.
\end{lemma}
 
\noindent
{\bf Proof.}
A direct inspection of  \eqref{def_S_Ij} shows that $S_{I,j}\in L^\infty(\R^+)$ for any $j\in \{1,\ldots,P\}$.
Next, we remark that for any $t\in \R^+$, $S_{I,j}$ can be recast as
\[
S_{I,j}(t)=\ds\int_{\mathbb R^d}\nabla_x\sigma_1(x-q^*_j)\left(\ds \int_{\mathbb R^n} \sigma_2(z) \Upsilon(t,x,z)\ud z\right)\ud x,\]
where $\Upsilon$ is solution of the Cauchy problem for the wave equation
$\partial^2_{tt} \Upsilon-c^2\Delta_z \Upsilon=0$, with $( \Upsilon,\partial_t \Upsilon)|_{t=0}=(\Psi_0^{\mathrm{lin}},\Psi_1^{\mathrm{lin}})$.
Since $z\mapsto  \Psi_0^{\mathrm{lin}}(x,z)$ and $z\mapsto  \Psi_1^{\mathrm{lin}}(x,z)$
are both compactly supported in $B(0,R_0)$, and $n\ge 3$ is odd, 
Huygens' principle implies that $z\mapsto \Upsilon(t,x,z)$ is supported in the shell $\{z\in \mathbb R^n,\ ct-R_0\leq |z|\leq ct+R_0\}$.
As a consequence, $\Upsilon(t,x,z)=0$ whenever
$|z|<ct-R_0$. Since $\sigma_2$ is supported in $B(0,R_2)$, it follows that the product $\sigma_2(z) \Upsilon(t,x,z)$ vanishes when $t\geq \frac{R_0+R_2}{c}$. In turn, $S_{I,j}$ is supported in $[0, \frac{R_0+R_2}{c}]$ for all $j\in \{1,\ldots,P\}$.
\QED

 The first  two terms in the right hand side of the momentum equation in \eqref{Lin_nonaut} can be cast as
\begin{equation}\label{triangle} \begin{array}{l}
 -\kappa \ds\sum_{k=1}^P
 \ds\int_{\mathbb R^d}
\nabla_x\sigma_1(y-q_k^*)\otimes 
 \nabla_x \sigma_1(y-q_j^*)
\ud y  \ q_j^{\mathrm{lin}}(t)
\\\qquad +
\ds\sum_{k=1}^P
 \ds\int_{\mathbb R^d}
\nabla_x\sigma_1(y-q_j^*)\otimes 
  \nabla_x \sigma_1^*(y-q_k^*)
  \ud y \ds\int_0^t \mathpzc p_c(t-s) q_k^{\mathrm{lin}}(s)
\ud s
\\
 =
-\alpha r_j \left( \kappa q_j^{\mathrm{lin}}(t)- \ds\int_0^t \mathpzc p_c(t-s) \langle q^{\mathrm{lin}}\rangle_j (s)\ud s\right)
. \end{array}\end{equation}
 Owing to Lemma~\ref{lemAV},
it becomes
\[
-\alpha r_j \ds\int_0^t \mathpzc p_c(t-s) \big( q_j^{\mathrm{lin}}(t)-  \langle q^{\mathrm{lin}}\rangle_j (s)\big)\ud s
-\alpha r_j \ds\int_t^\infty \mathpzc p_c(s)\ud s\ q_j^{\mathrm{lin}}(t).
\]
Still by virtue of Lemma~\ref{lemAV}, the last term vanishes when $t\geq \frac{2R_2}{c}$.
Therefore, we obtain 
 \[
 \dot p^{\mathrm{lin}}_j(t)=-\alpha r_j \ds\int_0^t \mathpzc p_c(t-s) \big( q_j^{\mathrm{lin}}(t)-  \langle q^{\mathrm{lin}}\rangle_j (s)\big)\ud s + R_j(t)
 \]
where $$R_j(t)=
-\alpha r_j \ds\int_t^\infty \mathpzc p_c(s)\ud s\ q_j^{\mathrm{lin}}(t)+
S_{I,j}(t)$$ (which now depends on $q^{\mathrm{lin}}_j$) vanishes when $t$ is large enough.
More generally, for any $k\in \Pi_j$, we get
 \[
 \dot p^{\mathrm{lin}}_k(t)=-\alpha r_j \ds\int_0^t \mathpzc p_c(t-s) \big( q_k^{\mathrm{lin}}(t)-  \langle q^{\mathrm{lin}}\rangle_j (s)\big)\ud s + R_k(t)
 .\]
Summing over $\Pi_j$, we obtain
\[
  \langle \dot p^{\mathrm{lin}}\rangle_j=-\alpha r_j \ds\int_0^t \mathpzc p_c(t-s) \big( \langle q^{\mathrm{lin}}\rangle _j(t)-  \langle q^{\mathrm{lin}}\rangle_j (s)\big)\ud s + \langle R\rangle _j(t)
 .\]
We now adapt the arguments in \cite{BdB}, expanding
 \[\langle q^{\mathrm{lin}}\rangle_j(t)-\langle q^{\mathrm{lin}}\rangle_j(s)=
 \ds\int_s^t \langle p^{\mathrm{lin}}\rangle_j(\sigma)\ud \sigma
 =
 (t-s)\langle p^{\mathrm{lin}}\rangle_j(t)- \ds\int_s^t (\sigma-s)\langle \dot p^{\mathrm{lin}}\rangle_j(\sigma)\ud \sigma.
 \]
 We arrive at 
 \begin{align*}
\langle \dot p ^{\mathrm{lin}}\rangle _j(t)&=-\alpha r_j \ds\int_0^t (t-s) \mathpzc  p_c(t-s)\ud s\ 
 \langle p^{\mathrm{lin}}\rangle_j(t)
 \\&\quad+ \alpha r_j \ds\int_0^t \mathpzc p_c(t-s) \left(\ds\int_s^t (\sigma-s)\langle \dot p^{\mathrm{lin}}\rangle_j(\sigma)\ud \sigma\right)\ud s + \langle R\rangle _j(t).\end{align*}
 With $n=3$, we appeal to Lemma~\ref{mtpc} which yields
\begin{align*}
\langle \dot p ^{\mathrm{lin}}\rangle _j(t)&=-\ds\frac{\upgamma \alpha r_j }c
 \langle p^{\mathrm{lin}}\rangle_j(t)
+ \alpha r_j \ds\int_0^t \mathpzc p_c(t-s) \left(\ds\int_s^t (\sigma-s)\langle \dot p^{\mathrm{lin}}\rangle_j(\sigma)\ud \sigma\right)\ud s \\
&\quad  +\alpha r_j \ds\int_t^\infty s \mathpzc  p_c(s)\ud s\ 
 \langle p^{\mathrm{lin}}\rangle_j(t)+
 \langle R\rangle _j(t).\end{align*}
 where the last two terms vanish when $t$ becomes large.
 Then, it remains to investigate the remainder term in order to prove
 
 \begin{proposition}\label{prop:meanLin} Let $n=3$. Let $\upgamma$ be defined in Lemma~\ref{mtpc}. Provided the wave speed $c$ is large enough, the mean momentum $\langle p^{\mathrm{lin}}\rangle_j$ for any $j\in \{1,...,L\}$ tends to 0 exponentially fast, and the mean position $\langle q^{\mathrm{lin}}\rangle_j$ tends to a constant. Precisely, for any $\epsilon\in (0,1)$, there exists $c_*(\epsilon)>0$ such that if $c\geq c_*(\epsilon)$, then 
 for any initial data $q^{\mathrm{lin}}_0,p_0^{\mathrm{lin}},\Psi^{\mathrm{lin}}_0,\Psi^{\mathrm{lin}}_1$, with 
 $\mathrm{supp}(\Psi^{\mathrm{lin}}_0,\Psi^{\mathrm{lin}}_1)\subset \mathbb R^d\times B(0,R_0)$ for some $0<R_0<\infty$,
we can find $K(\epsilon)$ and $q_{j,\infty}$ such that 
\[|\langle p^{\mathrm{lin}}\rangle_j(t)|\leq K(\epsilon)e^{-\upgamma \alpha r_j(1-\epsilon) t/c},\qquad
|\langle q^{\mathrm{lin}}\rangle_j(t)-q_{j,\infty} |\leq \frac{K(\epsilon)c}{\upgamma \alpha r_j(1-\epsilon)}e^{-\upgamma \alpha r_j(1-\epsilon) t/c}.\]
\end{proposition}

\noindent
{\bf Proof.}
 We shall need to keep track of the dependence 
 of certain coefficients with respect to the wave speed $c$. 
 In particular, it will be useful to 
 remark that, owing to the change of variable $s=ct$, 
 $ \|\mathpzc{p}_c\|_{L^1((0,\infty))}=\|\mathpzc{p}_1\|_{L^1((0,\infty))}$ does not depend on $c$.
 We shall also make use of the fact that $\mathpzc{p}_c$ and the source $S_{I,j}$ are both  compactly supported:
 let 
 \[ R_*=\max\Big(R_0+R_2,2R_2\Big)\]
 so that $\mathrm{supp}(\mathpzc{p}_c)\subset [0,\frac{R_*}{c}]$,  $\mathrm{supp}(S_{I,j})\subset [0,\frac{R_*}{c}]$.
 According to the previous observations, we have 
 \begin{align*}
&\langle \dot p ^{\mathrm{lin}}\rangle _j(t)=-\ds\frac{\upgamma \alpha r_j }c
 \langle p^{\mathrm{lin}}\rangle_j(t)
+ \alpha r_j \ds\int_0^t \mathpzc p_c(t-s) \left(\ds\int_s^t (\sigma-s)\langle \dot p^{\mathrm{lin}}\rangle_j(\sigma)\ud \sigma\right)\ud s 
+\Xi_j(t),\\& \mathrm{supp}(\Xi_j)\subset [0,R_*/c]. 
\end{align*}

We start with a rough estimate on the solution: from 
\begin{align*}
\langle q^{\mathrm{lin}}\rangle_j(t)=&\langle q^{\mathrm{lin}}\rangle_j(0)+\ds\int_0^t \langle p^{\mathrm{lin}}\rangle_j(s)\ud s,
\\
\langle p^{\mathrm{lin}}\rangle_j(t)=&\langle p^{\mathrm{lin}}\rangle_j(0)\\
&+\ds\int_0^t \left(-\alpha r_j \kappa \langle q^{\mathrm{lin}}\rangle_j(s)+
\alpha r_j \left(\ds\int_0^s \mathpzc{p}_c(s-\sigma) \langle q^{\mathrm{lin}}\rangle_j (\sigma)\ud \sigma\right)
+\langle  S_{I}\rangle_j (s)
\right)\ud s
\end{align*}
(where we used \eqref{triangle}),  we obtain 
\begin{align*}
&|\langle p^{\mathrm{lin}}\rangle_j(t)|\leq |\langle p^{\mathrm{lin}}\rangle_j(0)|+
\alpha r_j \kappa
\ds\int_0^t \left(|\langle q^{\mathrm{lin}}\rangle_j(0)|+\ds\int_0^s| \langle p^{\mathrm{lin}}\rangle_j(\sigma)|\ud \sigma\right)\ud s
\\\qquad &+
\alpha r_j 
\ds\int_0^t \Big(|\langle q^{\mathrm{lin}}\rangle_j(0)|+\ds\int_0^\sigma| \langle p^{\mathrm{lin}}\rangle_j(\tau)
|\ud \tau\Big)
\left(\ds\int_\sigma^t |\mathpzc{p}_c(s-\sigma)|
\ud s\right)
 \ud\sigma
+\ds\int_0^t
|\langle  S_{I}\rangle_j (s)|
\ud s\\
&\leq
 |\langle p^{\mathrm{lin}}\rangle_j(0)|+
2\alpha r_j \| \mathpzc{p}_c\|_{L^1}
|\langle q^{\mathrm{lin}}\rangle_j(0)| t 
+ \|\langle  S_{I}\rangle_j \|_{{L^1((0,t))}}
+2\alpha r_j \| \mathpzc{p}_c\|_{L^1} t\ds\int_0^t 
 |\langle p^{\mathrm{lin}}\rangle_j(s)| 
 \ud s.
\end{align*}
Since $\|\mathpzc{p}_c\|_{L^1} =\|\mathpzc{p}_1\|_{L^1}$, applying the
 Gr\"onwall's lemma leads to 
\begin{align*}
|\langle p^{\mathrm{lin}}\rangle_j(t)|&\leq 
e^{2\alpha r_j \| \mathpzc{p}_1\|_{L^1} t^2}
\Big( |\langle p^{\mathrm{lin}}\rangle_j(0)|+
2\alpha r_j \| \mathpzc{p}_1\|_{L^1}
|\langle q^{\mathrm{lin}}\rangle_j(0)| t 
+ \|\langle  S_{I}\rangle_j \|_{{L^1((0,t))}}\Big).
\end{align*}
Accordingly, we also have
\begin{align*}
|\langle q^{\mathrm{lin}}\rangle_j(t)|&\leq |\langle q^{\mathrm{lin}}\rangle_j(0)|
+ te^{2\alpha r_j \| \mathpzc{p}_1\|_{L^1} t^2}
\Big( |\langle p^{\mathrm{lin}}\rangle_j(0)|+
2\alpha r_j \| \mathpzc{p}_1\|_{L^1}
|\langle q^{\mathrm{lin}}\rangle_j(0)| t 
+ \|\langle  S_{I}\rangle_j \|_{{L^1((0,t))}}\Big).
\end{align*}

We are going now to establish a more precise estimate. Let us introduce a parameter $\theta>0$, that will be determined later on, and set
\[\begin{array}{l}
\tilde p_j(t)=e^{\theta t}\langle p^{\mathrm{lin}}\rangle_j(t),
\\[.4cm]
M_0(t)=\ds\sup_{0\leq s\leq t} | \tilde  p_j(s)|= \ds\sup_{0\leq s\leq t} e^{\theta t} \Big| \langle p^{\mathrm{lin}}\rangle_j(s)\Big|,\\[.4cm]
M_1(t)=\ds\sup_{0\leq s\leq t} \Big|\ds\frac{\ud}{\ud s} \tilde p_j(s)\Big|.
\end{array}\]
Since 
$\frac{\ud}{\ud t}\tilde p_j(t)-\theta \tilde p_j(t)=e^{\theta t}\langle \dot p^{\mathrm{lin}}\rangle_j(t)$
we get
\begin{equation}\label{tildep}
\ds\frac{\ud}{\ud t} \tilde p_j(t)
=
\Big(\theta-\ds\frac{\upgamma\alpha r_j}{c}\Big) \tilde p_j(t)
+  \alpha r_j e^{\theta t}\ds\int_0^t \mathpzc p_c(t-s) \left(\ds\int_s^t (\sigma-s)\langle \dot p^{\mathrm{lin}}\rangle_j(\sigma)\ud \sigma\right)\ud s 
 + e^{\theta t}\Xi_j(t).
\end{equation}
The last  term is supported in $[0,R_*/c]$ and is dominated by 
\begin{align*} 
&e^{\theta R_*/c}\|\Xi_j\|_{L^\infty([0,R_*/c))}\\
&\quad\leq e^{\theta R_*/c}\Big(
\alpha r_j \|\mathpzc{p}_1\|_{L^1} \Big(\ds\frac{R_*}{c}\| \langle p^{\mathrm{lin}}\rangle_j\|_{L^\infty([0,R_*/c])}
+\| \langle q^{\mathrm{lin}}\rangle_j\|_{L^\infty([0,R_*/c])}\Big)
+\|\langle S_{I}\rangle_j\|_{L^\infty([0,R_*/c])}
\Big).
\end{align*}
Owing to the basic estimate obtained above, we obtain
$$ e^{\theta R_*/c}\|\Xi_j\|_{L^\infty([0,R_*/c))}
\leq e^{\theta R_*/c} 
e^{2\alpha r_j\|  \mathpzc p_1\|_{L^1} R_*^2/c^2}
(A_0+B_0(R_*/c)^2):=\mu(c)$$
for constants $A_0,B_0$ that depend only on the initial data.
The penultimate term in \eqref{tildep} can be dominated by
\[\begin{array}{l}
\left|\alpha r_j\ds\int_0^t e^{\theta (t-s)}\mathpzc{p}_c(t-s)\left(
\ds\int_s^t
e^{-\theta (\sigma-s)}(\sigma-s)\Big( \frac{\ud}{\ud t}\tilde p_j(\sigma)-\theta \tilde p_j(\sigma)\Big)
\ud \sigma\right)\ud s
\right|
 \\[.4cm]
 \leq
  \alpha r_j 
 \ds\frac{R^2_* e^{\theta R_*/c}}{2c^2}\|\mathpzc{p}_1\|_{L^1}
 \big(M_1(t)+\theta M_0(t)\big),
\end{array}\]
where we used the fact that $\mathpzc{p}_c$ is supported in $[0,R_*/c]$.
Finally,  by choosing $\theta >\alpha r_j\frac \upgamma c $, 
  we arrive at 
\[\begin{array}{lll}
\Big|\ds\frac{\ud}{\ud t} \tilde p_j(t)\Big|
&\leq& 
\Big(\theta-\alpha r_j\ds\frac\upgamma c\Big)|\tilde p_j(t)|+
 \alpha r_j 
 \ds\frac{R^2_* e^{\theta R_*/c}}{2c^2}\|p_1\|_{L^1}
 \big(M_1(t)+\theta M_0(t)\big)+\mu(c)\mathbf{1}_{[0,R_*/c]}(t).
\end{array}\]
 It yields
\[
\Big(1-  \ds\frac{\alpha r_j R^2_* 
 \|\mathpzc{p}_1\|_{L^1}}{2c^2}e^{\theta R_*/c}\Big)M_1(t)
\leq 
\Big(\theta-\ds\frac{\upgamma\alpha r_j}{c}+\theta 
\ds\frac{\alpha r_j R^2_*
 \|\mathpzc{p}_1\|_{L^1}}{2c^2}e^{\theta R_*/c}
\Big)M_0(t) + \mu(c)\mathbf{1}_{[0,R_*/c]}(t).
\]

Let $\epsilon\in (0,1)$ be fixed once for all.
We set $\theta=(1+\epsilon)\frac{\alpha r_j\upgamma}{c}$.
Then, we can pick $c_*(\epsilon)$ large enough so that
\begin{equation}\label{eq_cond_c_eps}
 \ds\frac{\alpha r_j R^2_*
  \|\mathpzc{p}_1\|_{L^1}}{2c^2}e^{\theta R_*/c}\le \frac{\epsilon}{(1+\epsilon)}
\end{equation}
 holds for any $c\geq c_*(\epsilon)$. It follows that
 $$
 0<\frac12<1-  \frac{\alpha r_j R^2_*
 \|\mathpzc{p}_1\|_{L^1}}{2c^2}e^{\theta R_*/c}<1 
 $$
 and
\[\begin{array}{l}\Big(1- \ds \frac{\alpha r_j R^2_* 
 \|p_1\|_{L^1}}{2c^2}e^{\theta R_*/c}\Big)^{-1}
\Big(\theta-\ds\frac{\upgamma\alpha r_j}{c}+\theta 
\ds\frac{\alpha r_j R^2_* 
 \|\mathpzc{p}_1\|_{L^1}}{2c^2}e^{\theta R_*/c}
\Big)
\\[.4cm]\qquad=
\Big(1- \ds \frac{\alpha r_j R^2_*
\|\mathpzc{p}_1\|_{L^1}}{2c^2}e^{\theta R_*/c}\Big)^{-1}
\ds\frac{\upgamma\alpha r_j}{c}
\Big(\epsilon+(1+\epsilon)
\ds\frac{\alpha r_j R^2_*
 \|\mathpzc{p}_1\|_{L^1}}{2c^2}e^{\theta R_*/c}\Big)
\leq 
4\epsilon\ds\frac{\upgamma\alpha r_j}{c}
\end{array}\]
 for any $c\geq c_*(\epsilon)$.
Since 
\[
|\tilde p_j(t)|=\left|\tilde p_j(0)+\int_0^t \ds\frac{\ud}{\ud s}\tilde p_j(s)\ud s\right|
\leq |\tilde p_j(0)|+\ds\int_0^t M_1(s) \ud s
\]
we are  led to 
\[
M_0(t)\le |\tilde p_j(0)|+\ds\int_0^t M_1(s)\ud s
\leq 
|\langle p^{\mathrm{lin}}\rangle_j(0)|+
\ds\int_0^t \Big(4\epsilon\ds\frac{\upgamma\alpha r_j}{c}M_0(s)+4\mu(c)\mathbf{1}_{[0,R_*/c]}(s)\Big)\ud s.
\]
Gr\"onwall's lemma implies 
\[
e^{(1+\epsilon)\upgamma\alpha r_j  t/c }|\langle p^{\mathrm{lin}}\rangle_j(t)| 
\leq M_0(t)\leq K(\epsilon) e^{4\epsilon(\upgamma\alpha r_j)t/c}
\]
where the constant $K(\epsilon)$ depends on $\mu(c)\tfrac{R_*}{c}$ and on the initial state. By choosing $\tilde \epsilon=3\epsilon$, we obtain the desired result. 
In turn, since the mean momentum tends to 0 exponentially,
\[
\langle q^{\mathrm{lin}}\rangle_j(t) =\langle q^{\mathrm{lin}}\rangle_j(0)+\ds\int_0^t \langle p^{\mathrm{lin}}\rangle_j(s)\ud s\]
admits a finite limit, given by $q_{j,\infty}:=\langle q^{\mathrm{lin}}\rangle_j(0)+\ds\int_0^\infty \langle p^{\mathrm{lin}}\rangle_j(s)\ud s$. Moreover, 
\begin{equation*}
  |\langle q^{\mathrm{lin}}\rangle_j(t)-q_{j,\infty}|\le \int_t^\infty |\langle p^{\mathrm{lin}}\rangle_j(s)|\ud s \le \frac{K(\epsilon)c}{\upgamma \alpha r_j(1-\epsilon)}e^{-\upgamma \alpha r_j(1-\epsilon) t/c}.
\end{equation*}
\QED

\begin{rmk}
  Note that the choice of $c_*(\epsilon)$ can be done uniformly over $j$ by replacing $r_j$ by $P$ in \eqref{eq_cond_c_eps}. The rate of convergence of the mean momentum is still given by $e^{-\upgamma\alpha r_j(1-\epsilon)t/c}$ and the convergence is then faster when $r_j$ is bigger.
\end{rmk}

\section{The two-particles case}
\label{TwoPart}

Let us now focus on the case of two particles. This case is interesting since we are able to have a better characterization of the large time behavior of the solutions. 
The equilibrium solutions are defined by \eqref{EqSol}. To be more specific we distinguish two cases
\begin{itemize}
\item [\textcircled{{\small 1}}] Separate particles: we have $L=2$, $r_1=r_2=1$, $|q^*_1-q^*_2|>2R_1$,
\item [\textcircled{{\small 2}}] Packed particles: we have $L=1$, $r_1=r_2=2$,
$q^*_1=q^*_2$ (that can be supposed to be 0 without loss of generality).
\end{itemize}
In both cases, we have $p^*_1=p^*_2=0$.

\subsection{The linearized system: damping and oscillation}

Expliciting 
\eqref{Lin_nonaut}, for the two particles case, we get
\begin{align*}
&\dot q^{\mathrm{lin}}_j =p^{\mathrm{lin}}_j,
\\
&\dot q^{\mathrm{lin}}_j =-\alpha r\left(\kappa q^{\mathrm{lin}}_j-\ds\int_0^t \mathpzc{p}_c(t-s) \langle q^{\mathrm{lin}}\rangle_j(s)\ud s\right)+S_{I,j}
\end{align*}
where 
\begin{itemize}
\item  $r=1$, $\langle q^{\mathrm{lin}}\rangle_j=q^{\mathrm{lin}}_j$ for case \textcircled{{\small 1}},
\item  $r=2$, $\langle q^{\mathrm{lin}}\rangle_1=\langle q^{\mathrm{lin}}\rangle_2=\tfrac{q^{\mathrm{lin}}_1+q^{\mathrm{lin}}_2}2$ for case \textcircled{{\small 2}}.
\end{itemize}
For case \textcircled{{\small 1}}, the particles are uncloupled and the behavior of the solutions of the linearized system is completely characterized by Proposition~\ref{prop:meanLin}.

Let us push forward the discussion about case \textcircled{{\small 2}}.
It is convenient to  reformulate the problem by setting 
$$\delta^\pm q=\frac{q^{\mathrm{lin}}_1\pm q^{\mathrm{lin}}_2}{2},\qquad \delta^\pm p=\frac{p^{\mathrm{lin}}_1\pm p^{\mathrm{lin}}_2}{2}.$$ 
Of course,  $\delta^+q, \delta^+p$ coincide with the (common) mean value $\langle q^{\mathrm{lin}}\rangle_j, \langle p^{\mathrm{lin}}\rangle_j$
for $j\in \{1,2\}$, the behavior of which 
has been characterized  in Proposition~\ref{prop:meanLin}.
 Next, $(\delta^-q,\delta^-p)$ satisfies
  \[ \ds\frac{\ud}{\ud t} 
 \begin{pmatrix}
 \delta^- q\\ \delta^- p\end{pmatrix}=
  \begin{pmatrix} 0 & \mathbb I
 \\
 -2\kappa \alpha \mathbb I & 0 
 \end{pmatrix} \begin{pmatrix}
 \delta^- q\\ \delta^- p\end{pmatrix}
 +
 \begin{pmatrix}
0\\ \delta^-S_I\end{pmatrix}.
 \]
  According to Lemma~\ref{lem:SI}
  $\delta^- S_I$ is compactly supported. Up to this perturbation that vanishes in finite time, the dynamics is governed 
  by a matrix having 
   the purely imaginary eigenvalues $\pm i\sqrt{2\kappa \alpha}$. Therefore, 
$\delta q^-, \delta p^-$ have an oscillating behavior while   $\delta q^+, \delta p^+$ are subjected to damping. 

\subsection{The nonlinear problem}

We consider the ``well-prepared'' case for two particles, which means 
that the mean velocity is set to 0 from the beginning. For the sake of simplicity, we neglect 
the initial data for the wave equation, assuming $\Psi_0=0$, $\Psi_1=0$.
We have $\frac{\ud}{\ud t}q_1=p_1$, $\frac{\ud}{\ud t}q_2=p_2$ together with
\[\begin{array}{l}
\ds\frac{\ud}{\ud t}p_1(t)=\ds\int_0^t  \mathpzc{p}_c(t-s)\big( \nabla\Sigma_1(q_1(t)-q_1(s))
+\nabla\Sigma_1(q_1(t)-q_2(s))\big)\ud s,
\\\ds\frac{\ud}{\ud t}p_2(t)=\ds\int_0^t  \mathpzc{p}_c(t-s)\big( \nabla\Sigma_1(q_2(t)-q_1(s))
+\nabla\Sigma_1(q_2(t)-q_2(s))\big)\ud s.
\end{array}\]
Let us set 
$$\delta q^\pm=\ds\frac{q_1\pm q_2}{2}, \qquad \delta p^\pm=\ds\frac{p_1\pm p_2}{2},$$
so that 
\[
q_1=\delta^+q+\delta^-q,\quad p_1=\delta^+p+\delta^-p,\quad
q_2=\delta^+q-\delta^-q,\quad p_2=\delta^+p-\delta^-p
 .\]
 Therefore, we can equivalently express the problem as a system on the new variables:
$\frac{\ud}{\ud t}\delta^\pm q=\delta^\pm p$, together with
\begin{eqnarray}
\ds\frac{\ud}{\ud t}\delta^+p(t)=\ds\frac12\ds\int_0^t  \mathpzc{p}_c(t-s)\big( \nabla\Sigma_1(\delta^+q(t)-\delta^+q(s)
+\delta^-q(t)-\delta^-q(s))
\nonumber\\
\hspace*{3cm}+\nabla\Sigma_1(\delta^+q(t)-\delta^+q(s)+\delta^-q(t)+\delta^-q(s))
\nonumber\\
\hspace*{3cm}
+\nabla\Sigma_1(\delta^+q(t)-\delta^+q(s)-\delta^-q(t)-\delta^-q(s))
\nonumber\\
\hspace*{6cm}+\nabla\Sigma_1(\delta^+q(t)-\delta^+q(s)-\delta^-q(t)+\delta^-q(s))
\big)\ud s,
\label{eq1}
\\[.3cm]
\ds\frac{\ud}{\ud t}\delta^-p(t)=\ds\frac12\ds\int_0^t  \mathpzc{p}_c(t-s)\big( \nabla\Sigma_1(\delta^+q(t)-\delta^+q(s)
+\delta^-q(t)-\delta^-q(s))\nonumber\\
\hspace*{3cm}
+\nabla\Sigma_1(\delta^+q(t)-\delta^+q(s)+\delta^-q(t)+\delta^-q(s))
\nonumber\\
\hspace*{3cm}
-\nabla\Sigma_1(\delta^+q(t)-\delta^+q(s)-\delta^-q(t)-\delta^-q(s))
\nonumber\\
\hspace*{6cm}
-\nabla\Sigma_1(\delta^+q(t)-\delta^+q(s)-\delta^-q(t)+\delta^-q(s))\big)\ud s.\label{eq2}
\end{eqnarray}
Since $\sigma_1$ is even, $\Sigma_1$ is even too and $\nabla\Sigma_1$ is odd.
Hence, we observe that if $\delta^+q$ vanishes initially, it remains equal to 0 forever
(if $\delta^+q=0$ in the right hand side of the  first ODE of \eqref{eq1}, the right hand side vanishes).
Accordingly, 
when $(\delta^+q,\delta^+p)\big|_{t=0}=0$, the system restricts to an equation for $(\delta^-q,\delta ^-p)$:
\begin{equation}\label{reduit2}
\left\{\begin{array}{l}
\ds\frac{\ud}{\ud t}\delta^- q=\delta^- p,
\\[.3cm]
\ds\frac{\ud}{\ud t}\delta^- p=
\ds\int_0^t  \mathpzc{p}_c(t-s)\big( \nabla\Sigma_1(\delta^-q(t)-\delta^-q(s))+\nabla\Sigma_1(\delta^-q(t)+\delta^-q(s))
\big)\ud s.
\end{array}\right.
\end{equation}
From now on, we simplify the notation, using $(q,p)$ instead of $(\delta^-q,\delta^-p)$.

\begin{theo}\label{th_wp}
We suppose $n=3$ and $d=1$.
There exists $\delta_0,c_0>0$ such that, if the initial state is small enough $|(q_0,p_0)|\leq \delta_0$ and the wave speed is large enough
$c\geq c_0$, then, we can find a constant $K>0$ such that
the associated solution 
$t\mapsto (q,p)(t)$ of \eqref{reduit2}  satisfies
\[|(q,p)(t)|\leq \ds\frac{K}{\sqrt{1+t/c}}.\]
\end{theo}

\noindent{
{\bf{Proof.}}
We use the elementary relations
\begin{align*}
\nabla\Sigma_1(a\pm b)=&\nabla\ds\int \sigma_1(a\pm b -y)\sigma_1(y)\ud y
=\nabla\ds\int \sigma_1(a -y)\sigma_1(y\pm b)\ud y\\
=&-\ds\int \nabla\sigma_1(y-a)\sigma_1(y\pm b)\ud y
\end{align*}
with $a=q(t)$, $b=q(s)$. It leads to 
\begin{equation}\label{pdot}
\begin{array}{l}
\ds\frac{\ud}{\ud t} p
=
\ds\int_0^t \ds\int \mathpzc{p}_c(t-s) \big( -\nabla\sigma_1(y-q(t))\big)
\big\{\sigma_1(y+q(s))+\sigma_1(y-q(s))
\big\}\ud y
\ud s.
\end{array}
\end{equation}
In the integrand, the function $y\mapsto \{...\}$ is even
while $\nabla\sigma_1(-y-q(t))=-\nabla\sigma_1(y+q(t))$.
It follows that 
\begin{equation}\label{pdot2}
\ds\frac{\ud}{\ud t} p
=
-\frac{1}{2} \ds\int_0^t \ds\int \mathpzc{p}_c(t-s)\big( \nabla\sigma_1(y-q(t))
-\nabla\sigma_1(y+q(t))\big)
\big\{\sigma_1(y+q(s))+\sigma_1(y-q(s))
\big\}\ud y
\ud s.
\end{equation}
We split the rhs as $\text{\textcircled{\tiny{A}}}+ \text{\textcircled{\tiny{B}}}$
with 
\begin{align*}
\text{\textcircled{\tiny{A}}}(t)
=& \ds -\frac{1}{2} \ds\int_0^t \ds\int \mathpzc{p}_c(t-s)\big( \nabla\sigma_1(y-q(t))
-\nabla\sigma_1(y+q(t))\big)
\\
&\hspace*{4cm}\big(\sigma_1(y+q(t))+\sigma_1(y-q(t))
\big)\ud y
\ud s
\\
=& \ds -\frac{1}{2} \ds\int_0^t  \mathpzc{p}_c(t-s)
\ud s
\ds\int
\big( \nabla\sigma_1(y-q(t))\sigma_1(y+q(t))
-\nabla\sigma_1(y+q(t))\sigma_1(y-q(t))
\big)
\ud y
\end{align*}
since 
\[
\ds\int \sigma_1(y+a)\nabla\sigma_1(y+a)\ud y
=\ds\frac12\ds\int \nabla|\sigma_1(y+a)|^2\ud y
=0.\]
Hence, we get 
\[
\text{\textcircled{\tiny{A}}}(t)
=
\ds\int_0^t  \mathpzc{p}_c(t-s)
\ud s
\ds\int
 \nabla\sigma_1(y+2q(t))\sigma_1(y)\ud y
=
\ds\int_0^t  \mathpzc{p}_c(t-s)
\ud s
\nabla\Sigma_1(2q(t)).
\] 
Next, 
we have
\begin{align*}
\text{\textcircled{\tiny{B}}}(t)
=&\,
\ds \frac{1}{2} \ds\int_0^t \ds\int \mathpzc{p}_c(t-s)\big( \nabla\sigma_1(y-q(t))
-\nabla\sigma_1(y+q(t))\big)
\\
&\hspace*{3cm}
\big(\sigma_1(y+q(t))-\sigma_1(y+q(s))+\sigma_1(y-q(t))-\sigma_1(y-q(s))
\big)\ud y
\ud s
\\[.3cm]
=&\,
\ds \frac{1}{2} \ds\int_0^t \ds\int \mathpzc{p}_c(t-s)\big( \nabla\sigma_1(y-q(t))
-\nabla\sigma_1(y+q(t))\big)
\\
&\hspace*{3cm}\ds\int_s^t 
\underbrace{\ds\frac{\ud}{\ud\tau}
\big(\sigma_1(y+q(\tau))+\sigma_1(y-q(\tau))\big)}_{
=\big(\nabla\sigma_1(y+q(\tau))-\nabla\sigma_1(y-q(\tau))\big)\cdot p(\tau)}
\ud\tau \ud y
\ud s
\end{align*}
and using integration by parts, we are led to
\[
\begin {array}{lll}
\text{\textcircled{\tiny{B}}}(t)
&=&
\ds -\frac{1}{2} \ds\int_0^t \ds\int  (t-s)\mathpzc{p}_c(t-s)  
\big( \nabla\sigma_1(y-q(t))
-\nabla\sigma_1(y+q(t))\big)\\
&&\hspace*{6cm}
\otimes \big( \nabla\sigma_1(y-q(t))
-\nabla\sigma_1(y+q(t))\big) p(t)
\ud y \ud s
\\
[.3cm]
&&
\ds + \frac{1}{2} \ds\int_0^t \ds\int \mathpzc{p}_c(t-s)\big( \nabla\sigma_1(y-q(t))
-\nabla\sigma_1(y+q(t))\big) \times
\\
[.3cm]
&&\hspace*{2cm}\left(\ds\int_s^t (\tau -s)
\ds\frac{\ud}{\ud\tau}\Big[
\big(\nabla\sigma_1(y-q(\tau))-\nabla\sigma_1(y+q(\tau))\big) \cdot p(\tau)\Big]
\ud\tau \right)\ud y
\ud s.
\end{array}
\]
It can be interpreted as a Taylor expansion 
$f(t)-f(s)=f'(t)(t-s)-\int_s^t f''(\tau)(\tau-s)\ud \tau$ with 
$f(\tau)=\sigma_1(y+q(\tau))+\sigma_1(y-q(\tau))$.
We write this expression as
\begin{align*}
&\text{\textcircled{\tiny{B}}}(t)
\\
&= -\frac{1}{2} \int_0^t \tau\, \mathpzc{p}_c(\tau)\ud\tau
\bigg( \int \big(\nabla\sigma_1(y-q(t)) - \nabla\sigma_1(y+q(t))\big)
\\
&\qquad\qquad\qquad \otimes \big(\nabla\sigma_1(y-q(t)) - \nabla\sigma_1(y+q(t))\big) \ud y \bigg)\, p(t)
\\
&\qquad + \text{\textcircled{\tiny{C}}}(t),
\end{align*}
where the first term in the rhs recasts as
\[\begin{array}{l}
\ds -
\ds\int_0^t \tau \mathpzc{p}_c(\tau)\ud \tau
\left(\ds\int \big( \nabla\sigma_1(y)\otimes \nabla\sigma_1(y)
+\nabla\sigma_1(2q(t)-y)\big)\otimes  \nabla\sigma_1(y)
\ud y\right)p(t)
\\[.3cm]
\qquad = +\ds 
\ds\int_0^t \tau \mathpzc{p}_c(\tau)\ud \tau
(D^2\Sigma_1(0)-D^2\Sigma_1(2q(t)))p(t)
\end{array}
\]
and where we denote
\begin{multline*}
\text{\textcircled{\tiny{C}}}(t)=\frac{1}{2} \ds\int_0^t \ds\int \mathpzc{p}_c(t-s)\big( \nabla\sigma_1(y-q(t))
-\nabla\sigma_1(y+q(t))\big) 
\\
\left(\ds\int_s^t (\tau -s)
\ds\frac{\ud}{\ud\tau}\Big[
\big(\nabla\sigma_1(y-q(\tau))-\nabla\sigma_1(y+q(\tau))\big) \cdot p(\tau)\Big]
\ud\tau \right)\ud y
\ud s.
\end{multline*}
Using Lemma~\ref{lemAV} and \ref{mtpc}, we arrive at 
\begin{equation}\label{Lie_Rem}
\ds\frac{\ud}{\ud t}p(t)
=
\kappa \nabla\Sigma_1(2q(t)) +\ds\frac{\upgamma}{
c}
(D^2\Sigma_1(0)-D^2\Sigma_1(2q(t)))p(t)
+\text{\textcircled{\tiny{C}}}(t)
+\text{\textcircled{\tiny{D}}}(t)\end{equation}
where 
\[
\text{\textcircled{\tiny{D}}}(t)
=
-\ds\int_t^\infty \mathpzc{p}_c(\tau)\ud \tau \nabla\Sigma_1(2q(t))
-
 \ds\int_t^\infty \tau \mathpzc{p}_c(\tau)\ud \tau (D^2\Sigma_1(0)-D^2\Sigma_1(2q(t)))p(t)
\]
vanishes for $t\geq \frac{2R_2}{c}$.
As we are going to see, the first term provides an oscillating behavior, the second gives a \emph{non linear} damping
and we will have to manage with the remainder terms that are expected not to compromise this behavior.

From now on, we continue the discussion by restricting the discussion to the one dimensional case ($d=1$).
We suppose that the function $x\mapsto \sigma_1(x)$ is even, 
and monotonically decreasing on $[0,\infty)$,
so that 
$$\Sigma_1: x\longmapsto \ds\int_{\mathbb R} \sigma_1(x-y)\sigma_1(y)\ud y$$
is even too. 
Moreover, since $\sigma_1$ is even, $\sigma_1'$ is odd, and we have, for any $x\geq 0$,
\begin{align*}
\Sigma_1'(x)=&\ds\int_{\mathbb R} \sigma_1'(x-y)\sigma_1(y)\ud y
=\ds\int_{\mathbb R} \sigma_1(x-y)\sigma_1'(y)\ud y
=\ds\int_{0}^\infty \sigma_1'(y)\big(\sigma_1(x-y)-\sigma_1(x+y)\big)\ud y
\\[.3cm]
=&
\ds\int_{0}^x \sigma_1'(y)\big(\sigma_1(x-y)-\sigma_1(x+y)\big)\ud y
+
\ds\int_{x}^\infty \sigma_1'(y)\big(\sigma_1(y-x)-\sigma_1(y+x)\big)\ud y
\leq 0
\end{align*}
and $\Sigma_1$ is monotonically decreasing on $[0,\infty)$.
A typical shape is given in Fig.~\ref{FSig1}.
The derivatives of odd order of $\Sigma_1$ vanish at $x=0$, and 
$\Sigma_1''(0)=-\|\sigma_1'\|_{L^2}^2<0$: $\Sigma_1$
 reaches its maximum $\|\sigma_1\|_{L^2}^2$ at $x=0$.  
 We also note that $\Sigma_1^{(4)}(0)=\|\sigma_1''\|_{L^2}^2>0$.
 Then, we can pick some $\delta>0$ 
 so that 
 \[
 \textrm{for any $q\in [-\delta,\delta]$, 
 $\sigma_*\leq \Sigma_1(q)$,
 $0<\alpha_* \leq -\Sigma_1''(q)$ and $0<\eta_*\leq \Sigma_1^{(4)}(q)\leq \eta^*$,}\]
 for some positive constants $\alpha_*,\eta_*,\eta^*,\sigma_*$.
 
\begin{figure}
\centering
\begin{tikzpicture}[scale=3]
\draw [thin,->] (-1.3,0) -- (1.3,0);
 \draw [thin,->] (0,-.2) -- (0,.6);
  \draw [dotted] (.3,0) -- (.3,{exp(-1/(1-.3*.3))});
  \draw [dotted] (-.3,0) -- (-.3,{exp(-1/(1-.3*.3))});
\node at (.3,-.1) {$\delta$};
\node at (-.3,-.1) {$-\delta$};
\draw [domain=-0.999:0.999, thick] plot (\x,{exp(-1./(1-\x*\x)});
\end{tikzpicture}
\caption{Shape of the function $\Sigma_1$}
\label{FSig1}
\end{figure}
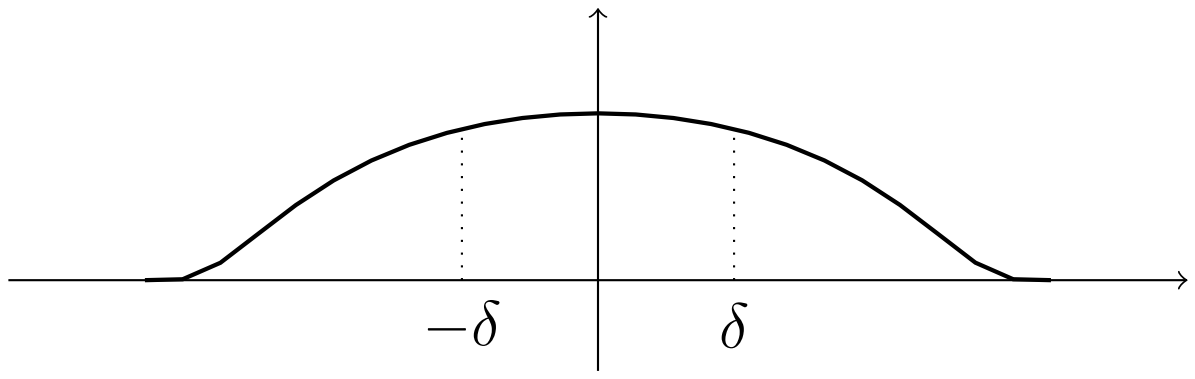

Coming back to \eqref{Lie_Rem} and  keeping only the (expected) leading terms, we get
\[
\ds\frac{\ud^2}{\ud t^2}q(t)
=-2\kappa \|\sigma_1'\|_{L^2}^2 \ q(t)
-\ds  \frac{\upgamma}{c} \|\sigma_1''\|_{L^2}^2 q^2(t)\ds\frac{\ud}{\ud t}q(t).
\]
We recognize the Li\'enard equation \cite{Lien}, for which the solutions behave like $\mathscr O(1/\sqrt t)$ for large times. 
We need to manage with estimates on the remainders in  \eqref{Lie_Rem}.
We proceed by averaging methods \cite[Chapter~11]{Ver}.

In what follows, we denote by $K_0, K_1...$  positive constants that depend on the parameters of the model, but $c$, possibly keeping the same notation while the precise value of the constant may vary from a line to another.
 Taking advantage of the one dimensional framework ($d=1$), we introduce new variables by setting
 \[
 \mathcal E=\ds\frac{p^2}{2}+ \frac{\kappa}{2}(\Sigma_1(0)-\Sigma_1(2q)),\
 \frac{p}{\sqrt{2}}=\sqrt{\mathcal E} \sin(\theta),\
 \sqrt{\frac{\kappa}{2}} \ \mathrm{sgn}(q)(\sqrt{\Sigma_1(0)-\Sigma_1(2q)})=\sqrt{\mathcal E} \cos(\theta)
 .\]
 
 As said above, owing to the monotonicity properties of $\Sigma_1$, we can define $\delta_*>0$ such that 
 \[\textrm{$0\leq \frac{\kappa}{2}(\Sigma_1(0)-\Sigma_1(2q))\leq \delta_*$ iff $|q|\leq \delta $.}\]
 By writing
\[ \Sigma_1(0)-\Sigma_1(2q)=-\ds\int_0^1 (1-\mu)\Sigma_1''(2\mu q)4q^2\ud\mu\] 
  in this domain, we equally have the estimates
 \begin{equation}\label{cvx1}
 |q|^2\leq \ds\frac{\Sigma_1(0)-\Sigma_1(2q)}{2\alpha_*}\leq \ds\frac{\mathcal E}{\alpha_*\kappa},
 \end{equation}
 and
 \begin{equation}\label{cvx2}
 \begin{array}{l}
\ds\frac{\eta_*}{\| \Sigma_1''\|_{L^\infty}^2}
(\Sigma_1(0)-\Sigma_1(2q))
\leq 
2 \eta_* |q|^2\leq \Sigma_1''(2q)-\Sigma_1''(0)=\ds\int_0^1 (1-\mu)\Sigma_1^{(4)}(2\theta q)4|q|^2\ud\mu
\\
[.3cm]
\hspace*{9cm}
\leq  
2\eta^* |q|^2 \leq\ds\frac{\eta^*}{\alpha_*} (\Sigma_1(0)-\Sigma_1(2q)).
 \end{array} \end{equation}

 By computing $\frac{\ud}{\ud t}\tan(\theta(t))$, we find
 \begin{equation}\label{sysEth}\begin{array}{l}
 \dot {\mathcal E}
 =\ds\frac{\upgamma}{c}(\Sigma_1''(0)-\Sigma_1''(2q))p^2 +(\text{\textcircled{\tiny{C}}}+ \text{\textcircled{\tiny{D}}})p
,\\[.3cm]
\dot \theta
=\ds \underbrace{\frac{\kappa \ \Sigma_1'(2q)}{\sqrt{2\mathcal{E}}\cos(\theta)}}_{=\frac{\kappa \ \mathrm{sgn}(q)\Sigma_1'(2q)}{\sqrt{\kappa} (\Sigma_1(0)-\Sigma_1(2q))}}+\frac{\gamma}{ c}\left(\Sigma_1''(0)-\Sigma_1''(2q)\right) \sin(\theta) \cos(\theta)+\frac{(\text{\textcircled{\tiny{C}}}+\text{\textcircled{\tiny{D}}})\cos(\theta)}{\sqrt{2\mathcal{E}}}
\end{array}
 \end{equation}
 
 \begin{rmk}
 The standard Li\'enard equation has the form $\ddot q=-\omega^2 q-\gamma q^2 \dot q$, with $\omega,\gamma>0$.
 The corresponding change of unknowns simply
 reads $\mathcal E=\frac{p^2}{2}+\frac{\omega^2 q^2}{2}$, $p=\sqrt {2\mathcal E} \sin(\theta)$, 
 $\omega q= \sqrt {2\mathcal E} \cos(\theta)$.
 Hence, we have $\dot {\mathcal  E}=-\frac\gamma{\omega^2}\sin^2(2\theta) \mathcal E^2$ and 
 $\dot \theta=\omega(-1- \frac{2\mathcal E\gamma}{\omega^3}\sin^3(\theta)\cos(\theta))$.
 The energy $\mathcal E$ is non increasing: its derivative can vanish, but, for small enough energy, we have $\dot \theta\leq \frac{-\omega}{2}<0$ and the energy cannot exhibit a flat profile over long time intervals.
 Quantifying this intuition gives the decay rate in $\mathscr O(1/\sqrt t)$. Here, the contribution of 
 $\text{\textcircled{\tiny{C}}}+ \text{\textcircled{\tiny{D}}}$ compromise the monotonicity, and one has to evaluate more finely the variation of $\mathcal E$ and $\theta$.
 \end{rmk}
 
We start 
 from an initial data such that  $$ \mathcal 
 E(0)=\mathcal E_0< \ds\frac{\delta_*} 4
 .$$
Coming back to \eqref{pdot}, we see that
 $$|\dot{p}(t)| \leqslant \|\sigma_1'\|_{L^\infty}\|\sigma_1\|_{L^1} \| \mathpzc{p}_c\|_{L^1}$$
which allows us to derive  the rough estimate 
\[\begin{array}{lll}
 |p(t)|&\leq& |p(0)|+\ds\int_0^t |\dot p(s)|\ud s
%
\leq  \sqrt{2\mathcal E_0}
 +
\|\sigma_1'\|_{L^\infty}\|\sigma_1'\|_{L^1} \| \mathpzc{p}_c\|_{L^1} t.\end{array}\]
We bear in mind that, with the adopted scaling,  $\| \mathpzc{p}_c\|_{L^1} =\| \mathpzc{p}_1\|_{L^1}$ 
and the estimate actually does not depend on $c$.
Combining this estimate to  
 \[
 \kappa(\Sigma_1(0)-\Sigma_1(2q(t)))=
 \kappa(\Sigma_1(0)-\Sigma_1(2q(0)))
 -2\kappa\ds\int_0^t \Sigma_1'(2q(\tau))p(\tau)\ud \tau 
 \]
 shows that 
 \[
 \mathcal E(t)\leq \mathcal{E}_0+\mathfrak m_1t+\mathfrak m_2t^2,
 \]
 with positive constants $\mathfrak m_1,\mathfrak m_2$ 
 depending on the parameters  of the model, but $c$, and the initial energy. 
 In particular, using this estimate with  $t\in [0,\frac{2R_2}c]$, and by taking $c$ large enough, 
 we obtain 
 \[
 \mathcal E(t)\leq 2\mathcal E_0<\ds\frac{\delta_*}2 
 \]
 on this time interval.
 Therefore, we can consider what happens beyond $\frac{2R_2}c$ and disregard the $\text{\textcircled{\tiny{D}}}$ term in \eqref{sysEth}.
 We set 
 \[
 T=\sup\{t\geq 0: \forall s \in [0,t], \ \mathcal E(s)\leq 4\mathcal E_0<\delta_*
  \}>\ds\frac{2R_2}{c}.\]
 On the time interval $[\frac{2R_2}{c},T]$, we can make full use of the basic estimates \eqref{cvx1} and \eqref{cvx2}.

We have \[\begin{array}{l}
 \text{\textcircled{\tiny{C}}}(t)=
 \ds \frac{1}{2} \ds\int_0^t \ds\int \mathpzc{p}_c(t-s)\big( \sigma_1'(y-q(t))
-\sigma_1'(y+q(t))\big)
\\
[.3cm]
\hspace*{4cm}\left(\ds\int_s^t (\tau -s)
\ds\frac{\ud}{\ud\tau}\Big[
\big(\sigma_1'(y-q(\tau))-\sigma_1'(y+q(\tau))\big) p(\tau)\Big]
\ud\tau \right)\ud y
\ud s
\\[.3cm]
= \ds \frac{1}{2} \ds\int 
_0^t \ds\int \mathpzc{p}_c(t-s)\big( \sigma_1'(y-q(t))
-\sigma_1'(y+q(t))\big)
\\
[.3cm]
\hspace*{2cm}
\bigg(\ds\int_s^t (\tau -s)
\Big[-\big(\sigma_1''(y-q(\tau))+\sigma_1''(y+q(\tau))\big) |p(\tau)|^2
\\[.3cm]
\hspace*{6cm}
+
\big(\sigma_1'(y-q(\tau))-\sigma_1'(y+q(\tau))\big)\dot p(\tau)
\Big]
\ud\tau \bigg)\ud y
\ud s.
\end{array}\]
We are going to use the estimates 
\[
|\sigma_1'(y-q(\tau))-\sigma_1'(y+q(\tau))|\leq 2\|\sigma_1''\|_{L^\infty} |q(\tau)|\]
and
\begin{align*}
\ds\int | \sigma_1'(y-q(t))
-\sigma_1'(y+q(t))|\ud y
=&\ds\int \left| \ds\int_{-1}^{+1} \ds\frac{\ud}{\ud r} \sigma_1'(y+r q(t))
\ud r\right|\ud y
\\&\leq 
\ds\int_{-1}^{+1} \left(\ds\int  \sigma_1''(y+r q(t)) \ud y\right)
\ud r
|q(t)| =2\|\sigma_1''\|_{L^1} |q(t)|.
\end{align*}
Hence, we obtain
\[
    |\text{\textcircled{\tiny{C}}}(t)| \leqslant  2|q(t)| \ \|\sigma_1''\|_{L^1} 
    \|\sigma_1''\|_{L^\infty}
    \ds\int 
    _0^t 
    |\mathpzc{p}_c(t-s)| \left(\ds\int_s^t (\tau -s)
\Big( |p(\tau)|^2
+
|q(\tau)| |\dot{p}(\tau)|
\Big)
\ud\tau \right)
\ud s.
\]

Going back to \eqref{pdot2}, we obtain  the following refined estimate on $\dot p$:
\[
|\dot{p}(t)| \leqslant 
2\|\mathpzc{p}_c\|_{L^1}\|\sigma_1\|_{L^1}\|\sigma_1''\|_{L^\infty}\ |q(t)|.
\]
We also use the basic estimates 
\[|p(t)| \leqslant \sqrt{2 \mathcal{E}(t)},\qquad 
|q(t)|\leqslant \sqrt{\frac{\mathcal{E}(t)}{\alpha_* \kappa}},
\]
(see \eqref{cvx1}) which hold  on $[0,T]$.
For further purposes, let us introduce, for  $t \in [\frac{2R_2}{c},T]$:
\[
\mathcal E_c(t)=\sup\left\{\mathcal E(\tau), \ \tau \in  \Big[t-\frac{2R_2}{c}, t\Big] \subset [0,T] \right\}.\] 
By continuity, $\mathcal{E}_c(t)=\mathcal E(s_t)$ is reached for some $s_t\in [t-\frac{2R_2}{c},t]$ ($s_t$ denoting the largest time where $\mathcal E_c(t)$ is reached).
Bearing in mind that $\mathpzc{p}_c$ is supported in $[0,\frac{2R_2}{c}]$, 
we then estimate, for $t\in [\frac{2R_2}{c},T]$,
\[\begin{array}{lll}
|\text{\textcircled{\tiny{C}}}(t)| &\leqslant&
 2 \sqrt{\ds\frac{\mathcal{E}(t)}{\alpha_* \kappa}} \ \|\sigma_1''\|_{L^1} \|\sigma_1''\|_{L^\infty} 
 \\[.3cm]&&
\qquad\times\ds\int_{t-\frac{2R_2}{c}}^t |\mathpzc{p}_c(t-s)|  \left(\ds\int_s^t (\tau -s)
\Big(1+\ds\frac{\|\mathpzc{p}_c\|_{L^1}\|\sigma_1\|_{L^1}\|\sigma_1''\|_{L^\infty}
}{\alpha_* \kappa}
\Big)2\mathcal{E}_c(t)
\ud\tau \right)
\ud s \\
[.3cm]&\leqslant& K_0 \mathcal{E}_c(t) \sqrt{\mathcal{E}(t)}  \ds\int_{0}^{\frac{2R_2}{c}} s^2 |\mathpzc{p}_c(s)| \ud s.
\end{array}\]

Finally, by coming back to the definition of the kernel $\mathpzc{p}_c$, we arrive at
\begin{equation}\label{estCp}
\big|\text{\textcircled{\tiny{C}}}(t)\big|
\leq \ds\frac{K_0}{c^2} \mathcal{E}_c(t) \sqrt{\mathcal{E}(t)} ,\qquad 
\big|\text{\textcircled{\tiny{C}}}(t) p(t)\big|
\leq 
\ds\frac{K_0}{c^2} \mathcal{E}_c(t) \mathcal{E}(t), \end{equation}
for $t\in [\frac{2R_2}{c},T]$.
Moreover, for $t\in [0, T]$, we have $  \mathcal{E}(t)\leq \mathcal{E}_c(t)\leq 4 \mathcal{E}_0 < \delta_*$
and, as already mentioned,  we can appeal to the inequalities \eqref{cvx1} and \eqref{cvx2}.
Integrating the first equation in \eqref{sysEth} yields, for $t \in [\frac{2R_2}{c},T]$
and assuming $c\geq 1$ if necessary,
\[\begin{array}{lll}
\mathcal E(t)-\mathcal{E}_c(t)
=
\mathcal E(t)-\mathcal E(s_t)
&=&\ds\int_{s_t}^t \left(
+\ds\frac{\upgamma}{c}(\Sigma_1''(0)-\Sigma_1''(2q(\tau)))p^2(\tau)
+p(\tau) \text{\textcircled{\tiny{C}}}(\tau)
\right)\ud \tau 
\\[.3cm]&\geq& 
-\ds\int_{s_t}^t
\left(\ds\frac{2\upgamma \eta^*}{c} q^2(\tau) p^2(\tau) 
+|p(\tau) \text{\textcircled{\tiny{C}}}(\tau)|
\right)\ud \tau
\\[.3cm]
\\
&\geq&
-\ds\frac{4\upgamma \eta^*}{\alpha_*\kappa c} \ds\int_{s_t}^t \mathcal E(\tau)^2 \ud \tau -\ds\frac {K_0}{c^2}\ds\int_{s_t}^t \mathcal E(\tau)\mathcal E_c(\tau) \ud \tau
\\[.3cm]
\\
&\geq&
-\ds  \frac{K_1(t-s_t)}{c} \mathcal{E}_0\mathcal{E}_c(t) \geq \ds -\frac{K_1}{c^2} \mathcal{E}_0\mathcal{E}_c(t) 
\end{array}\]
We deduce that
\begin{equation}\label{ratio}
\ds\frac{\mathcal E(t)}{\mathcal E_c(t)}\geq 
\underbrace{1 -\frac{K_1}{c^2} \mathcal{E}_0}_{:=\rho(c)}
\end{equation}
where we can pick $c$ large enough for this quantity to lie in $(0,1)$.
Note that how large  $c$ is chosen depends on $\mathcal E_0$ thus on the initial state.
As a consequence, we obtain
\[
\left|\ds\frac{ \text{\textcircled{\tiny{C}}}p(t)}{\mathcal {E}^2(t)}\right|\leq
\ds\frac{K_0}{c^2}\ds\frac{\mathcal{E}_c(t)}{\mathcal E(t)}\leq \ds\frac{K_0}{c^2\rho(c)}
\]
which, again, can be made arbitrarily small by choosing $c$ large enough.
Next, the first term in the equation for $\theta$ in \eqref{sysEth}
reads
\[
\ds\frac{\kappa \ \mathrm{sgn}(q)}{\sqrt{\kappa}(\Sigma_1(0)-\Sigma_1(2q))}
\ds\int_0^1 \Sigma_1''(2\mu q) 2q\ud \mu
=
\ds\frac{2 \kappa|q|}{\sqrt{\kappa}(\Sigma_1(0)-\Sigma_1(2q))}
\ds\int_0^1 \Sigma_1''(2\mu q) \ud \mu<0,
\]
and it thus lies in the interval $$
[-a,-b],\qquad a=\ds \sqrt{\frac{2\kappa}{\alpha_*}} \ \|\Sigma_1''\|_{\infty} , \qquad b=
\ds \sqrt{\frac{2\kappa}{\|\Sigma_1''\|_{\infty}}} \ \alpha_*.$$
On the same token, we obtain for the second term
\[
\ds\frac{1}{\sqrt 2}
\left|\ds\frac{\upgamma}{c}\sin(\theta)\cos(\theta)\ds (\Sigma_1''(0)-\Sigma_1''(2q))\right|
\leq \frac{\sqrt 2\eta_* \upgamma}{ c} |q|^2 \leq 
\frac{\sqrt 2\eta_* \upgamma}{ c} \ \delta_*^2\]
while the third term satisfies
$$\left| \frac{\text{\textcircled{\tiny{C}}}\cos(\theta)}{\sqrt{\mathcal{E}}} \right| \leq
\frac{K_0}{c^2} \mathcal{E}_0$$
We deduce that
\begin{equation}\label{thdot}
-2a\leq \dot \theta\leq -b/2
\end{equation}
holds, provided $c$ is large enough.

We go back to the equation for $\mathcal E$ in \eqref{sysEth}
that can be cast as
\[\dot {\mathcal E}=\ds\frac{\upgamma}{c}\ds\frac{\Sigma_1''(0)-\Sigma_1''(2q)}{\kappa(\Sigma_1(0)-\Sigma_1(2q))}
\underbrace{\frac{\kappa}{2}(\Sigma_1(0)-\Sigma_1(2q))
\ds \frac{p^2}{2}}_{
=\cos^2(\theta)\sin^2(\theta)\mathcal E^2
} + \text{\textcircled{\tiny{C}}}p,
\]
for $t\geq \frac{2R_2}{c}$. We bear in mind that
\[
0<\ds\frac{\eta_*}{\|\Sigma_1''\|_{L^\infty}}\leq
\ds\frac{\Sigma_1''(2q)-\Sigma_1''(0)}{\Sigma_1(0)-\Sigma_1(2q)}
=
\ds\frac{1}{\Sigma_1(0)-\Sigma_1(2q)}\ds\int_0^1 (1-\mu)\Sigma^{(4)}(2\mu q) 4q^2\ud \mu
\leq 
\ds\frac{\eta^*}{\alpha}
.\]
Therefore, using \eqref{estCp} and \eqref{ratio}, we are led to
\begin{equation}\label{estE1}
\dot {\mathcal E} \leq \mathcal{E}^2\left(
-\frac{\upgamma}{4c}\frac{\eta_*}{\|\Sigma_1''\|_{L^\infty}} \sin^2(2\theta)
+\ds\frac{\text{\textcircled{\tiny{C}}}(t) p(t)}{\mathcal{E}^2(t)}
\right)
\leq 
-\frac{\mathcal{E}^2 }{c}  \left(
K_2\sin^2(2\theta)
-\ds\frac{K_0}{c\rho(c)} 
\right)
.\end{equation} 
The dissipation term is proportional to $\sin^2(2\theta)$ and it might vanish; however
the angle $\theta$ does not lie too long in such a region, by virtue of \eqref{thdot}.
To be more specific, we can
 estimate the time between two consecutive zeroes of $\tau\mapsto \sin^2(2\theta(\tau))$, that is when $\theta(\tau)=\frac{k\pi}{2}$, for some $k\in \mathbb Z$: 
with $\sin(2\theta(\tau_j))=0$, $\tau_2>\tau_1$, $\theta(\tau_1)-\theta(\tau_2)=\frac\pi 2$,
we integrate  inequality~\ref{thdot}
between $\tau_1$ and $\tau_2$, to obtain 
\begin{equation}\label{ecart}
\ds\frac\pi{4a}\leq \tau_2-\tau_1\leq \ds\frac\pi b.
\end{equation}
From \eqref{estE1}, we deduce that, for $\frac{2R_2}{c}\leq t_1<t_2\leq T$,
\begin{equation}\label{estE2}\begin{array}{lll}
\ds\frac{1}{\mathcal E(t_2)}-\ds\frac{1}{\mathcal E(t_1)}
&\geq& \ds \frac{1}{c}\ds \int_{t_1}^{t_2} \left( K_2\sin^2(2\theta(\tau)) 
-\ds\frac{K_0}{c\rho(c)} 
\right)\ud \tau
\\[.3cm]
&=&  \ds \frac{(t_2-t_1)}{c}\left(
 \ds\frac{K_2}{t_2-t_1} \left( \int_{t_1}^{t_2}  \sin^2(2\theta(\tau)) \frac{\dot{\theta}(\tau) \ud \tau}{\dot{\theta}(\tau)} \right)
-\ds\frac{K_0}{c\rho(c)} 
\right)\\[.3cm]
&\geq& \ds \frac{(t_2-t_1)}{c}\left(
 \ds \underbrace{\frac{2K_2}{b(t_2-t_1)}\ds\int_{\theta(t_2)}^{\theta(t_1)}  \sin^2(2\theta) \ud \theta}_{\geqslant 0}
-\ds\frac{K_0}{c\rho(c)} 
\right)\\[.5cm]
&\geq& 
 -\ds\frac{K_0(t_2-t_1)}{c^2\rho(c)}.
\end{array}\end{equation}
We argue by contradiction: let us assume that $T$ is finite: $\mathcal E(T)=4\mathcal E_0$.
With $t_2=T$ and $t_1=\frac{2R_2}{c}$, we get
\[
-\ds\frac{1}{4\mathcal E_0}=\frac{1}{4 \mathcal{E}_0}-\frac{1}{2 \mathcal{E}_0}\geq 
\underbrace{\ds\frac{1}{\mathcal E(T)}}_{=\frac{1}{4\mathcal{E}_0}}-\ds\frac{1}{\mathcal E(2R_2/c)}\geq -\frac{K_0}{c^2\rho(c)}\Big(T-\frac{2R_2}{c}\Big)\ds  
,\]
and therefore
\[
T\geq \ds \frac{c^2\rho(c)
}
{4K_0\mathcal{E}_0 }+\ds\frac{2R_2}{c}.\]
Still by choosing $c$ large enough, this can be made greater than $\pi(1/b-1/4a)$ so that 
$[2R_2/c,T]$ contains at least  two consecutive  times $\tau_1,\tau_2$ where $\sin^2(2\theta(\tau))$ vanishes.
Owing to \eqref{ecart}, there is a finite number of such zeros $\tau_1,..,.\tau_{N+1}$ in $[\frac{2R_2}{c},T]$, see Fig.~\ref{FTfin}.

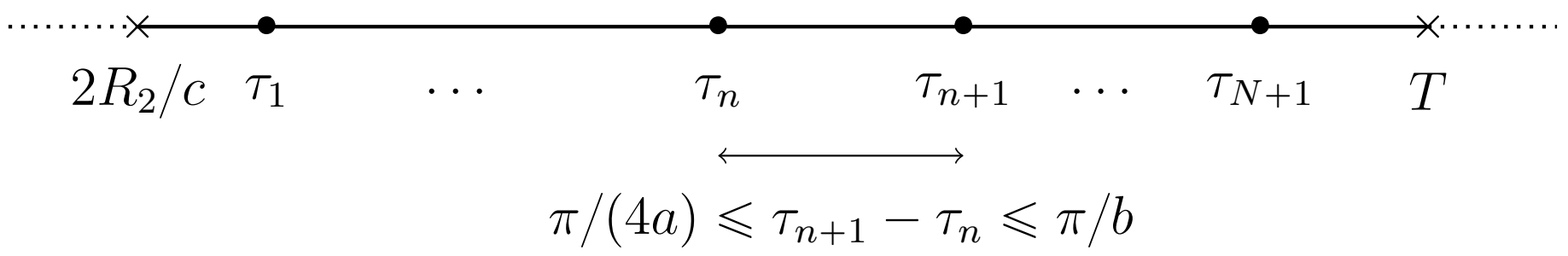
\begin{figure}[!htpb]
\centering
\begin{tikzpicture}[scale=1]
\draw [thick] (0,0) -- (10,0);
\draw [thick, dotted] (-1,0) -- (0,0);
\draw [thick, dotted] (10,0) -- (11,0);
\draw (0,0) node[] {\small{$\times$}};
\draw (10,0) node[] {\small{$\times$}};
\draw (1,0) node[] {\small{$\bullet$}};
\draw (8.7,0) node[] {\small{$\bullet$}};
\node at (1,-.5) {$\tau_1$};
\node at (8.7,-.5) {$\tau_{N+1}$};
\node at (0,-.5) {$2R_2/c$};
\node at (10,-.5) {$T$};
\draw (4.5,0) node[] {\small{$\bullet$}};
\draw (6.4,0) node[] {\small{$\bullet$}};
\node at (4.5,-.5) {$\tau_n$};
\node at (6.4,-.5) {$\tau_{n+1}$};
\node at (2.5,-.5) {$\dots$};
\node at (7.5,-.5) {$\dots$};
\draw [thin,<->] (4.5,-1) -- (6.4,-1);
\node at (5.45,-1.5) {$\pi/(4a)\leq \tau_{n+1}-\tau_n\leq\pi/b$};
\end{tikzpicture}
\caption{Zeroes of $\sin(2\theta(\tau))$} \label{FTfin}
\end{figure}

Coming back to \eqref{estE2}, and using \eqref{thdot} and \eqref{ecart}, we are led to 
\begin{align*}
-\ds\frac{1}{4\mathcal E_0}\geq&\, 
\ds\frac{1}{\mathcal E(T)}-\ds\frac{1}{\mathcal E(2R_2/c)}\geq 
\ds\frac{1}{c}\ds \int_{\frac{2R_2}{c}}^{T} \left( K_2\sin^2(2\theta(\tau)) 
-\ds\frac{K_0}{c\rho(c)} 
\right)\ud \tau
\\[.3cm]
\geq&\, 
\ds\frac{K_2}{c}\ds\sum_{n=1}^N \ds\int_{\tau_n}^{\tau_{n+1}}  \sin^2(2\theta(\tau)) \ud \tau
-\ds\frac{K_0(T-2R_2/c)}{c^2\rho(c)} 
\\[.3cm]
=&\, 
\ds\frac{K_2}{c}\ds\sum_{n=1}^N \frac{\tau_{n+1}-\tau_n}{\tau_{n+1}-\tau_n} \ds\int_{\tau_n}^{\tau_{n+1}}  \frac{\sin^2(2\theta(\tau))}{\dot{\theta}(\tau)} \dot{\theta}(\tau) \ud \tau
-\ds\frac{K_0(T-2R_2/c)}{c^2\rho(c)} 
\\[.3cm]
\geq&\, 
\ds\frac{2K_2}{\pi c}\ds\int_{0}^{\pi/2}  \sin^2(2\theta)\ud \theta\ds\sum_{n=1}^N (\tau_{n+1}-\tau_n)
-\ds\frac{K_0(T-2R_2/c)}{c^2\rho(c)} 
\\[.3cm]
\geq&\, 
\ds\frac{K_2}{2c}(\tau_{N+1}-\tau_1)
-\ds\frac{K_0(T-2R_2/c)}{c^2\rho(c)} 
\\[.3cm]
\geq &\, 
\ds \frac{K_2}{2c}\ds \Big(T- \frac{2 R_2}{c}\Big) \left(1-\ds \frac{2}{c\rho(c)}\ds\frac{K_0}{K_2}
\right)- \ds\frac{K_0}{2c}\Big((T-\tau_{N+1})+\Big(\tau_1-\frac{2R_2}{c}\Big)\Big)
\\[.3cm]
\geq&\, 
 \ds\frac{K_2}{K_0}\ds\frac{c\rho(c)}{8\mathcal E_0} \left(1-\ds \frac{2}{c\rho(c)}
 \ds\frac{K_0}{K_2}
\right)
-\ds \frac{K_0\pi}{c b}.
\end{align*}
(We used the fact that, by construction, $0\leq T-\tau_{N+1}\leq \frac\pi b$ and $0\leq \tau_1-\frac{2R_2}{c}\leq \frac\pi b$.)
However, $c$ can be chosen large enough so that the right hand side is positive, which contradicts that $T$ is finite.

Therefore $T=+\infty$ and, by \eqref{thdot}, there exists an increasing  sequence $\big(\tau_n\big)_{n\in \mathbb N}$ of consecutive zeroes of $\sin^2(2\theta(\tau))$ with $\lim_{n\to \infty}\tau_n=+\infty$.
Adapting the same reasoning as before from \eqref{estE2}, we obtain
\[\begin{array}{lll}
\ds\frac{1}{\mathcal E(\tau_{n})}-\ds\frac{1}{\mathcal E(\tau_1)}&=& 
\ds \sum_{k=1}^{n-1} \ds \frac{1}{\mathcal{E}(\tau_{k+1})}- \ds \frac{1}{\mathcal{E}(\tau_{k})}
\\[.3cm]
&\geq&
\ds\frac{1}{c}\ds \sum_{k=1}^{n-1} \ds \int_{\tau_k}^{\tau_{k+1}} \left( K_2 \sin^2(2\theta(\tau)) 
-\ds\frac{K_0}{c\rho(c)} 
\right)\ud \tau
\\[.3cm]
&\geq&
\ds\frac{K_2}{c}\ds\sum_{k=1}^{n-1} (\tau_{k+1}-\tau_k) 
 \ds\frac{1}{\tau_{k+1}-\tau_k} \left( \int_{\tau_k}^{\tau_{k+1}}  \sin^2(2\theta(\tau)) \frac{\dot{\theta}(\tau) \ud \tau}{\dot{\theta}(\tau)} \right)
\\[.1cm]
&&
-\ds \sum_{k=1}^{n-1} (\tau_{k+1}-\tau_k) \ds\frac{K_0}{c^2\rho(c)}
\\[.5cm]
&\geq&
\ds\frac{2K_2}{\pi c} \left( \ds \int_{0}^{\frac{\pi}{2}} \sin^2(2\theta) \ud \theta \right) \ds \sum_{k=1}^{n-1} (\tau_{k+1}-\tau_k) - \ds\frac{K_0(\tau_{n}-\tau_1)}{c^2\rho(c)} 
\\[.5cm]
&\geq&
(\tau_{n}-\tau_1)
\underbrace{\ds \frac{K_2}{2c} \left(1 - \ds\frac{2}{c\rho(c)}
\ds\frac{K_0}{K_2}
\right)}_{:=\Omega(c)}
\end{array}\]
with $\mathcal E(\tau_1)\leq 4\mathcal E_0$. For $c$ large enough, the coefficient $\Omega(c)$ in the right hand side  is strictly positive.
It yields
\begin{equation}\label{disc_est}\mathcal E(\tau_{n})\leqslant 
\frac{4\mathcal{E}_0}{1+4\mathcal E_0\Omega(c)(\tau_{n}-\tau_1)}.
\end{equation}

We now consider any time $t>\frac{2R_2}{c}$.
It can be associated to successive zeroes $\tau_n<\tau_{n+1}$ of $\theta$,   
such that $\tau_n\leq t<\tau_{n+1}$.
Owing to \eqref{estE2}, we have
\[
\ds\frac{1}{\mathcal E(t)}\geq \ds\frac{1}{\mathcal E(\tau_n)} - \ds\frac{K_0(t-\tau_n)}{c^2\rho(c)}.\]
Moreover, by \eqref{ecart}, we know that 
$0\leq t-\tau_n\leq \tau_{n+1}-\tau_{n}\leq \frac\pi b$.
Hence, coming back to \eqref{disc_est}, we obtain
\[\begin{array}{lll}
\ds\frac{1}{\mathcal E(t)}&\geq&
\ds\frac{1}{4\mathcal E_0}+\Omega(c)(\tau_{n}-\tau_1)
- \ds\frac{K_0\pi}{b c^2\rho(c)}
\\[.3cm]&\geq&
\ds\frac{1}{4\mathcal E_0}+
\ds\frac{\Omega(c)}{2}t
+\ds\frac{\Omega(c)}{2}\tau_n
-\ds\frac{\Omega(c)}{2}(t-\tau_n)
-\Omega(c)\tau_1
 -\ds\frac{K_0\pi}{b c^2\rho(c)}
 \\[.3cm]&\geq&
\ds\frac{1}{4\mathcal E_0}+
\ds\frac{\Omega(c)}{2}t
+\Omega(c)\left(
\ds\frac{\tau_n}{2}
-\ds\frac{\pi }{2b}
-\tau_1
 -\ds\frac{K_0\pi}{b c^2\rho(c)\Omega(c)}\right).
\end{array}\]
We can find $T_c$ large enough, so that when $t\geq T_c$, $\tau_n$ is sufficiently large to ensure that 
the last term between parenthesis is positive.
Accordingly, we get
\[
\mathcal E(t)\leq \ds\frac{4\mathcal E_0}{1+2\mathcal E_0\Omega(c)t}\]
for any $t\geq T_c$.
\QED
 
 \appendix
 
 \section{Orbital stability for the case with a confining potential}
\label{App}

The estimates obtained in Section~\ref{En_Min} cannot permit us to establish an orbital stability statement.
Such a statement relies on 
coercivity estimates which hold here ``up to space-translation'', hence providing a  quite poor information
since the localisation of the particles
is not characterized at all.
This difficulty disappears when considering the problem with a confining potential.
Namely, \eqref{N_particles} becomes
\begin{equation}\label{N_particles_pot}
\ddot q_{i}(t)=-
\ds\iint_{\mathbb R^d\times \mathbb R^n} 
\sigma_1(q_i(t)-y)\sigma_2(z) \nabla_y \psi(t,y,z)\ud y\ud z-\nabla_x V(q_i(t)),
\end{equation}
still with \eqref{N_wave} for defining the self-consistent potential.
Here the external potential $V: x\in \mathbb R^N\to [0,\infty)$ is required to satisfy
\begin{equation}
\label{hyp_potext}
\left\{
\begin{array}{l}
V\in C^2,\qquad V(0)=0,\qquad \nabla V(0)=0,\\
 \textrm{there exists $C_*>0$ such that for any $D^2V(x)\xi\cdot \xi\geq C_*|\xi|^2$ for any $x,\xi \in \mathbb R^n$}. 
\end{array}
\right.
\end{equation} 
Accordingly, $V$ admits a  unique minimum located at the origin.

Searching for equilibrium solutions, we still obtain
\[
p_j^*=0,\quad \varpi^*=0,\quad
-\Delta_z \psi^*(x,z)=-\sigma_2(z)\left(\sum\limits_{k=1}^{P}\sigma_{1}(x-q^*_{k})\right)\]
but \eqref{eq_defequilibrium} becomes
\[-\kappa\ds\sum_{k=1}^P \nabla\Sigma_1(q_j^*-q^*_k)
-\nabla V(q_j^*)=0.
\]
Clearly, having $q^*_1=...=q^*_P=0$ satisfies this relation. We claim that there is no other solution
of this equilibrium equation.
Note also that there is no isolated particle far from 0: if there exists $j\in \{1,..,P\}$ such that 
for any $k\in \{1,..,P\}$, we have $| q^*_j- q^*_k|\geq 2R_1$, then $\nabla\Sigma_1( q^*_j- q^*_k)=0$ and 
it follows that $\nabla V( q_j^*)=0$ hence $ q^*_j=0$.
As before, we adopt the notation \eqref{defdjk}
where we now deal with pseudo-particles $\tilde q^*_1,...,\tilde q^*_{P'}$
with $\tilde q^*_j\neq \tilde q^*_k$ for any $j,k\in \{1,...,P'\}$, $1<P'\leq P$, and we now have
\begin{equation}\label{equil:pot2}
\kappa\ds\sum_{k=1}^{P'} d_{j,k} (\tilde q_j^*-\tilde q^*_k)
+\nabla V(\tilde q_j^*)=0,\end{equation}
with $d_{j,k}$ defined by \eqref{defdjk}.
We remind the reader that $r_j$ stands for the number of particles represented by $\tilde q^*_j$.
We multiply \eqref{equil:pot2} by $r_j \tilde q^*_j$:
\[\begin{array}{l}
\ds\sum_{j,k=1}^{P'} r_j d_{j,k} (\tilde q^*_j-\tilde q^*_k)\cdot \tilde q^*_j 
=\ds\sum_{k,j=1}^{P'} r_k d_{k,j} (\tilde q^*_k-\tilde q^*_j)\cdot \tilde q^*_k 
\\
\qquad =
\ds\frac12 \ds\sum_{j,k=1}^{P'} r_j d_{j,k} |\tilde q^*_j|^2 
+\ds\frac12 \ds\sum_{j,k=1}^{P'} r_k d_{k,j} |\tilde q^*_k|^2 
-2 \ds\sum_{j,k=1}^{P'} \ds\frac{r_j d_{j,k}+r_k d_{k,j}}{4}
 \tilde q^*_j\cdot \tilde q^*_k
 \\
 \qquad =-
 \ds\sum_{k=1}^{P'} r_j \nabla V(\tilde q^*_j)\cdot \tilde q^*_j\leq -C_* \ds\sum_{k=1}^{P'} r_j |\tilde q^*_j|^2\leq 0,
\end{array}\]
by virtue of \eqref{hyp_potext}.
We set $M_{j,k}=\frac{r_j d_{j,k}+r_k d_{k,j}}{4}=M_{k,j}\geq 0$.
We arrive at 
\begin{align*}
\ds\sum_{j,k=1}^{P'} M_{j,k} |\tilde q^*_j-\tilde q^*_k|^2
+\ds\frac14  &\ds\sum_{j,k=1}^{P'} \big(r_j d_{j,k} - r_k d_{k,j}\big) |\tilde q^*_j|^2 
+\ds\frac14  \ds\sum_{j,k=1}^{P'} \big(  r_k d_{k,j}-r_j d_{j,k}\big) |\tilde q^*_k|^2 \\
&\leq -C_* \ds\sum_{k=1}^{P'} r_j |\tilde q^*_j|^2\leq 0.
\end{align*}
However, by definition, we have
\[\ds\sum_{k=1}^{P'} r_j d_{j,k} 
=\ds\sum_{k=1}^{P'}  
 r_k d_{k,j}
\]
and the last two sums in the left hand side vanish.
It follows that
\[0\leq \ds\sum_{j,k=1}^{P'} M_{j,k} |\tilde q^*_j-\tilde q^*_k|^2\leq -C_* \ds\sum_{k=1}^{P'} r_j |\tilde q^*_j|^2\leq 0.\]
We conclude that all particles are packed at the origin.

The energy functional now takes into account the external potential: 
\[\begin{array}{lll}
\mathscr H(X)&=&\ds\sum_{j=1}^P \left(\ds\frac{|p_j|^2}{2}+ V(q_j)\right)
+\ds\frac{1}{2}\ds\iint_{\mathbb R^d\times \mathbb R^n} |\varpi|^2\ud z\ud x+ \ds\frac12\ds\iint_{\mathbb R^d\times\mathbb R^n} |\nabla_z \psi|^2\ud z\ud x
\\&&+\ds\sum_{j=1}^P\ds\iint_{\mathbb R^d\times\mathbb R^n} 
\sigma_1(x-q_j)\sigma_2(z)\psi(x,z)
\ud z\ud x
\end{array}\]
with $X=(q_1,...,q_P,p_1,...,p_P,\psi,\varpi)$.
The equilibrium is given by $X^*=(0,0,\psi^*(x,z)=-P\Gamma(z ) \sigma_1(x),0)$.
It minimizes the energy with
\[
\mathscr H(X^*)=-\ds\frac{P^2}{2}\kappa \|\sigma_1\|^2_{L^2(\mathbb R^n)},\]
and we readily check that 
\[\begin{array}{lll}
\nabla\mathscr H(X^*)(X)
&=&\ds\iint_{\mathbb R^d\times\mathbb R^n}\nabla_z\psi^*\cdot\nabla_z\psi\ud z\ud x
+
P \ds\iint_{\mathbb R^d\times\mathbb R^n}\sigma_1(x)\sigma_2(z)\psi(x,z)\ud z \ud x
\\&&
-
\ds\sum_{j=1}^P \ds\iint_{\mathbb R^d\times\mathbb R^n}\nabla_x\sigma_1(x)\cdot q_j \sigma_2(z)\psi^*(x,z)\ud z \ud x
\\
&=&
P\ds\iint_{\mathbb R^d\times\mathbb R^n}(\Delta_z\Gamma +\sigma_2(z))\sigma_1(x)\psi(x,z)\ud z\ud x
\\&&+P\ds\sum_{j=1}^P q_j \cdot 
\ds\iint_{\mathbb R^d\times\mathbb R^n}\sigma_1 (x) \nabla_x\sigma_1(x) \sigma_2(z)\Gamma(z)\ud z \ud x=0
.\end{array}\]
Next, with $X=(q,p,\psi,\varpi)$, we compute
\begin{align*}
&\ds\frac12 D^2\mathscr H(X^*)(X,X)
\\
&=\ds\frac12 \ds\sum_{j=1}^P\left(|p_j|^2+D^2V(0)q_j\cdot q_j\right)
+
\ds\frac12 \ds\iint_{\mathbb R^n\times \mathbb R^d} |\nabla_z\psi|^2\ud z\ud x+
\ds\frac1{2}\ds\iint_{\mathbb R^n\times \mathbb R^d} |\varpi|^2\ud z\ud x
\\
&\quad-\ds\sum_{k=1}^P 
\ds\iint_{\mathbb R^d\times\mathbb R^n}\sigma_2(z) \nabla_x\sigma_1(x)\cdot q_k \psi\ud z\ud x
-\ds\frac12\ds\sum_{k=1}^P 
\ds\iint_{\mathbb R^d\times\mathbb R^n}\sigma_2(z) \nabla_x\sigma_1(x) \cdot q_k  \
\nabla_x\psi^*\cdot q_k\ud z\ud x.
\end{align*}
The last term is nothing but
\[
\ds\frac P2\ds\sum_{k=1}^P 
\ds\iint_{\mathbb R^d\times\mathbb R^n} \Gamma(z)\sigma_2(z) \nabla_x\sigma_1(x) \cdot q_k  \
\nabla_x\sigma_1(x)\cdot q_k\ud z\ud x
=\ds\frac{\alpha \kappa P}2\ds\sum_{k=1}^P |q_k|^2.
\]
It is convenient to rewrite the penultimate term by means of the Fourier transform
\[
-\ds\sum_{k=1}^P 
\ds\iint_{\mathbb R^d\times\mathbb R^n}\ds\frac{\xi\widehat \sigma_2(\xi)}{|\xi|^2}
\cdot \xi\widehat {\psi}(x,\xi)
 \nabla_x\sigma_1(x)\cdot q_k \ds\frac{\ud \xi}{(2\pi)^n}\ud x,\]
 and to combine it with the integral of $|\nabla_z\psi|^2$.
Hence, 
we arrive at 
\begin{equation}\label{disspD2_1}\begin{array}{lll}
\ds\frac12 
D^2\mathscr H(X^*)(X,X)
&=&\ds\frac12 \ds\sum_{j=1}^P\left(|p_j|^2+D^2V(0)q_j\cdot q_j\right)
+
\ds\frac1{2}\ds\iint_{\mathbb R^d\times\mathbb R^n} |\varpi|^2\ud z\ud x
\\
&&
+
\ds\frac12\ds\iint_{\mathbb R^d\times\mathbb R^n}
\Big|
\xi \widehat \psi(x,\xi) - \ds\frac{\xi\widehat \sigma_2(\xi)}{|\xi|^2}\ds\sum_{k=1}^P \nabla_x\sigma_1(x)\cdot q_k
\Big|^2 \ds\frac{\ud \xi}{(2\pi)^n}\ud x
\\&&
+\ds\frac{\alpha \kappa P} 2 \ds\sum_{k=1}^P |q_k|^2
-
\ds\frac12
\ds\iint_{\mathbb R^d\times\mathbb R^n}\ds\frac{|\widehat \sigma_2(\xi)|^2}{|\xi|^2}
\Big|\ds\sum_{k=1}^P
 \nabla_x\sigma_1(x)\cdot q_k\Big|^2 \ds\frac{\ud \xi}{(2\pi)^n}\ud x
 .
\end{array}\end{equation}
The last two terms in \eqref{disspD2_1}  combine as
\begin{equation}\label{det_dissipD2}
\begin{array}{l}
\ds\frac{\kappa}{2}
\left( \alpha  P\ds\sum_{k=1}^P |q_k|^2
- \ds\int_{\mathbb R^d} \Big|\ds\sum_{k=1}^P \nabla_x\sigma_1(x)\cdot q_k\Big|^2\ud x
\right) 
\\
\qquad
\geq 
\ds\frac{\kappa P}{2}
\left( \alpha  \ds\sum_{k=1}^P |q_k|^2
- \ds\sum_{k=1}^P \ds\int_{\mathbb R^d} \big| \nabla_x\sigma_1(x)\cdot q_k\big|^2\ud x
\right) 
\\
\qquad
\geq 
\ds\frac{\kappa P}{2}
\left( \alpha  \ds\sum_{k=1}^P |q_k|^2
- \ds\sum_{k=1}^P\left( \ds\int_{\mathbb R^d}  \nabla_x\sigma_1(x) \otimes \nabla \sigma_1(x) \ud x\right)q_k\cdot q_k\right)
=0,
\end{array}
\end{equation}
by using the Cauchy-Schwarz inequality.
Coming back to \eqref{hyp_potext}, we end up with the coercivity estimate
\[\begin{array}{lll}
D^2\mathscr H(X^*)(X,X)
&=&\ds\frac12\ds\sum_{j=1}^P\left(|p_j|^2+D^2V(0)q_j\cdot q_j\right)
+
\ds\frac1{2}\ds\iint_{\mathbb R^d\times\mathbb R^n} |\varpi|^2\ud z\ud x
\\&&+
\ds\frac12\ds\iint_{\mathbb R^d\times\mathbb R^n}
\Big|
\xi \widehat \psi(x,\xi) - \ds\frac{\xi\widehat \sigma_2(\xi)}{|\xi|^2}\ds\sum_{k=1}^P \nabla_x\sigma_1(x)\cdot q_k
\Big|^2 \ds\frac{\ud \xi}{(2\pi)^n}\ud x
\\
&\geq &\ds\frac12\ds\sum_{j=1}^P\left(|p_j|^2+C_*|q_j|^2\right)
+
\ds\frac1{2}\ds\iint_{\mathbb R^d\times\mathbb R^n} |\varpi|^2\ud z\ud x
\\&&+
\ds\frac12\ds\iint_{\mathbb R^d\times\mathbb R^n}
\left|
\nabla_z\left( \psi(x,z) -  \nabla_x\psi^*(x,z)\cdot\ds\frac1P\ds\sum_{k=1}^P q_k
\right)\right|^2 \ud z\ud x
\geq 0.
\end{array}\]
Let us introduce the shorthand notation $N^*(X)=\sqrt{D^2\mathscr H(X^*)(X,X)}$
and the natural energy space
$$\mathbb H=(\mathbb R^d)^{2P}\times L^2(\mathbb R^d\times\mathbb R^n)
\times L^2(\mathbb R^d;\dot H^1(\mathbb R^n)).$$

\begin{lemma}
The application $X\mapsto N^*(X)\geq 0$ defines a norm on the space $\mathbb H$ and there exists $\underline k,\overline k>0$ such that 
$\underline k\| X\|_{\mathbb H} \leq N^*(X)\leq  \overline k
\| X\|_{\mathbb H}$
holds for any $X\in \mathbb H$, where 
$$\|
 X\|_{\mathbb H}^2=
 \ds\sum_{j=1}^P( |p_j|^2 +|q_j|^2)+\ds\iint_{\mathbb R^n\times\mathbb R^d} ( |\varpi|^2+|\nabla_z \psi|^2) \ud z\ud x.$$
\end{lemma}

\noindent
{\bf Proof.}
 The application $X\mapsto N^*(X)$ satisfies
 $N^*(X)\geq 0$, $N^*(sX)=sN^*(X)$ for any 
$X\in \mathbb H$ and $s\geq 0$. We have $N^*(0)=0$ 
and, conversely, having $N^*(X)=0$ implies that 
$p_1=...=p_P=0$, $q_1=...=q_P=0$, $\varpi=0$ and 
$$\ds\iint_{\mathbb R^d\times\mathbb R^n}
\left|
\nabla_z\left( \psi(x,z) -  \nabla_x\psi^*(x,z)\cdot\ds\frac1P\ds\sum_{k=1}^P q_k
\right)\right|^2 \ud z\ud x=0.
$$
Since the $q_k$'s vanish it reduces to $\|\nabla_z\psi\|_{L^2(\mathbb R^n\times \mathbb R^d)}^2=0$, hence $\psi=0$ too, and finally $X=0$.
The triangle inequality is a direct consequence of the Cauchy-Schwarz inequality.

For proving the equivalence of norms, only the term 
$$\ds\iint_{\mathbb R^d\times\mathbb R^n}
\left|
\nabla_z\left( \psi(x,z) -  \nabla_x\psi^*(x,z)\cdot\ds\frac1P\ds\sum_{k=1}^P q_k
\right)\right|^2 \ud z\ud x
$$
deserves to be discussed.
 It is dominated by 
 \[\begin{array}{l}
2 \ds\iint_{\mathbb R^d\times\mathbb R^n}
|\nabla_z\psi(x,z)|^2 \ud z\ud x 
+
2
 \ds\iint_{\mathbb R^d\times\mathbb R^n}
 |\nabla_z\Gamma(z)|^2\
 \left| 
 \nabla_x\sigma_1(x)\cdot\ds\sum_{k=1}^P q_k\right|^2
 \ud z\ud x
\\
\qquad\leq
2\|\nabla_z\psi\|_{L^2(\mathbb R^n\times \mathbb R^d)}^2
+2\kappa P  \ds\iint_{\mathbb R^d}|\nabla_x\sigma_1|^2\ud x \ds\sum_{k=1}^P |q_k|^2
 ,\end{array}\]
 which permits us to identify $\overline k$.
 
 Next, we argue by contradiction, assuming  that we can find a sequence $\big(X_\nu\big)_{\nu\in \mathbb N}$
 such that $\|X_\nu\|_{\mathbb H}=1$ but $\lim_{\nu\to \infty} N^*(X_\nu)=0$.
 The latter implies that $p_{j,\nu}$, $q_{j,\nu}$ and $\varpi_\nu$ tend to 0 as $\nu\to \infty$
 and thus we deduce that $\lim_{\nu\to\infty}\|\nabla_z\psi\|_{L^2}=1$.
 It also implies that 
 $$\ds\lim_{\nu \to \infty}
\ds\iint_{\mathbb R^d\times\mathbb R^n}
\left|
  \nabla_z\cdot\left(\nabla_x\psi^*(x,z)\cdot\ds\sum_{k=1}^P q_{k,\nu}\right)
\right|^2 \ud z\ud x=0,
 $$
 by a direct application of Lebesgue's theorem.
 Then, we equally have
 \[
 \begin{array}{lll}
 \| \nabla_z\psi_\nu\|_{L^2}&\leq& 
\left(\ds\iint_{\mathbb R^d\times\mathbb R^n}
\left|\nabla_z\left( \psi_\nu(x,z)-
  \nabla_x\psi^*(x,z)\cdot\ds\sum_{k=1}^P q_{k,\nu}
\right)\right|^2 \ud z\ud x\right)^{1/2}
\\&&\qquad+
\left(\ds\iint_{\mathbb R^d\times\mathbb R^n}
\left|
 \nabla_z\cdot\left( \nabla_x\psi^*(x,z)\cdot\ds\sum_{k=1}^P q_{k,\nu}\right)
\right|^2 \ud z\ud x\right)^{1/2}
 \xrightarrow[\nu\to \infty]{}0,\end{array}\]
which leads to a contradiction.
\QED

We now prove the stability of the equilibrium state.
We know that  $$\underline k \|\xi\|^2_{\mathbb H}\leq D^2\mathscr H(X^*)(\xi,\xi)=|N^*(\xi)|^2\leq \overline k \|\xi\|^2_{\mathbb H}$$ holds for any $\xi\in \mathbb H$.
By continuity of $X\mapsto D^2\mathscr H(X)$, there exists $\delta >0$ such that 
if $\| X-X^*\|_{\mathbb H}\leq \delta$ then, for any $\xi\in \mathbb H$, we still have
\[ \frac{\underline k } 2\|\xi\|_{\mathbb H}^2\leq  \frac1 2 |N^*(\xi)|^2 \leq D^2\mathscr H(X)(\xi,\xi)
\leq 2 |N^*(\xi)|^2\leq 2\overline k \|\xi\|^2_{\mathbb H}
 .\] 
Accordingly,
\[\mathscr H(X)-\mathscr H(X^*)=\ds\int_0^1(1-\theta)D^2\mathscr H(X^*+\theta(X-X^*))(X-X^*,X-X^*)\ud \theta
\]
satisfies 
\[ \frac{\underline k } 2\|X-X^*\|_{\mathbb H}^2\leq 
\mathscr H(X)-\mathscr H(X^*)\leq 2\overline k \|X-X^*\|^2_{\mathbb H}\]
 as far as $\| X-X^*\|\leq \delta$.

Let $t\mapsto X(t)$ stands for the solution of the particle-wave system with confining potential that starts from the initial date $X_0$.
We are going to use energy conservation
\[
\mathscr H(X(t))-\mathscr H(X^*)=\mathscr H(X_0)-\mathscr H(X^*).\]
Let $0<\epsilon<\delta$. We suppose $\|X_0-X^*\|< \epsilon$. Set 
\[T=\sup\{t\geq 0, \|X(s)-X^*\|_{\mathbb H}\leq \epsilon,\ 0\leq s\leq t\}.\]
By continuity $T>0$. If $T$ is finite, then $\|X(T)-X^*\|_{\mathbb H}=\epsilon$ and we can write
\[
\frac{\underline k } 2\|X(T)-X^*\|_{\mathbb H}^2\leq
\mathscr H(X(T))-\mathscr H(X^*)\leq 2\overline k \|X_0-X^*\|^2_{\mathbb H}.\]
Imposing
$
\|X_0-X^*\|^2_{\mathbb H}\leq \frac{\underline k}{8\overline k}\epsilon^2 
$
yields a contradiction.
We conclude with the following stability statement.

\begin{proposition}
Let $V$ satisfy \eqref{hyp_potext}.
For any $\epsilon\in (0,\delta)$, there exists $\eta>0$ such that if 
the initial data for \eqref{N_particles_pot} coupled to \eqref{N_wave} satisfies
$\|X_0-X^*\|_{\mathbb H}\leq \eta$, then 
$\|X(t)-X^*\|_{\mathbb H}\leq \epsilon$ holds for any $t\geq 0$.
\end{proposition}

  \section*{Acknowledgements}
 This work is supported by the ANR project MaDynOS, ANR-24-CE40-3535.
 It also takes advantage of the support of the CNRS international Research Project KinEq., which has permitted fruitful discussions with P.-E.~Jabin and A.~Mellet.


\end{document}